\documentclass[a4paper,10pt]{article}
\usepackage[paper=a4paper, left=3cm, right=3cm, top=2.5cm, bottom=3cm]{geometry} 
\usepackage{graphicx, amssymb, amsmath, epsfig}
\usepackage{color}
\usepackage{booktabs}
\usepackage{boxedminipage}
\usepackage[ruled]{algorithm2e}
\usepackage{algorithmic}
\usepackage[utf8]{inputenc}
\usepackage{epstopdf}       
\usepackage{float}
\usepackage{upgreek}

\graphicspath{{./}{figs/}}

\title{POD reduced order modeling for evolution equations utilizing arbitrary finite element discretizations}
\author{Carmen Gr\"a\ss{}le and Michael Hinze} 
\date{}
\begin{document}

\maketitle

\begin{abstract}
The main focus of the present work is the inclusion of spatial adaptivity for the snapshot computation in the offline phase of model order reduction utilizing Proper Orthogonal Decomposition (POD-MOR) for nonlinear parabolic evolution problems. We consider snapshots which live in different finite element spaces, which means in a fully discrete setting that the snapshots are vectors of different length. From a numerical point of view, this leads to the problem that the usual POD procedure which utilizes a singular value decomposition of the snapshot matrix, cannot be carried out. In order to overcome this problem, we here construct the POD model / basis using the eigensystem of the correlation matrix (snapshot gramian), which is motivated from a continuous perspective and is set up explicitly e.g. without the necessity of interpolating snapshots  into a common finite element space. It is an advantage of this approach that the assembling of the matrix only requires the evaluation of inner products of snapshots in a common Hilbert space. This allows a great flexibility concerning the spatial discretization of the snapshots. The analysis for the error between the resulting POD solution and the true solution reveals that the accuracy of the reduced order solution can be estimated by the spatial and temporal discretization error as well as the 
POD error. Finally, to illustrate the feasibility our approach, we present a test case of the Cahn-Hilliard system utilizing $h$-adapted hierarchical meshes and two settings of a linear heat equation using nested and non-nested grids.\\

\noindent {\scriptsize \textit{Acknowlegdements:} We like to thank Christian Kahle for providing many C++ libraries which we could use for the coding. The authors gratefully acknowledge the financial support by the Deutsche Forschungsgemeinschaft through the priority program SPP1962 entitled ``Non-smooth and Complementarity-based Distributed Parameter Systems: Simulation and Hierarchical Optimization''.} \hspace{2cm} 

\end{abstract}

\section{Introduction}
\noindent For many problem settings, model order reduction utilizing Proper 
Orthogonal Decomposition (POD) has proven to be a powerful tool in order to reduce large scale systems to surrogate models of low dimension while preserving a good approximation quality. The range of applications of POD model order reduction (POD-MOR) comprises a broad scope, including (amongst others) linear and nonlinear parabolic equations \cite{KV01}, optimal control 
of partial differential equations \cite{GV10, Vol01, KV99, HV05, AH01} and fluid dynamics \cite{Lum67, HK00, LMQR14}. A general introduction to POD and 
reduced order modeling can be found in \cite{HLBR12, Pin09, Vol13}, for example. The key idea of the POD technique is to apply a Galerkin ansatz, in which the ansatz functions, i.e. the POD basis functions, contain information about the underlying dynamical system. Following the approach of snapshot based POD in \cite{Sir87}, the system information is retrieved from snapshots of the solution trajectory at several time instances, which are generated in a simulation.\\ 
In practice, many simulations require adaptive strategies for the spatial discretization in order to be implementable. For example, in the simulation of Cahn-Hilliard systems based on diffuse interface approaches, many degrees of freedom are required at the interfacial regions in order to well resemble the steep gradients, whereas in the pure phases only little numbers of degrees of freedom are needed, see \cite{HHT11} for example. Utilizing a uniform mesh for such problems would drastically enlarge the computational effort and storage problems would occur.\\
In \cite{BR14} adaptive finite element methods and POD reduced order modeling are considered as two different techniques in order to reduce the complexity of the numerical solution of optimal control problems. The combination of both approaches contains a major challenge. The inclusion of spatial adaptivity in the context of model reduction means, in a discrete formulation, that the snapshots at each time instance may have different lengths due to their different spatial resolutions. For this reason, the snapshot matrix cannot be set up directly and the usual POD procedure cannot be carried out. In this paper, we consider the combination of adaptive finite element snapshots with POD model order reduction. Our perspective is based on the continuous setting, which allows us to derive a reduced order model utilizing snapshots which only need to lie in a common Hilbert space. Thus, for the actual numerical implementation, this allows us to use spatially adapted finite element snapshots or even a blend of snapshots stemming from different discretization schemes. \\
The inclusion of spatial adaptivity in the POD framework is advantageous from two perspectives: on the one hand, the use of adaptive finite elements for snapshot generation can remarkably reduce the offline computation time in comparison to the use of a uniform mesh with resolution of the finest level of the adaptive grids. On the other hand, we expect to speed up the computations when solving the POD surrogate model in contrast to utilizing adaptive finite elements, since we solve POD reduced systems of low order.\\
In order to overcome the difficulties arising from combining POD with adaptive finite elements, different concepts are proposed. In \cite{FPNPGAG09}, the use of dynamically adaptive meshes is combined with POD-MOR for an unstructured ocean model. A fixed reference mesh is utilized, onto which the spatial adaptive snapshots are interpolated. This allows snapshots of the same lengths 
at each time level and the usual POD procedure can then be performed on this fixed reference mesh. In \cite[Ch. 2.4.3]{Las14} an interpolation approach is outlined. The idea is to interpolate given snapshots of arbitrary spatial discretizations by (piecewise) polynomials. For the fully discrete POD setting, the spatial discretization points are chosen appropriately for the numerical integration of the polynomials. In the context of reduced basis methods, adaptive wavelet discretizations are used in \cite{ASU15} in the offline snapshot computation phase. In \cite{Yan16}, a reduced basis method is developed which is certified by a residual bound relative to the infinite-dimensional weak solution. Different adaptive strategies are considered 
for both the finite element and the reduced basis level. Furthermore, in \cite{HKP14} three numerical concepts to treat the moving free boundary 
for the calculation of the snapshots  are compared: first a Landau-type transformation, second a control volume approach and third a moving mesh approach. In contrary to the first two concepts, the number of grid points in the moving mesh approach (r-adaptivity) is kept fixed, but they are moved according to the evolution of the free boundary. POD is applicable with only minor modifications. Recently and in parallel to our work, the combination of POD model reduction with adaptive finite element snapshots is realized in \cite{URL16} by constructing common finite element spaces. Two options are considered: either all snapshots are expressed in terms of a common finite element basis or pairs of snapshots are expressed in terms of a common finite element basis of these pairs. Moreover, error statements for a parametrized elliptic boundary value problem are proved. In the numerical examples, $h$-adaptive finite elements with fixed polynomial degree are utilized.\\
The aim of this work is to derive a POD reduced order model for a general semilinear evolution equation which can be set up and solved for arbitrary finite element discretizations without the necessity of e.g. using interpolation with respect to the spatial variable. This approach is motivated from a continuous perspective, where snapshots from different finite element spaces belong to a same Hilbert space. The assembly of the snapshot matrix can be avoided by directly setting up an eigenvalue problem for which we only need the calculation of the inner product of the finite element ansatz functions. This can be computed for arbitrary finite element discretizations and suffices to set up a POD-ROM. We provide an algorithm which enables the evaluation of the inner products of snapshots stemming from arbitrary finite element discretizations. In the numerical examples, the method is demonstrated for both $h$-adapted snapshots and snapshots with non-nested meshes. Moreover, the treatment of the nonlinearity is discussed.\\
The paper is organized as follows: the general problem setting of an 
abstract semilinear parabolic evolution problem is described in Section 2. We recall the POD method in real Hilbert spaces in Section 3 and set up a POD eigenvalue problem which can be assembled for any finite element discretization. The POD reduced order model for arbitrary finite element discretizations is discussed in Section 4. In order to validate the quality 
of the POD solution, we investigate the error between the POD solution and the true solution in Section 5. Numerical examples are presented in Section 6 to illustrate our approach.\\

\section{Abstract semilinear parabolic evolution problem}

\subsection{Problem setting}

\noindent Let us specify the abstract semilinear parabolic evolution problem which we consider in the following. Let $(V, \langle \cdot , \cdot \rangle_V)$ and $(H, \langle \cdot , \cdot \rangle_H)$ be real separable Hilbert spaces such that there exists a dense and continuous embedding $V \hookrightarrow H$. The dual space $H'$ can be identified with $H$ by the Riesz representation theorem and the Gelfand triple $(V,H,V')$ is formed by $V \hookrightarrow H = H' \hookrightarrow V'$. Since $V$ is continuously embedded in $H$, there 
exists a constant $c_v > 0$ such that 
\begin{equation}\label{norm_estimation}
 \| v \|_H^2 \leq c_v \| v \|_V^2 \quad \text{ for all } v \in V.
\end{equation}
For a given symmetric, $V$-elliptic bilinear form $a: V \times V \to \mathbb{R}$, we assume boundedness, i.e. 
\begin{equation}\label{bilinear}
 \exists \beta > 0: \quad |a(u,v)| \leq \beta \| u \|_V \| v \|_V \quad \text{ for all } u,v \in V
\end{equation}
and coercivity, i.e. 
$$\exists \kappa > 0: \quad a(u,u) \geq \kappa \| u \|_V^2 \quad \text{ for all } u \in V.$$
Let $\mathcal{A}: V \to V'$ be the bounded linear operator associated with the bilinear form $a$, i.e. $\mathcal{A} \in \mathcal{L}(V,V')$ and 
$$ a(u,v) = \langle \mathcal{A} u, v \rangle_{V',V} = \langle u \mathcal{A}^\star v \rangle_{V,V'} \quad \text{for all } u, v \in V,$$
where $\langle \cdot, \cdot \rangle_{V',V}$ denotes the dual pairing of $V'$ and $V$. Moreover, we denote by $\mathcal{N}: V \to V'$ a nonlinear operator.
We are concerned with the following Cauchy problem for a semilinear evolution problem. Let $T >0$ be a fixed end time. For a given initial function $g \in H$ and external force $f \in L^2(0,T;V')$ we consider the problem: find $y \in W(0,T) := \{ v \in L^2(0,T; V), v_t \in L^2(0,T;V')\}$ with 
\begin{equation}\label{P}
  \left\{
 \begin{array}{rcl}
   \frac{d}{dt} \langle y(t), v \rangle_H + a(y(t),v) + \langle \mathcal{N}(y(t)), v \rangle_{V',V} & = & \langle f(t), v \rangle_{V',V},\\   
   \langle y(0), v \rangle_H & = & \langle g, v \rangle_H, \\   
 \end{array}
 \right.
 \end{equation}
 for all $v \in V$ and for almost all $t \in (0,T]$. Note that it holds $\frac{d}{dt} \langle y(t), v \rangle_H = \langle y'(t),v \rangle_{V',V}$ for 
 all $y \in L^2(0,T;V)$ with $y' \in L^2(0,T;V')$ and all $v \in V$ in the sense of distributions in $(0,T)$.\\
 
 \noindent \textbf{Assumption 2.1.} For every $f \in L^2(0,T; V')$ and $g \in H$ there exists a unique weak solution of \eqref{P} with 
 $$ y \in L^2(0,T;V) \cap C([0,T];H) \cap H^1(0,T;V').$$\\
 
 \noindent \textbf{Remark 2.2.} Under monotonicity, boundedness and Lipschitz 
 continuity assumptions on the nonlinear operator $\mathcal{N}$, existence and uniqueness results for a general abstract evolution equation of type \eqref{P} are proved in \cite[Th. 4.1]{Yag10} or \cite[Ch. 6]{Paz83}, for example.\\

\noindent \textbf{Example 2.3 (Semilinear heat equation).} Let $\Omega \subset 
 \mathbb{R}^k, k \in \{2,3\}$ be a bounded open domain with Lipschitz continuous boundary $\partial \Omega$ and let $ T > 0$ be a fixed end time. We set $Q := (0,T) \times \Omega$ and $\Sigma := (0,T) \times \partial \Omega$ and $c \geq 0$. For a given forcing term $f \in L^2(Q)$ and initial condition $g \in L^2(\Omega)$, we consider the semilinear heat equation with homogeneous Dirichlet boundary condition:

\begin{equation}\label{heat}
  \left\{
 \begin{array}{r c l l}
   y_t (t,x) - \Delta y (t,x) + c y^3(t,x)& = & f(t,x) & \text{in } Q,\\   
   y(t,x) & = & 0 & \text{on } \Sigma, \\  
   y(0,x) & = & g(x) & \text{in } \Omega.
 \end{array}
 \right.
 \end{equation}
 
 \noindent The existence of a unique solution to \eqref{heat} is proved in \cite{RZ99}, for example. We can write \eqref{heat} as an abstract evolution problem of type \eqref{P} by deriving a variational formulation for \eqref{heat} with $V=H_0^1(\Omega)$ as the space of test functions, $H=L^2(\Omega)$ and integrating over the space $\Omega$. The bilinear form
 $a: V \times V \to \mathbb{R}$ is introduced by
 $$a(u,v) = \int_\Omega \nabla u \cdot \nabla v$$
  and the operator $\mathcal{N}: V \to V'$ is defined as $\mathcal{N}(y) = y^3$. For $c \equiv 0$, the heat equation \eqref{heat} is linear.\\
 
  \noindent \textbf{Example 2.4 (Cahn-Hilliard equations).} Let $\Omega, T, Q$ and $ \Sigma$ be defined as in Example 2.3. The Cahn-Hilliard system was proposed in \cite{CH58} as a model for phase separation in binary alloys. Introducing the chemical potential $w$, the Cahn-Hilliard equations can be formulated in the common setting as a coupled system for the phase field $c$ and the chemical potential $w$:
 
\begin{equation}\label{CHcoupled} 
  \left\{
\begin{array}{rcll}
c_t(t,x) + y \cdot \nabla c(t,x) & = & m \Delta w(t,x) & \text{in } Q,\\
w(t,x) & = & -\sigma \varepsilon \Delta c(t,x) + \frac{\sigma}{\varepsilon} W'(c(t,x)) & \text{in } Q,\\
\nabla c (t,x) \cdot \nu_{\Omega} &  =  & \nabla w (t,x) \cdot \nu_{\Omega}  = 0  & \text{on } \Sigma,\\
c(0,x) & = & c_0(x) & \text{in } \Omega.
\end{array}
\right.
\end{equation}

\noindent By $\nu_{\Omega}$ we denote the outward normal on $\partial \Omega$, $m \geq 0$ is a constant mobility, $\sigma > 0$ denotes the surface tension and $0 < \varepsilon \ll 1$ represents the interface parameter. Note that the convective term $y \cdot \nabla c$ describes the transport with (constant) velocity $y$. The transport term represents the coupling to the Navier-Stokes 
equations in the context of multiphase flow, see e.g. \cite{HH77} and \cite{AGG12}. The phase field function $c$ describes the phase of a binary material with components $A$ and $B$. It is $c \equiv -1$ in the pure $A$-phase and $c \equiv +1$ in the pure $B$-phase. The interfacial region is described by $c \in (-1,1)$ and its thickness is finite and of order $\mathcal{O}(\varepsilon)$. The function $W(c)$ is a double well potential. A typical choice for $W$ is the polynomial free energy function 
\begin{equation}\label{poly}
 W^p(c) = (1-c^2)^2/4, 
\end{equation}
which has exactly two minimal points at $c = \pm 1$, i.e. at the energetically favorable state. Another choice for $W$ is the relaxed double obstacle free energy
\begin{equation}\label{rdofe}
  W^{\text{rel}}_s(c) = \frac{1}{2} (1-c^2) + \frac{s}{2}(\max(c-1,0)^2+\min(c+1,0)^2),
\end{equation}
with relaxation parameter $s \gg 0$, which is introduced in \cite{HHT11} as the Moreau-Yosida relaxation of the double obstacle free energy
\begin{equation*}
W^\infty(c) = \begin{cases}
 \frac{1}{2}(1-c^2),  & \text{if }c \in [-1,1],\\
  +\infty, & \text{else. }
\end{cases}
\end{equation*}
For more details on the choices for $W$ we refer to \cite{Abe07} and \cite{BE91}, for example. Concerning existence, uniqueness and regularity of a solution to \eqref{CHcoupled}, we refer to \cite{BE91}.

In order to derive a variational form of type \eqref{P}, we write \eqref{CHcoupled} as a single fourth-order parabolic equation for $c$ by

\begin{equation}\label{CHsingle} 
  \left\{
\begin{array}{rcll}
c_t(t,x) + y \cdot \nabla c  & = & m\Delta (-\sigma\varepsilon \Delta c(t,x) + \frac{\sigma}{\varepsilon} W'(c(t,x))) & \text{in } Q,\\
\nabla c (t,x) \cdot \nu_{\Omega} &  =  & \nabla (-\sigma\varepsilon \Delta c(t,x) + \frac{\sigma}{\varepsilon} W'(c(t,x))) \cdot \nu_{\Omega} = 0& \text{on } \Sigma,\\
c(0,x) & = & c_0(x) & \text{in } \Omega.
\end{array}
\right.
\end{equation}

\noindent We choose $V= \{ v \in H^1(\Omega): \frac{1}{|\Omega|} \int_\Omega v = 0 \}$ equipped with the inner product $(u,v)_V:= \int_\Omega \nabla u \nabla v$, so that the dual space of $V$ is given by $V'=\{f \in (H^1(\Omega))': \langle f,1 \rangle = 0\}$ such that $V \hookrightarrow H=V'$ and $\langle ., .\rangle$ denotes the duality pairing. We note that $(V,(.,.)_V)$ is a Hilbert space. We define the $V'-$inner product for $f,g\in V'$ as $(f,g)_{V'} := \int_\Omega \nabla (-\Delta)^{-1} f \cdot \nabla (-\Delta)^{-1} g$ where $(-\Delta)^{-1}$ denotes the inverse of the negative Laplacian with zero Neumann boundary data. Note that $(f,g)_{V'} = (f,(-\Delta)^{-1}g)_{L^2(\Omega)} = ((-\Delta)^{-1}f,g)_{L^2(\Omega)}$. We introduce the bilinear form $a: V \times V \to \mathbb{R}$ by 
 $$ a(u,v) = \sigma \varepsilon (\nabla u, \nabla v )_{L^2(\Omega)} + \frac{1}{m} (y \cdot \nabla u, v)_{V'}$$
 and define the nonlinear operator $\mathcal{N}$ by $\mathcal{N}(c) = \frac{\sigma}{\varepsilon} W'(c).$
 The evolution problem can be written in the form
 \begin{equation}\label{CH-weak-abstract}
  \frac{1}{m}(c_t(t),v)_{V'} + a(c(t),v) + \langle \mathcal{N}(c(t)),v \rangle = 0 \quad \forall v \in V \text{ and a.a. } t \in (0,T]. \footnote{We acknowledge a hint of Harald Garcke who pointed this form of the weak formulation of \eqref{CHsingle} to us.}
 \end{equation}
 We note that this fits our abstract setting formulated in \eqref{P} with the Gelfand triple $V\hookrightarrow H \equiv V' \hookrightarrow V'$.

 \subsection{Temporal and spatial discretization}

\noindent In order to solve \eqref{P} numerically, we apply the implicit 
Euler method for temporal discretization. Of course, other time integration 
schemes are possible. For a given $n \in \mathbb{N}$ let 
\begin{equation}\label{timegrid}
 0 = t_0 < \dotsc < t_n = T
\end{equation}
denote an arbitrary grid in the time interval $[0,T]$ with time step sizes 
$$\Delta t_j := t_j - t_{j-1}$$ 
\noindent for $j = 1, \dotsc, n$. We set $I_j = [t_{j-1}, t_j]$ for each time interval $j = 1, \dotsc, n$. The resulting time-discrete system consists in finding a sequence $\{ \bar{y}_j\}_{j=0}^n \subset V$ satisfying the following system of equations 

\begin{equation}\label{P_timediscrete}
  \left\{
 \begin{array}{rcll}
  \langle \displaystyle\frac{ \bar{y}_j-  \bar{y}_{j-1}}{\Delta t_j} , v \rangle_H + a( \bar{y}_j,v) + \langle \mathcal{N}( \bar{y}_j), v \rangle_{V',V} & = & \displaystyle\frac{1}{\Delta t_j} \int_{t_{j-1}}^{t_j} \langle f_j, v \rangle_{V',V}dt, & \text{ for } j = 1, \dotsc, n\\   
   \langle  \bar{y}_0, v \rangle_H & = & \langle g, v \rangle_H, &\\   
 \end{array}
 \right.
 \end{equation}
 
 \noindent for all $v \in V$ and $f_j$ denotes $f_j = f(t_j) \in V'$ for $j = 1, \cdots, n$.\\         
 
 For the spatial discretization we utilize adaptive finite elements. At each 
 time point $t_j$, $j = 0, \dotsc, n,$ we introduce a regular triangulation $\mathcal{T}_j$ of $\bar{\Omega}$ and define an $N_j$-dimensional conformal subspace $V_j$ of $V$ by 
 $$V_j := \text{span} \{v_1^j , \dotsc, v_{N_j}^j\} \subset V $$
 with nodal basis $\{ v_i^j\}_{i=1}^{N_j}$, i.e. $v_i^j (P_k^j) = \delta_{ik}$ for $i,k = 1, \dotsc, N_j$ with the nodes $\{P_k^j\}_{k=1}^{N_j}$ of the 
 underlying triangulation $\mathcal{T}_j$. Therefore, at each time level $j=0, \dotsc , n$, the utilized finite element spaces $\{ V_j\}_{j=0}^n$ can differ both in the underlying triangulation of the domain $\bar{\Omega}$ and in the polynomial degree. This means that the solutions can be computed utilizing $h-$, $p-$ and $r-$adaptivity, where $h-$adaptivity denotes local refinement and coarsening of the triangulation according to certain error indicators, $p-$adaptivity means increasing and decreasing the polynomial degree according to the smoothness of the solution and $r-$adaptivity, or moving mesh methods, relocates the mesh points to concentrate them in specific regions. We apply a Galerkin scheme for \eqref{P_timediscrete}. Thus we look for a sequence $\{y_j\}_{j=0}^n $ with $y_j \in V_j$ which fulfills
\begin{equation}\label{P_fullydiscrete}
  \left\{
 \begin{array}{rcll}
  \langle \displaystyle\frac{ y_j-\mathbb{I}^j  y_{j-1}}{\Delta t_j} , v \rangle_H + a(y_j,v) + \langle \mathcal{N}( y_j), v \rangle_{V',V} & = & \langle f_j, v \rangle_{V',V}, & \text{ for } j = 1, \dotsc, n\\   
   \langle y_0, v \rangle_H & = & \langle g, v \rangle_H, &\\   
 \end{array}
 \right.
 \end{equation}
 
 \noindent for all $v \in V_j$, where $\mathbb{I}^j: C(\bar{\Omega }) \to V_j$ denotes the Lagrange interpolation. Since $y_j \in V_j$ holds, we make the Galerkin ansatz
 \begin{equation}\label{yFE}
 y_j  = \displaystyle\sum_{i=1}^{N^j} \mathsf{y}_i^j v_i^j   \in 
 V_j \subset V 
 \end{equation}
 for $j = 0, \dotsc, n$ with appropriate mode coefficients $\{ \mathsf{y}_i^j\}_{i=1}^{N^j}$.\\

\section{POD method utilizing snapshots with arbitrary finite element discretizations}

\subsection{POD method in real Hilbert spaces}

\noindent The aim of this work is to propose a POD-ROM which uses the correlation matrix in order to construct the POD basis and POD surrogate model and avoids the necessity of e.g. interpolating the snapshots into a common finite element space. For this reason, the POD method is explained from an infinite-dimensional perspective in this section, where we use a finite number of snapshots which lie in a common Hilbert space. The POD method in Hilbert spaces is explained in \cite{KV02} and \cite{Vol01}, for example. Here, we recall the main aspects.\\

\noindent Assume we are given snapshots 
$$ y_0 \in V_0, \dotsc, y_n \in V_n$$ 
of \eqref{P}, which can be finite element samples of the solution trajectory 
$\mathcal{V} = \text{span}\{ y(t) \; | \; t \in [0,T] \}$ 
for \eqref{P} on the given timegrid $\{t_j\}_{j=0}^n$ introduced 
in \eqref{timegrid}. For each time level $j = 0, \dotsc , n$ the
snapshots belong to different subspaces $V_0, \dotsc , V_n \subset V$. Note that by construction we have $ \mathcal{V} := \text{span}\{ y_j\}_{j=0}^n \subset V$.\\

\noindent The idea of the POD method is to describe the space $\mathcal{V}$ by 
means of few orthonormal functions $\{ \psi_i \}_{i=1}^\ell \subset V$, 
with $\ell \leq d := \dim \mathcal{V}$, such that error between the snapshots $\{ y_j \}_{j=0}^n$ and the projection of the snapshots onto the subspace $V^\ell = \text{span} \{\psi_1, ..., \psi_\ell\} \subset V$ is minimized in the following sense: 
\begin{equation}\label{minPOD}
\underset{\psi_1 , ... , \psi_\ell \in V}{\text{min }}
 \displaystyle\sum_{j=0}^n \alpha_j \left\| y_j - \displaystyle\sum_{i=1}^\ell 
 \langle y_j, \psi_i \rangle_X \; \psi_i \right\|_X^2  \hspace*{0.2cm} \text{ s.t. }  \langle \psi_i , \psi_j \rangle_X = \delta_{ij} \quad \text{ for } 1 \leq i, j \leq \ell
 \end{equation}
 
\noindent where $X$ denotes either $V$ or $H$ and with e.g. nonnegative 
trapezoidal weights $\{ \alpha_j \}_{j=0}^n$,
\begin{equation}\label{weights_alpha}
 \alpha_0 = \frac{\Delta t_1}{2}, \; \alpha_j = \frac{\Delta t_j + \Delta t_{j+1}}{2} \text{ for } j = 1, \dotsc, n-1 \; \text{ and } \alpha_{n} = \frac{\Delta t_{n}}{2}.             
\end{equation}
A solution to \eqref{minPOD} is called a rank-$\ell$ POD basis. For this 
equality constrained minimization problem \eqref{minPOD}, first-order necessary optimality conditions can be derived. For this purpose, we introduce the bounded linear operator $\mathcal{Y}: \mathbb{R}^{n+1} \to V$ by 
$$\mathcal{Y} \phi = \displaystyle\sum_{j=0}^n  \sqrt{\alpha_j} \phi_j y_j \quad \text{for } \phi = ( \upphi_0, \dotsc , \upphi_n) \in \mathbb{R}^{n+1}.$$
\noindent Since the image $\mathcal{Y}(V) = \text{span} \{ y_0, \dotsc, y_n\}$ has finite dimension, the operator $\mathcal{Y}$ is compact. Its Hilbert space adjoint $\mathcal{Y}^\star: V \to \mathbb{R}^{n+1}$ satisfies 
$\langle \mathcal{Y} \phi, \psi \rangle_X = \langle \phi , \mathcal{Y}^\star \psi \rangle_{\mathbb{R}^{n+1}}$ for $\phi \in \mathbb{R}^{n+1}$ and $\psi \in V$ and is given by 
$$\mathcal{Y}^\star \psi = \begin{pmatrix}
                          \langle \psi, \sqrt{\alpha_0} y_0 \rangle_X \\
                          \vdots \\
                         \langle \psi, \sqrt{\alpha_n} y_n \rangle_X \\
                         \end{pmatrix}
 \quad \text{for } \psi \in V.$$
 
 \noindent Then, the action of $\mathcal{K} := \mathcal{Y}^\star \mathcal{Y} : \mathbb{R}^{n+1} \to \mathbb{R}^{n+1}$ 
 is given by $$\mathcal{K} \phi = \begin{pmatrix}
                                        \displaystyle\sum_{j=0}^n \sqrt{\alpha_j} \langle \sqrt{\alpha_0} y_0, y_j \rangle_X \upphi_j \\
                                        \vdots \\
              \displaystyle\sum_{j=0}^n \sqrt{\alpha_j} \langle \sqrt{\alpha_n} y_n, y_j \rangle_X \upphi_j\\
                                       \end{pmatrix} \quad \text{for } 
                                       \phi =(\upphi_0, \dots, \upphi_n) \in \mathbb{R}^{n+1}.$$

 \noindent $\mathcal{K}$ can be represented as the symmetric matrix 
 \begin{equation}\label{mathcalK}
  \mathcal{K} =     \begin{pmatrix}
                             \sqrt{\alpha_0} \sqrt{\alpha_0}  \langle y_0, y_0 \rangle_X & \hdots &  \sqrt{\alpha_0} \sqrt{\alpha_n} \langle y_0, y_n \rangle_X    \\  
                           & & \\
                           \vdots & & \vdots \\
                           & & \\
                \sqrt{\alpha_n} \sqrt{\alpha_0} \langle y_n, y_0 \rangle_X & \hdots &  \sqrt{\alpha_n} \sqrt{\alpha_n} \langle y_n, y_n \rangle_X    \\  
                                  \end{pmatrix} \in \mathbb{R}^{(n+1) \times (n+1)}.
 \end{equation}
 
 \noindent We introduce the operator $\mathcal{R}:= 
 \mathcal{Y}\mathcal{Y}^\star: V \to V$, whose action is given by
 $$\mathcal{R}\psi = \displaystyle\sum_{j=0}^n \alpha_j \langle \psi, y_j \rangle_X y_j \quad \text{ for } \psi \in V.$$

\noindent It can be shown that the operator $\mathcal{R}$
is bounded, nonnegative and self-adjoint. Since the image $\mathcal{R}(V) = 
\text{span}\{y_0,\dotsc, y_n\}$ has finite dimension, the operator $\mathcal{R}$ is compact. Therefore the Hilbert-Schmidt theorem (cf. \cite[Th. VI.16]{RS80}, for instance) can be applied which ensures the existence of a complete orthonormal basis $\{ \psi_i \}_{i=1}^\infty$ for $V$ and a 
sequence of corresponding nonnegative eigenvalues $\{\lambda_i\}_{i=1}^\infty$ 
with $$\mathcal{R} \psi_i = \lambda_i \psi_i \quad \text{with } \lambda_1 \geq \hdots \geq \lambda_d > 0 \text{ and } \lambda_i = 0 \text{ for all } i > d.$$

\noindent Likewise, one can compute the eigenvalues $\{ \lambda_i\}_{i=1}^d$ 
of $\mathcal{K}$, which coincide with the eigenvalues for $\mathcal{R}$ except for possibly zero. The corresponding orthonormal eigenvectors $\{\phi_i\}_{i=1}^d \subset \mathbb{R}^{n+1}$ of $\mathcal{K}$ are 
$$\phi_i = \displaystyle\frac{1}{\sqrt{\lambda_i}} (\mathcal{Y}^\star \psi_i) = 
 \displaystyle\frac{1}{\sqrt{\lambda_i}} \begin{pmatrix}
                                          \langle \psi_i, \sqrt{\alpha_0} y_0 \rangle_X\\
                                          \vdots\\
                                          \langle \psi_i,  \sqrt{\alpha_n} y_n \rangle_X\\
               \end{pmatrix} \in \mathbb{R}^{n+1} \quad \text{for } i = 1, \dotsc, d.$$

\noindent Thus, the functions $\{\psi_i\}_{i=1}^d$ can be determined via
 \begin{equation}\label{podBasis} \psi_i = \displaystyle\frac{1}{\sqrt{\lambda_i}} \mathcal{Y} \phi_i = 
         \displaystyle\frac{1}{\sqrt{\lambda_i}} \displaystyle\sum_{j=0}^n 
         \sqrt{\alpha_j} (\upphi_i)_j y_j \in V \quad \text{for } 
         i=1,..., d, \end{equation}

 \noindent where $(\upphi_i)_j$ denotes the $j$-th component of $\phi_i \in 
 \mathbb{R}^{n+1}$ for $i=1,...., d$.\\
 
 The following theorem states the necessary optimality conditions for problem \eqref{minPOD} and presents the POD projection error.\\
 
 \noindent \textbf{Theorem 3.1.} Let $\{\lambda_i\}_{i=1}^d$ denote the positive eigenvalues of $\mathcal{R}$, and let $\{\psi_i\}_{i=1}^d \subset V$ denote the corresponding eigenfunctions of $\mathcal{R}$. For every $\ell \in \mathbb{N}$ with $\ell \leq d$, a solution to \eqref{minPOD} is given by the eigenfunctions $\{\psi_i\}_{i=1}^\ell$ corresponding to the $\ell$ largest eigenvalues $\{\lambda_i\}_{i=1}^\ell$. Moreover the projection error is 
\begin{equation}\label{projection_error}
 \displaystyle\sum_{j=0}^n \alpha_j \left\| y_j - \displaystyle\sum_{i=1}^\ell 
 \langle y_j, \psi_i \rangle_X \; \psi_i 
 \right\|_X^2 = \displaystyle\sum_{i = \ell + 1}^d \lambda_i.
\end{equation}

 \textit{Proof.} Since $\text{span}\{ y_j\}_{j=0}^n \subset V$, the 
 proof runs analogously to the proof in \cite[Th. 3]{Vol01}.\\

\noindent The basis $\{\psi_i\}_{i=1}^\ell$ can alternatively be computed 
via singular value decomposition (SVD). The SVD of the operator $\mathcal{Y}$ is given by 
$$ \mathcal{Y} = \sum_{i=1}^d \sigma_i \sqrt{\alpha_i} \langle \cdot , \phi_i 
\rangle_{\mathbb{R}^{n+1}} \psi_i ,$$
where $\sigma_1 \geq \dotsc \geq \sigma_d > 0$ is the ordered sequence of singular values of $\mathcal{Y}$ with $\sigma_i = \sqrt{\lambda_i}$ for 
$ i = 1, \dotsc, d$. For more details we refer to \cite[Th. VI.17]{RS80}, for 
instance.\\
  
 \subsection{Numerical realization of the POD method utilizing 
 snapshots from arbitrary finite element spaces}

\begin{figure}[htbp]
  \centering
\includegraphics[scale=0.4]{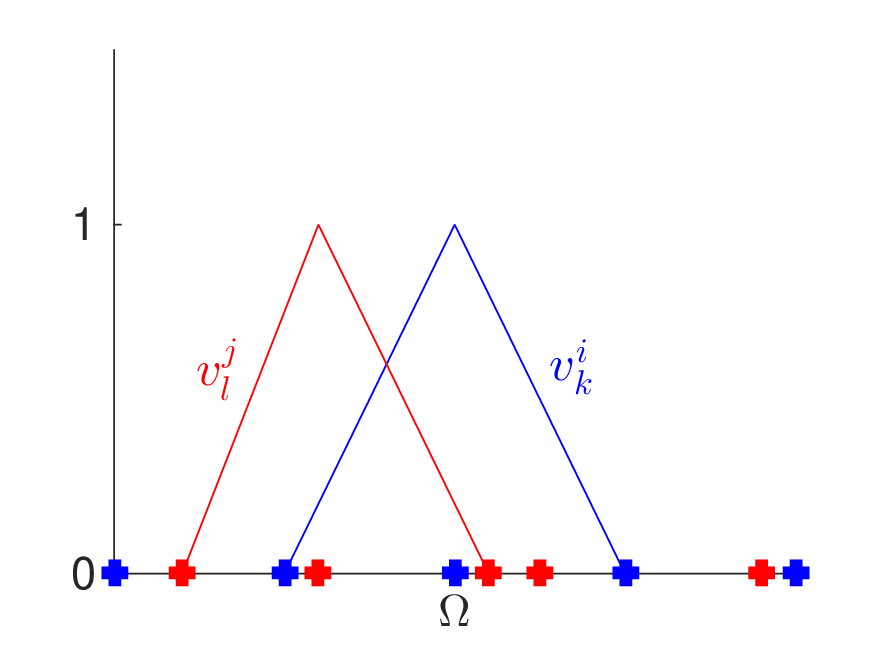} \includegraphics[scale=0.4]{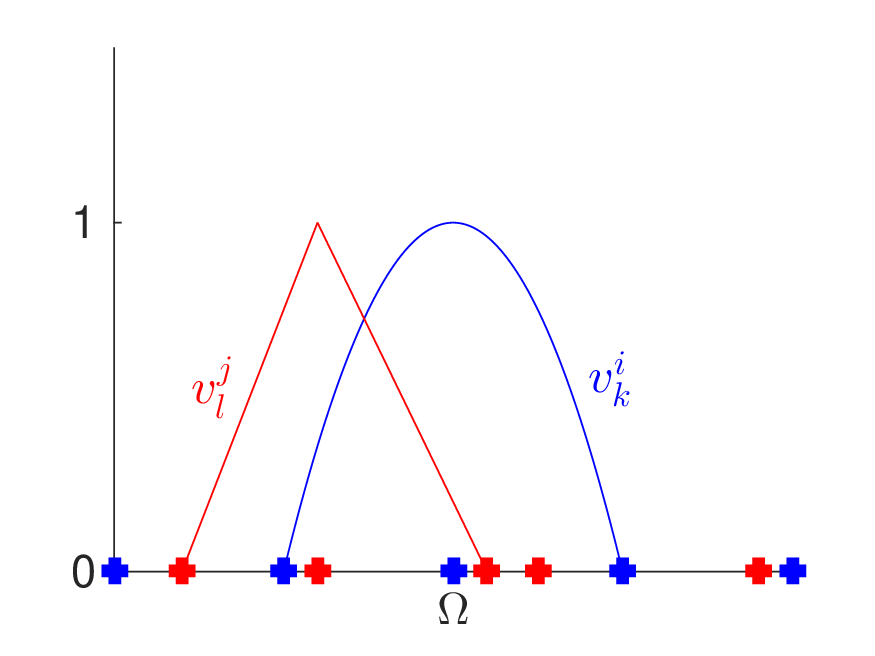} 
 \caption{\em 1D finite element basis functions $v_l^j$ and $v_k^i$ on their corresponding 
 grids. Left: both piecewise linear finite element ansatz functions. Right: piecewise 
 linear and cubic finite element ansatz functions} 
 \label{fig:twobases}
 \end{figure}

 \noindent Let us now turn our perspective to the numerical realization
 of computing a POD basis for snapshots which live in arbitrary finite element 
 spaces. For each time level $j=0, \dotsc , n$, the snapshots $\{ y_j\}_{j=0}^n$ shall be taken from different finite element spaces $\{ V_j\}_{j=0}^n $ which lie in a common Hilbert space $V$. In the fully discrete formulation of the POD method we are given the evaluation of the snapshots on their corresponding grids, i.e. we are given the vectors 
 $$\mathsf{y}^0 \in \mathbb{R}^{N_0} , \dotsc , \mathsf{y}^n \in \mathbb{R}^{N_n}$$ 
 of different lengths with $\mathsf{y}^j = (\mathsf{y}_1^j, \dotsc \mathsf{y}_{N_j}^j)^T \in \mathbb{R}^{N_j}$, for $j=0, \dotsc, n$. This is why we are not able to set up the discrete counterpart to the operator $\mathcal{R}$, which is an $N \times N$ matrix for uniform spatial discretization with $N$ nodes. Moreover, the representation of the 
 POD basis as a linear combination of the snapshots is no longer possible.\\
 To overcome this obstacle, our aim is to set up a reduced order model which can be formulated for arbitrary finite element discretizations. For this reason, we turn our attention to the matrix $\mathcal{K} \in \mathbb{R}^{(n+1) \times (n+1)}$. This matrix dimension only depends on the number of snapshots and can be computed for any underlying finite element discretization: the $ij-$th component $\mathcal{K}_{ij}$, for $i,j=0 \dotsc, n$, is given by
 $$  \sqrt{\alpha_i} \sqrt{\alpha_j} \langle y_j, y_j \rangle_X= \sqrt{\alpha_i} \sqrt{\alpha_j} \langle \displaystyle\sum_{k=1}^{N_i} \mathsf{y}_k^i v_k^i, \sum_{l=1}^{N_j} \mathsf{y}_l^j v_l^j  \rangle_X =   \sqrt{\alpha_i} \sqrt{\alpha_j} \sum_{k=1}^{N_i} \sum_{l=1}^{N_j} \mathsf{y}_k^i \mathsf{y}_l^j \langle v_k^i , v_l^j \rangle_X.$$
 Thus for any $i,j= 0, \dotsc, n$, we are able to compute the inner product 
 $\langle v_k^i , v_l^j \rangle_X$ for $k=1, \dotsc, N_i$ and $l=1, \dotsc, 
 N_j$, compare Figure \ref{fig:twobases}. We discover, then, that the calculation of the matrix $\mathcal{K}$ as well as the determination of its eigenvectors can be done for any underlying finite element discretization. Thus, the eigenvectors $\{ \phi_i\}_{i=1}^d$ of $\mathcal{K}$ are the right singular vectors of $\mathcal{Y}$ and contain the space independent time information. This fact will be used in the following to build up the reduced order model. We note that the calculation of the matrix $\mathcal{K}$ can be done for arbitrary finite element spaces, i.e. all kinds of adaptivity ($h$-, $p$- and $r$-adaptivity) can be considered and we do not need a common reference mesh. The complexity of this methods lies in the computation of the inner products of the intersections of the finite elements, whereas in \cite{URL16} the challenge lies in lifting the snapshots from different finite element spaces up to a common finite element space.

 \noindent Moreover, the calculation of the matrix $\mathcal{K} \in \mathbb{R}^{(n+1) \times (n+1)} $ (method 
 of snapshots) is favorable, since we assume the temporal dimension $n$ to be 
 far smaller than the spatial dimension(s). Due to the symmetry 
 of $\mathcal{K}$ it suffices to compute the entries on and upon the diagonal 
 of the matrix. The computations of the matrix entries of $\mathcal{K}$ can be 
 done fully in parallel. Hence, if sufficiently good hardware is available, setting 
 up the matrix $\mathcal{K}$ can be done very fast.\\

 \noindent \textbf{Example 3.2.} We choose $V = H^1(\Omega), H= L^2(\Omega)$ and
 set $X = H = L^2(\Omega)$. The triangulations of 
 $\bar{\Omega}$ for each time level $j=0, \dotsc, n$ are denoted by 
 $\{\mathcal{T}_j\}_{j=0}^n$ and the 
 finite element spaces are defined by
 $$ V_j = V(\mathcal{T}_j) = \{v \in C^0(\bar{\Omega}): v|_T \in \mathbb{P}_r(T), 
  \forall T \in \mathcal{T}_j \} \subset X, \quad j=0, \dotsc ,n,$$

\noindent where $\mathbb{P}_r$ denotes the space of polynomials of degree $r \in \mathbb{N}$. 
The computation of the $ij$-th entry 
 $\mathcal{K}_{ij} =  \sqrt{\alpha_i} \sqrt{\alpha_j} \langle y_i, y_j \rangle_{L^2(\Omega)}$
 of the matrix $\mathcal{K}$ is calculated by 
 \begin{center}
 $\begin{array}{l c l}
   \sqrt{\alpha_i} \sqrt{\alpha_j} \langle y_i, y_j \rangle_{L^2(\Omega)} 
& = &    \sqrt{\alpha_i} \sqrt{\alpha_j} \displaystyle\int\limits_\Omega y_i y_j dx \\
 & = & \sqrt{\alpha_i} \sqrt{\alpha_j} \displaystyle\sum_{k=1}^{N^i} \displaystyle\sum_{l=1}^{N^j} \mathsf{y}_k^i \mathsf{y}_l^j \displaystyle\int\limits_\Omega v_k^i v_l^j dx \\
 & = &  \sqrt{\alpha_i} \sqrt{\alpha_j} \displaystyle\sum_{k=1}^{N^i} \displaystyle\sum_{l=1}^{N^j} \mathsf{y}_k^i \mathsf{y}_l^j \left( \displaystyle\sum_{T \in \mathcal{T}_i} \; \displaystyle\sum_{\bar{T} \in \mathcal{T}_j} \; \displaystyle\int\limits_{T \cap \bar{T}} v_k^i v_l^j dx \right). \\
 \end{array}$
 \end{center}

 \noindent Computations become simpler when using nested grids. In this case, 
  the intersection of two arbitrary $n$-dimensional simplices coincides either 
  with the smaller simplex, or is a common edge simplex, or has no overlap.\\

\begin{figure}[H]
 \centering
 \includegraphics[scale=0.15]{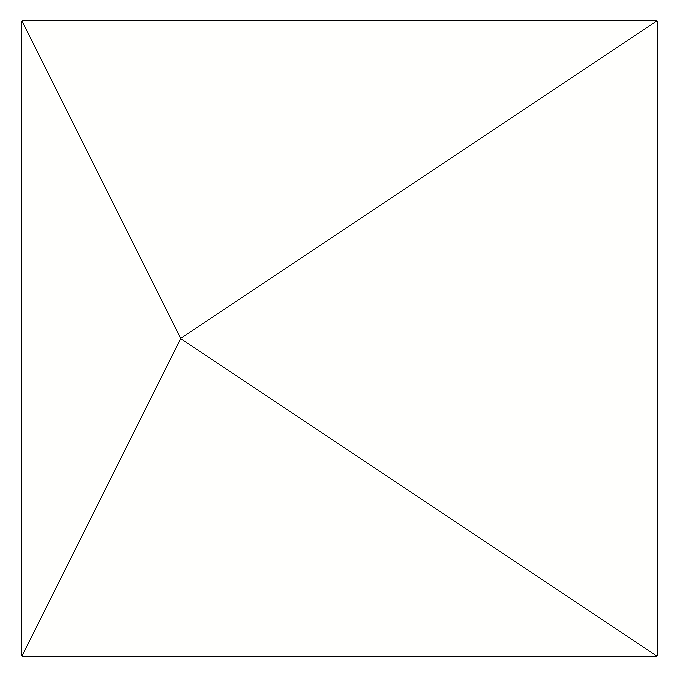} \includegraphics[scale=0.15]{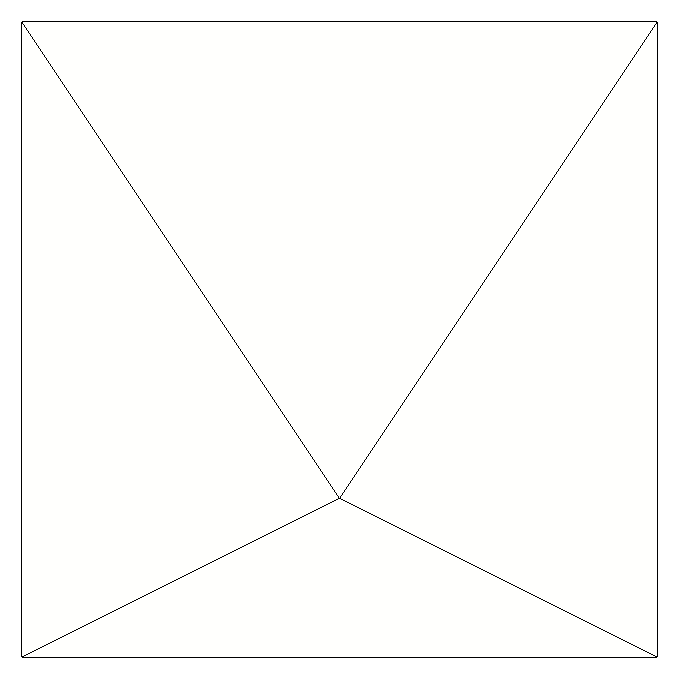} \includegraphics[scale=0.15]{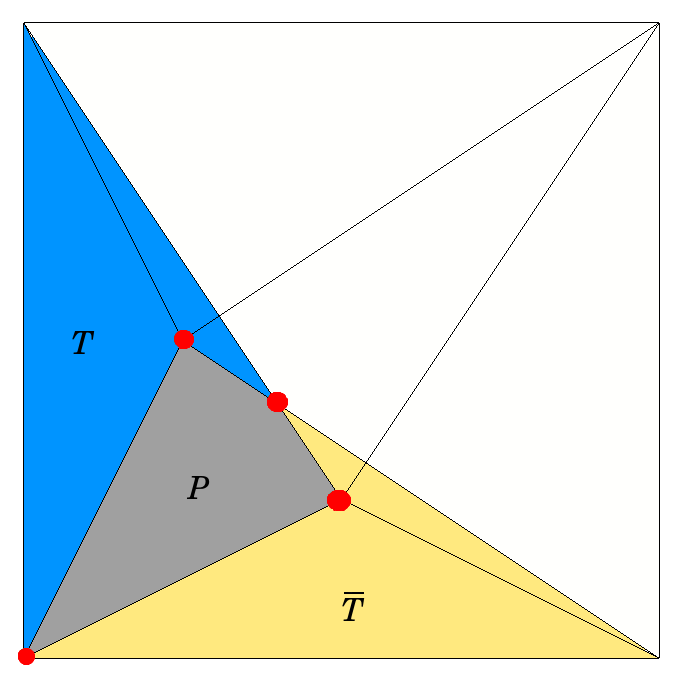}
 \caption{\em Two finite element meshes (left, middle) and the overlapping mesh (right)}
  \label{fig:meshesradaptive}
 \end{figure}

 \noindent  Example 3.2 reveals the challenge to deal with integrals of type
 \begin{equation}\label{int_polygon}
  \int_{P} y_i y_j dx
 \end{equation}
 over cut finite elements (polyhedra) $P = T \cap \bar{T}$ with $T \in \mathcal{T}_i$ and $\bar{T} \in \mathcal{T}_j$, which involves the integration of functions defined on different (non-matching) meshes, compare Figure \ref{fig:meshesradaptive}. 
 
 \noindent Our numerical realization is strongly based on \cite{MLL13}, where a similar challenge is investigated in the context of multimesh methods. There, the major challenges are identified as (i) collision detection (find intersecting simplices), (ii) mesh intersection (detect intersection interface) and (iii) integration on complex polyhedra. In the numerical 
 Example 6.3 we make use of built-in FEniCS \cite{ABHJKLRW15, LMW12} tools to handle the issues (i) and (ii). For (iii) the integration over the cut elements, a subtriangulation can be an option. But as pointed out in \cite{MLL13} in the context of a three-dimensional example, a subtetrahedralization of an arbitrary polyhedron is challenging and additional vertices might have to be added. Therefore, an alternative approach is outlined which is based on a boundary representation of the integrals, compare \cite{LR82, Mir96} for example. Exemplarily, in the case of a two-dimensional domain with piecewise linear finite element discretization, the integrals of type \eqref{int_polygon} can be computed as respectively weighted sum of the integrals
\begin{equation}\label{basic_integrals}
 \int_P x_0^2 dx, \quad \int_P x_1^2, \quad \int_P x_0 x_1 dx, \quad \int_P x_0 dx, \quad \int_P x_1 dx, \quad \int_P 1 dx,
\end{equation}
with $x=(x_0,x_1)$ using Stoke's formula
\begin{equation}\label{stokes}
 \int_P f(x) dx = \frac{1}{2+q} \sum_{i=1}^m \frac{b_i}{\|a_i\|} \int_{E_i} f(x) d\sigma_i,
\end{equation}
where $\{E_i\}_{i=1}^m$ denote the edges of the polyhedron $P \subset \mathbb{R}^2$, $a_i^T x = b_i$ is the hyperplane in which $E_i$ lies, and $f$ is a polynomial of degree $q$. Note that the line integrals on the edges $E_i$ can be computed using standard Gauss quadrature, for example.

\section{POD reduced order modeling}

\subsection{POD reduced order modeling for arbitrary finite element discretizations}

\noindent In this section, we stay in the infinite-dimensional setting of the POD method and set up the POD reduced order model utilizing snapshots with arbitrary finite element discretizations. This perspective allows us to determine the mode coefficients of the POD Galerkin ansatz for arbitrary underlying finite element discretizations. Suppose for given snapshots $y_0 \in V_0, \dotsc, y_n \in V_n$ we have computed the matrix $\mathcal{K}$, \eqref{mathcalK}, with $\mathcal{K}_{ij} =  \sqrt{\alpha_i} \sqrt{\alpha_j} \langle y_i, y_j \rangle_X$, for 
$i,j = 0, \dotsc , n$ as well as its $\ell$ largest eigenvalues $\{ \lambda_i\}_{i=1}^\ell$ and corresponding eigenvectors $\{ \phi_i \}_{i=1}^\ell \subset \mathbb{R}^{n+1}$ of low rank $\ell$ with $\ell \leq n+1$,
according to the strategy presented in Section 3. The POD basis $\{\psi_i\}_{i=1}^\ell$ is then given by \eqref{podBasis}, i.e
$$\psi_i = \frac{1}{\sqrt{\lambda_i}} \mathcal{Y} \phi_i \quad \text{ for } i=1, \dotsc, \ell.$$
This POD basis is utilized in order to compute a reduced order model for \eqref{P}.
For this reason, we make the POD Galerkin ansatz
\begin{equation}\label{PODGalerkin}
y^\ell (t) = \displaystyle\sum\limits_{i=1}^\ell \eta_i (t) \psi_i = 
\displaystyle\sum\limits_{i=1}^\ell \eta_i (t) \frac{1}{\sqrt{\lambda_i}} 
\mathcal{Y} \phi_i \quad \text{for all } t \in [0,T],
\end{equation}
as an approximation for $y(t)$, with the Fourier coefficients 
$$\eta_i(t) = \langle y^\ell(t), \psi_i \rangle_X = \langle y^\ell(t), 
\frac{1}{\sqrt{\lambda_i}} \mathcal{Y} \phi_i \rangle_X$$ 
for $i = 1, \dotsc , \ell$. Inserting $y^\ell$ into \eqref{P} and choosing $V^\ell = \text{span} \{\psi_1, \dotsc, \psi_\ell \} \subset V$ as the test space leads to the system
 \begin{equation}\label{ROM_sys_yell}
  \left\{
 \begin{array}{rcl}
  \displaystyle\frac{d}{dt} \langle y^\ell (t), \psi \rangle_H + a(y^\ell(t), \psi) + \langle \mathcal{N}(y^\ell(t)), \psi \rangle_{V',V}& = & \langle f(t), \psi \rangle_{V',V}\\
    \langle y^\ell(0), \psi \rangle_H & = & \langle g, \psi \rangle_H
    \end{array}
 \right.
 \end{equation}
 for all $\psi \in V^\ell$ and for almost all $t \in (0,T]$. Utilizing the ansatz \eqref{PODGalerkin}, we can write \eqref{ROM_sys_yell} as an $\ell$-dimensional ordinary differential equation system for the POD mode coefficients $\{\eta_i (t) \}_{i=1}^\ell \subset \mathbb{R}$:
 \begin{equation}\label{ROM_sys}
  \left\{
 \begin{array}{rcl}
   \displaystyle\sum_{j=1}^\ell {\dot \eta}_j (t) \langle \psi_i , \psi_j \rangle_H + \displaystyle\sum_{j=1}^\ell \eta_j (t) \langle \mathcal{A} \psi_j, \psi_i \rangle_{V',V} + \langle \mathcal{N}(y^\ell(t)), \psi_i \rangle_{V',V} & = & \langle f(t), \psi_i \rangle_{V',V}  \\
    & &  \hspace{1cm}  \text{ for } t \in (0,T], \\   
   \displaystyle\sum_{j=1}^\ell \eta_j(0) \langle \psi_i, \psi_j \rangle_H & = & \langle g, \psi_i \rangle_H , \\   
 \end{array}
 \right.
 \end{equation}
 for $i = 1, \dotsc , \ell$. Note that $\langle \psi_i , \psi_j \rangle_H = \delta_{ij}$ if we choose $X=H$ in the context of Section 3.\\

 \noindent Since we want to construct a reduced order model which can be built and solved for arbitrary finite element discretizations, we rewrite system \eqref{ROM_sys} utilizing the identity \eqref{podBasis}. Then, the system \eqref{ROM_sys} can be written as
 \begin{equation}\label{ROM_sys_phi}
  \left\{
 \begin{array}{rcl}
   \displaystyle\sum_{j=1}^\ell {\dot \eta}_j (t) \displaystyle\frac{1}{\sqrt{\lambda_i}} \displaystyle\frac{1}{\sqrt{\lambda_j}} \langle \mathcal{Y} \phi_i , \mathcal{Y} \phi_j \rangle_H + \displaystyle\sum_{j=1}^\ell \eta_j (t)  \frac{1}{\sqrt{\lambda_i}} \frac{1}{\sqrt{\lambda_j}}\langle 
   \mathcal{A} \mathcal{Y} \phi_j, \mathcal{Y} \phi_i \rangle_{V',V} & &  \\
   + \displaystyle\frac{1}{\sqrt{\lambda_i}} \langle \mathcal{N}(y^\ell(t)), \mathcal{Y} \phi_i \rangle_{V',V} & = & \displaystyle\frac{1}{\sqrt{\lambda_i}} \langle f(t), \mathcal{Y} \phi_i \rangle_{V',V} \\
   & & \hspace{0.7cm} \text{ for } t \in (0,T], \\   
    \displaystyle\sum_{j=1}^\ell \eta_j(0) \frac{1}{\sqrt{\lambda_i}} 
   \frac{1}{\sqrt{\lambda_j}} \langle \mathcal{Y} \phi_i, \mathcal{Y} \phi_j \rangle_H  & = & \displaystyle\frac{1}{\sqrt{\lambda_i}} \langle g, \mathcal{Y} \phi_i \rangle_H.  \\   
 \end{array}
 \right.
 \end{equation}
 \noindent In order to write \eqref{ROM_sys_phi} in a compact matrix-vector 
 form, let us introduce the diagonal matrix
 $$D \in \mathbb{R}^{\ell \times \ell} \text{ with } D = \text{diag}\left(\frac{1}{\sqrt{\lambda_1}}, \dotsc, \frac{1}{\sqrt{\lambda_\ell}}\right).$$ 
 From the first $\ell$ eigenvectors $\{ \phi_i\}_{i=1}^\ell$ 
 of $\mathcal{K}$ we build the matrix
 $$\Phi \in \mathbb{R}^{(n+1) \times \ell} \text{ by } \Phi = [\phi_1 \; | \; \dotsc \; | \; \phi_\ell].$$
 Then, the system \eqref{ROM_sys_phi} can be written as the system
 
 \begin{equation}\label{ROM}
  \left\{
 \begin{array}{rcll}
   D \Phi^T \mathcal{K} \Phi D \; {\dot \eta}(t) + D \Phi^T \mathcal{Y}^\star \mathcal{A} \mathcal{Y} \Phi D \; \eta (t) + DN(\eta(t)) & = & DF(t) & \text{ for } t \in (0,T], \\   
   D \Phi^T \mathcal{K} \Phi D \;  \eta (0) & = & D \bar{\eta}_0, & \\   
 \end{array}
 \right.
 \end{equation}
 
 \noindent for the vector-valued mapping  $\eta (t) = (\eta_1 (t), \dotsc, \eta_\ell(t))^T : [0,T] \to \mathbb{R}^\ell$. Note that the right hand side $F(t)$ and the initial condition $ \bar{\eta}_0$ are given by
$$ \left( F(t) \right)_i =  \langle f(t), \mathcal{Y} \phi_i \rangle_{V',V} = 
\langle \mathcal{Y}^\star f(t), \phi_i \rangle_{\mathbb{R}^{n+1}} $$
and 
$$ \left( \bar{\eta}_0 \right)_i = \langle g, \mathcal{Y} \phi_i \rangle_H =  \langle \mathcal{Y}^\star g, \phi_i \rangle_{\mathbb{R}^{n+1}},$$
for $i = 1, \dotsc, \ell$, respectively. Their calculation can be done explicitly 
for any arbitrary finite element discretization. For a given function $w \in V$ 
(for example $w=f(t)$ or $w = g$) with finite element discretization $w = \sum_{i=1}^{N_w} \mathsf{w}_i \varphi_i$, nodal basis $\{\varphi_i\}_{i=1}^{N_w} \subset V$ and appropriate mode coefficients $\{ \mathsf{w}_i \}_{i=1}^{N_w}$ we can compute 
 $$\left( \mathcal{Y}^\star w\right)_j = \langle w, y_j \rangle_X = 
 \langle \sum_{i=1}^{N_w} \mathsf{w}_i \varphi_i , \sum_{k=1}^{N_j} \mathsf{y}_k^j 
 v_k^j \rangle_X = \sum_{i=1}^{N_w} \sum_{k=1}^{N_j} \mathsf{w}_i \mathsf{y}_k^j 
 \langle \varphi_i , v_k^j \rangle_X, \quad \text{ for } j = 0, \dotsc, n $$ 
 where $ y_j \in V_j$ denotes the $j$-th snapshot. Again, for any $i = 1, \dotsc, 
 N_w$ and $k = 1, \dotsc, N_j$, the computation of the inner product $\langle \varphi_i , v_k^j \rangle_X$ can be done explicitly.\\

 \noindent Obviously, for linear evolution equations the POD reduced order model \eqref{ROM} can be set up and solved utilizing snapshots with arbitrary finite element discretizations. The computation of the nonlinear component $N(\eta(t))$ needs particular attention. In the following Section 4.2 we discuss the options to treat the nonlinearity.\\

 \noindent \textbf{Remark 4.1.} Of course, the derivation of a POD surrogate model \eqref{ROM} for \eqref{P} as explained above, is also applicable for other classes of differential equations like elliptic PDEs, for example.\\

\noindent {\textbf{Time-discrete reduced order model}}\\

\noindent In order to solve the reduced order system \eqref{ROM_sys_yell} numerically, we apply the implicit Euler method for time discretization and use for simplicity the same temporal grid $\{t_j\}_{j=0}^n$ \eqref{timegrid} as for the snapshots. It is also possible to use a different time grid, cf. \cite{KV02}. 
The time-discrete reduced order model reads 
 
\begin{equation}\label{ROM_timediscrete}
  \left\{
 \begin{array}{rcll}
  \left\langle \displaystyle\frac{y_j^\ell-y_{j-1}^\ell}{\Delta t_j} , \psi \right\rangle_H + a(y_j^\ell,\psi) + \langle \mathcal{N}(y_j^\ell), \psi \rangle_{V',V} & = & \displaystyle\frac{1}{\Delta t_j} \displaystyle\int_{t_{j-1}}^{t_j} \langle f_j, \psi \rangle_{V',V} & \text{ for all } \psi \in V^\ell,\\   
   \langle y_0^\ell, \psi \rangle_H & = & \langle g, \psi \rangle_H, &\\   
 \end{array}
 \right.
 \end{equation}
 \noindent for $j=1, \dotsc, n$ or equivalently
 \begin{equation}\label{ROM_infPOD_timediscrete}
  \left\{
 \begin{array}{rcll}
 D\Phi^T \mathcal{K} \Phi D \; \left( \displaystyle\frac{\eta^j - \eta^{j-1}}{\Delta t_j} \right) + D \Phi^T \mathcal{Y}^\star \mathcal{A} \mathcal{Y} \Phi D \; \eta^j + DN(\eta^j) & = & DF_j & \text{ for } j = 1, \dotsc, n,\\
 D \Phi^T \mathcal{K} \Phi D \; \eta^0 & = & D \bar{\eta}_0. &
 \end{array}
 \right.
 \end{equation}

 \subsection{Discussion of the nonlinear term $\boldsymbol{DN(\eta(t))}$}
 
\noindent Let us now consider the computation of the nonlinear term $DN(\eta(t)) \in 
\mathbb{R}^\ell$ of the POD-ROM \eqref{ROM}. It holds true
\begin{equation*}
 \begin{array}{r c l}
    (DN(\eta(t)))_k & = & \langle \mathcal{N}(y^\ell(t)), \psi_k \rangle_{V',V} \\[0.1cm]
    & = & \langle \mathcal{N}(\sum_{i=1}^\ell \eta_i(t) \psi_i), \psi_k \rangle_{V',V}
 \end{array}
\end{equation*}
for $k=1, \cdots, \ell$. It is well-known that the evaluation of nonlinearities in the reduced order model context in computationally expensive. To make this clear, let us assume, we are given a uniform finite element discretization with $N$ degrees of freedom. Then, in the fully discrete setting, the nonlinear term has the form 
$$ \Psi^T W \mathcal{N}(\Psi \eta(t)), $$
where $\Psi = \left[ \psi_1 \; | \; \cdots \; | \; \psi_\ell \right] \in \mathbb{R}^{N \times \ell}$ is the matrix in which the POD modes are stored columnwise and $W \in \mathbb{R}^{N \times N}$ is a weighting matrix related to the 
utilized inner product. Hence, the treatment of the nonlinearity requires the expansion of $y^\ell(t) = \Psi \eta(t) \in \mathbb{R}^N$ in the full space, then the 
nonlinearity can be evaluated and finally the result is projected back to the POD space. Obviously, this means that the reduced order model is not fully independent of the high order dimension and efficient simulation cannot be guaranteed. Therefore, it is convenient to seek for hyper reduction, i.e. for a treatment of the nonlinearity where the model evaluation cost is related to the low dimension $\ell$. A possible remedy which is commonly used, is given discrete by empirical interpolation methods, for which we refer to \cite{CS10} for DEIM and to \cite{DG15} for Q-DEIM. Another option is given by \cite{AK16} which investigates nonlinear model reduction via dynamic mode decomposition. Furthermore, in \cite{Wan15} nonlinear model reduction is realized by replacing the nonlinear term by its interpolation into the finite element space. The treatment of the nonlinearity is done in \cite{AWWB08} by the missing point estimation method and in \cite{NPP08} by best the points interpolation method.\\
All these methods have in common, that a common reference mesh is needed. Since we 
want to circumvent the interpolation of the snapshots onto a common grid and we 
want to set up the reduced order model for arbitrary finite element discretizations, we have to go a different way.\\

\noindent One option is to utilize EIM \cite{BMNP04} which is a hyper reduction technique formulated in a continuous setting. Alternatively, we can linearize and project the nonlinearity onto the POD space. For this approach, let us consider the linear reduced order system given by:

\begin{equation}\label{P_linearized}
  \left\{
 \begin{array}{rcl}
   \frac{d}{dt} \langle y^\ell(t), \psi \rangle_H + a(y^\ell(t),\psi) + \langle \mathcal{N}(\bar{y}(t)), \psi\rangle_{V',V}
    & = & \langle f(t), v \rangle_{V',V},\\   
   \langle y^\ell(0), \psi \rangle_H & = & \langle g, v \rangle_H, \\   
 \end{array}
 \right.
 \end{equation}
 For a given state $\bar{y}(t)\in V$, this linear evolution problem \eqref{P_linearized} can be set up and solved explicitly without spatial interpolation. In the numerical examples in Section 6, we take the finite element solution as given state in each time step, i.e. $\bar{y}(t_j) = y_j$.\\
Furthermore, the linearization of the reduced order model \eqref{ROM_sys_yell} can be considered:
\begin{equation}\label{P_linearized2}
  \left\{
 \begin{array}{rcl}
   \frac{d}{dt} \langle y^\ell(t), \psi \rangle_H + a(y^\ell(t),\psi) + \langle \mathcal{N}(\bar{y}(t)) + \mathcal{N}_y(\bar{y}(t))(y^\ell-\bar{y})(t), \psi \rangle_{V',V}    & = & \langle f(t), v \rangle_{V',V},\\   
   \langle y^\ell(0), \psi \rangle_H & = & \langle g, v \rangle_H, \\   
 \end{array}
 \right.
 \end{equation}
where $\mathcal{N}_y$ denotes the Fréchet derivative. This linearized problem is of interest in the context of optimal control, where it occurs in each iteration level within the SQP method, see \cite{HPUU09}, for example. Choosing the finite element solution as given state in each time instance leads to:
\begin{equation*}
\begin{array}{l c l}
 \langle \mathcal{N}(y_j),\psi_i \rangle_{V',V} & = & \frac{1}{\sqrt{\lambda_i}} \sum_{k=0}^n \sqrt{\alpha_k}(\upphi_i)_k \langle \mathcal{N}(y_j), y_k \rangle_{V',V}\\[0.2cm]
 \langle \mathcal{N}_y(y_j) y^\ell(t_j),\psi_i\rangle_{V',V} & = & \langle \mathcal{N}_y(y_j)(\sum_{k=1}^\ell \eta_k(t_j)\psi_k),\psi_i \rangle_{V',V} \\[0.2cm]
  & = & \sum_{k=1}^\ell \eta_k(t_j) \frac{1}{\sqrt{\lambda_k}} \frac{1}{\sqrt{\lambda_i}} \sum_{p=0}^n \sum_{r=0}^n \sqrt{\alpha_p} \sqrt{\alpha_r}  (\upphi_k)_p (\upphi_i)_r \langle \mathcal{N}_y(y_j) y_p, y_r \rangle_{V',V} \\[0.2cm]
 \langle \mathcal{N}_y(y_j)y_j,\psi_i \rangle_{V',V} & = & \frac{1}{\sqrt{\lambda_i}} \sum_{k=0}^n \sqrt{\alpha_k}  (\upphi_i)_k  \langle \mathcal{N}_y(y_j) y_j, y_k \rangle_{V',V}
\end{array} 
\end{equation*}
Finally, we approximate the nonlinearity $DN(\eta^j)$ in \eqref{ROM_infPOD_timediscrete} by
$$ (DN(\eta^j))_k \approx \langle \mathcal{N}(y_j) + \mathcal{N}_y(y_j)(y^\ell(t_j)- y_j),\psi_k \rangle_{V',V} $$
which can be written as
$$ DN(\eta^j) \approx D \Phi^T \mathsf{N}^j + D \Phi^T \mathtt{N}_y^j \Phi D \eta^j - D\Phi^T \mathsf{N}_y^j$$
where 
$$\mathsf{N}^j=  \begin{pmatrix}
\langle \mathcal{N}(y_j), \sqrt{\alpha_0}  y_0\rangle_{V',V}\\
\vdots\\
\langle \mathcal{N}(y_j), \sqrt{\alpha_n}  y_n\rangle_{V',V}\\
        \end{pmatrix} \in \mathbb{R}^{n+1}, \mathsf{N}_y^j=  \begin{pmatrix}
\langle \mathcal{N}_y(y_j) y_j, \sqrt{\alpha_0} y_0 \rangle_{V',V}\\
\vdots\\
\langle \mathcal{N}_y(y_j) y_j, \sqrt{\alpha_n} y_n \rangle_{V',V}\\
        \end{pmatrix} \in \mathbb{R}^{n+1}$$ and\\[0.2cm]
        $\mathtt{N}_y^j=  \begin{pmatrix}
\langle \mathcal{N}_y(y_j)\sqrt{\alpha_0} y_0, \sqrt{\alpha_0} y_0 \rangle_{V',V} & \hdots & \langle \mathcal{N}_y(y_j)\sqrt{\alpha_n} y_n, \sqrt{\alpha_0} y_0 \rangle_{V',V}\\
\vdots & & \vdots\\
\langle \mathcal{N}_y(y_j) \sqrt{\alpha_0} y_0, \sqrt{\alpha_n} y_n \rangle_{V',V} & \hdots & \langle \mathcal{N}_y(y_j)\sqrt{\alpha_n} y_n, \sqrt{\alpha_n} y_n \rangle_{V',V}\\
        \end{pmatrix} \in \mathbb{R}^{(n+1) \times (n+1)}$\\[0.4cm]
        
 \noindent For weakly nonlinear systems this approximation may be sufficient, depending on the problem and its goal. A great advantage of linearizing the semilinear partial differential equation is that only linear equations need to be solved which leads to a further speedup. However, if a more precise approximation is desired or necessary, we can think of approximations including higher order terms, like quadratic approximation, see e.g. \cite{Che99} and \cite{RW03}, or 
 Taylor expansions, see e.g. \cite{Phi00}, \cite{Phi03} and \cite{FZCZF04}. Nevertheless, the efficiency of higher order approximations is limited due to 
 growing memory and computational costs.


 \subsection{Expressing the POD solution in the full spatial domain}

 \noindent Having determined the solution $\eta(t)$ to \eqref{ROM}, we can set up 
 the reduced solution $y^\ell(t)$ in a continuous framework:
 \begin{equation}\label{solROM}
y^\ell (t) = \sum_{i=1}^\ell \eta_i(t)
\left(\frac{1}{\sqrt{\lambda_i}} \sum_{j=0}^n  \sqrt{\alpha_j} 
 (\upphi_i)_j  y_j \right).
\end{equation}
Now, let us turn to the fully discrete formulation of \eqref{solROM}. 
For a time-discrete setting, we introduce for simplicity the same temporal 
grid $\{t_j\}_{j=0}^n$ as for the snapshots. Let us recall the spatial discretization 
of the snapshots \eqref{yFE} utilizing arbitrary finite elements
$$ y_j  = \displaystyle\sum_{k=1}^{N_j} \mathsf{y}_k^j v_k^j \quad \text{ for } j = 0, \dotsc, n.$$

\noindent Let $\{ Q_r^j \}_{r=1}^{M_j}$ denote an arbitrary
set of grid points for the reduced system at time level $t_j$. 
The fully discrete POD solution can be computed by evaluation:
\begin{equation}\label{yPOD}
 y^\ell(t_j, Q_r^j ) = \sum_{i=1}^\ell \eta_i (t_j) \left( 
\frac{1}{\sqrt{\lambda_i}} \sum_{s=0}^n  \sqrt{\alpha_s} (\upphi_i)_s  (\sum_{k=1}^{N_s} 
\mathsf{y}_k^s v_k^s(Q_r^j)) \right) 
\end{equation}
for $j = 0, \dotsc, n$ and $r = 1, \dotsc, M_j$. This allows us to use any grid 
for expressing the POD solution in the full spatial domain. For example, we can use the same node points at time level $j$ for the POD simulation as we have used for the snapshots, i.e. for $j = 0, \dotsc, n$ it holds $M_j = N_j$ and $ Q_r^j= P_k^j$ for all $r,k = 1, \dotsc, N_j$. Another option can be to choose 
$$\{Q_r^j\}_{r=1}^{M_j} = \displaystyle\bigcup_{j=0}^n \bigcup_{k=1}^{N_j} \{P_k^j\},$$ 
i.e. the common finest grid. Obviously, a special and probably the easiest case
concerning the implementation is to choose snapshots which are expressed with respect to the same finite element basis functions and utilize the common finest grid for the simulation of the reduced order system, which is proposed by \cite{URL16}. After expressing the adaptively sampled snapshots with respect to a common finite element space, the subsequent steps coincide with the common approach of taking snapshots which are generated without adaptivity. Then, expression \eqref{yPOD} simplifies to
\begin{equation}\label{POD_full}
  y^\ell(t_j, P_r^j) = \sum_{i=1}^\ell \eta_i(t_j) \left( 
\frac{1}{\sqrt{\lambda_i}} \sum_{s=0}^n  \sqrt{\alpha_s} (\upphi_i)_s 
\mathsf{y}^s \right).
\end{equation}

\section{Error analysis for the reduced order model}

\noindent For the validation of the approximation quality of the POD reduced order
model, we are interested in analyzing the error between the POD solution and the true solution. Our aim is to estimate the expression 
\begin{equation*}
 \sum_{j=0}^n \alpha_j \| y(t_j) - y_j^\ell \|_H^2
\end{equation*}
where $\{y(t_j)\}_{j=0}^n \subset V$ denotes the true solution for \eqref{P} at time instances $\{t_j\}_{j=0}^n$ and $\{y_j^\ell\}_{j=0}^n$ is the solution to the time-discrete reduced order model \eqref{ROM_timediscrete}, i.e. 
\begin{equation}\label{ROM_sol}
 y_j^\ell = \sum_{i=1}^\ell \eta_i^j \psi_i,
\end{equation}
 with the POD basis $\{\psi_i\}_{i=1}^\ell$ computed from the snapshots $\{ y_j \}_{j=0}^n$ from arbitrary finite element spaces, i.e. $ y_j \in V_j$ for $j = 0, \dotsc, n$, as explained in Section 3. The weights $\{\alpha_j\}_{j=0}^n$ are trapezoidal time weights, see \eqref{weights_alpha}.
We choose $X=V$ in the context of Section 3. Let us introduce the orthogonal projection $\mathcal{P}^\ell: V \to V^\ell$ by 
\begin{equation*}
 \mathcal{P}^\ell v = \displaystyle\sum_{i=1}^\ell \langle v, \psi_i \rangle_V \psi_i \quad \text{ for all } v \in V.
\end{equation*}

\noindent It holds true $\| \mathcal{P}^\ell\|_{\mathcal{L}(V)} = 1$, where $\mathcal{L}(V)$ is the space of linear bounded operators from $V$ to $V$. The subsequent calculations follow closely the proofs of \cite[Thm. 4.7]{KV02} \cite[Thm. 3.2.5]{Vol13} and \cite[Th. 7]{KV01}. In our situation, we compute the POD basis corresponding to the fully discrete snapshots $\{ y_j \}_{j=0}^n$ utilizing adaptive finite element spaces, whereas in \cite{KV02}, \cite{Vol13} and \cite{KV01}, the POD basis is computed from snapshots corresponding to the solution trajectory at the given time instances. For this reason, in the following estimation an additional term corresponding to the error for the spatial discretization will appear.\\

\noindent We make use of the decomposition
\begin{equation}\label{decomposition}
 y(t_j) - y_j^\ell = y(t_j) -  y_j + y_j - \mathcal{P}^\ell y_j + \mathcal{P}^\ell y_j - \mathcal{P}^\ell y(t_j) + \mathcal{P}^\ell y(t_j) - y_j^\ell = \eta_j +  \varrho_j  + \zeta_j + \vartheta_j
\end{equation}
for $j=0, \dotsc , n$, where $\eta_j := y(t_j) - y_j, \varrho_j := y_j - \mathcal{P}^\ell y_j, \zeta_j = \mathcal{P}^\ell y_j - \mathcal{P}^\ell y(t_j)$ and
$\vartheta_j := \mathcal{P}^\ell y(t_j) - y_j^\ell$. The term $\eta_j$ is the discretization error. We utilize the decomposition
\begin{equation*}
  \eta_j = y(t_j) - y_j = y(t_j) -  \bar{y}_j + \bar{y}_j - y_j = E_t^j + E_h^j
 \end{equation*}
 where $ \bar{y}_j$ denotes the solution to the time-discrete problem \eqref{P_timediscrete}. By $E_t^j := y(t_j) - \bar{y}_j$ we denote the global time discretization error and $E_h^j :=   \bar{y}_j - y_j$ is the global spatial discretization error. It is
\begin{equation*}
 \| E_h^j \|_H \leq \max_{j=0, \dotsc, n} \| E_h^j\|_H =: \varepsilon_h
\end{equation*}
and 
\begin{equation*}
 \| E_t^j \|_H \leq \max_{j=0, \dotsc, n} \| E_t^j\|_H =: \varepsilon_t.
\end{equation*}
Since we use the implicit Euler method for time integration, it is $\varepsilon_t = 
\mathcal{O}(\Delta t)$ with $\Delta t:= \max_{j = 0, \dotsc, n} \Delta t_j$. Therefore, we can estimate 
\begin{equation}\label{rest_eta}
 \sum_{j=0}^n \alpha_j \|y(t_j) - y_j\|_H^2 \leq \sum_{j=0}^n \alpha_j \| E_t^j + E_h^j \|_H^2 \leq 2 \sum_{j=0}^n \alpha_j ((\Delta t)^2 + \varepsilon_h^2) \leq 2 T ( (\Delta t)^2 + \varepsilon_h^2).
\end{equation}
\noindent Moreover, we have 
\begin{equation}\label{zeta_j}
 \sum_{j=0}^n \alpha_j \| \zeta_j \|_H^2 =  \sum_{j=0}^n \alpha_j \| \mathcal{P}^\ell  y_j - \mathcal{P}^\ell y(t_j) \|_H^2 \leq \| \mathcal{P}^\ell \|_{\mathcal{L}(H)}^2 \sum_{j=0}^n \alpha_j \| \eta_j \|_H^2.
\end{equation}

\noindent The term $\varrho_j$ is the projection error of the snapshot $y_j$ projected onto the POD space $V^\ell$. Using \eqref{norm_estimation}, the weighted sum of all projection errors is given by the sum of the neglected eigenvalues \eqref{projection_error}, i.e. 
\begin{equation}\label{varrho}
 \sum_{j=0}^n \alpha_j \| \varrho_j \|_H^2 = \sum_{j=0}^n \alpha_j \|  y_j - 
 \sum_{i=1}^\ell \langle  y_j, \psi_i \rangle_V \psi_i \|_H^2  \leq c_v \sum_{j=0}^n \alpha_j \|y_j - \sum_{i=1}^\ell \langle y_j, \psi_i\rangle_V \psi_i \|_V^2 \leq c_v \sum_{i=\ell+1}^d \lambda_i.
\end{equation}
\noindent It remains to estimate the term $\vartheta_j$ which is the error between the projection of the true solution $y(t_j)$ at time instance $t_j$ onto the POD 
space $V^\ell$ and the time-discrete ROM solution $y_j^\ell$ to \eqref{ROM_sol}. 
With the use of the notation $\bar{\partial} \vartheta_j = (\vartheta_j - \vartheta_{j-1})/\Delta t_j$ for $j = 1, \dotsc, n$, we get

\begin{equation*}
  \left.
 \begin{array}{rcl}
 \langle \bar{\partial} \vartheta_j, \psi \rangle_H & = & \langle \mathcal{P}^\ell \left(\displaystyle\frac{y(t_j) - y(t_{j-1})}{\Delta t_j}\right) - \displaystyle\frac{y_j^\ell - y_{j-1}^\ell}{\Delta t_j}, \psi \rangle_H \\   
  & = & \langle \mathcal{P}^\ell \left(\displaystyle\frac{y(t_j) - y(t_{j-1})}{\Delta t_j}\right) +  \mathcal{N}(y_j^\ell) - f_j, \psi \rangle_H + a(y_j^\ell, \psi) \\
  & = & \langle \mathcal{P}^\ell \left(\displaystyle\frac{y(t_j) - y(t_{j-1})}{\Delta t_j} \right) - \displaystyle\frac{y(t_j) - y(t_{j-1})}{\Delta t_j} + \mathcal{N}(y_j^\ell) - \mathcal{N}(y(t_j)), \psi \rangle_H + a(y_j^\ell-y(t_j), \psi) \\
  & = & \langle z_j   + \mathcal{N}(y_j^\ell) - \mathcal{N}(y(t_j)), \psi \rangle_H + a(y_j^\ell-y(t_j), \psi)
 \end{array}
 \right.
 \end{equation*}

\noindent for $\psi \in V^\ell$ with $z_j := \mathcal{P}\left(\displaystyle \frac{y(t_j) - y(t_{j-1})}{\Delta t_j}\right) - \displaystyle\frac{y(t_j) - y(t_{j-1})}{\Delta t_j}$. With the choice $\psi = \vartheta_j$ and the use of the identity
$$ 2 \langle u-v, u \rangle = \| u \|^2 - \|v \|^2 + \|u-v\|^2 $$
we obtain
\begin{equation*}
 \| \vartheta_j\|_H^2 \leq \| \vartheta_{j-1} \|_H^2 + 2 \Delta t_j \left( \beta \| y_j^\ell - y(t_j)\|_H
 + \| z_j \|_H + \| \mathcal{N}(y_j^\ell) - \mathcal{N}(y(t_j))\|_H \right) \| \vartheta_j \|_H,
\end{equation*}
where we have utilized \eqref{bilinear}. We assume that $\mathcal{N}$ is Lipschitz continuous, i.e. there exists $L>0$ such that 
$$ \| \mathcal{N}(y_j^\ell) - \mathcal{N}(y(t_j)) \|_H \leq L \| y_j^\ell - y(t_j) \|_H \quad \text{ for } j = 1, \dotsc, n.$$
Applying Young's inequality we find
\begin{equation*}
\| \vartheta_j \|_H^2 \leq \| \vartheta_{j-1} \|_H^2 + \Delta t_j ( c_1 \| \varrho_j \|_H^2 + c_2 \| \vartheta_j \|_H^2 + \| z_j \|_H^2 + c_1\| \eta_j \|_H^2 + c_1 \| \zeta_j \|_H^2) 
\end{equation*}
with the constants $c_1 := \beta + L$ and $c_2 := 5 (\beta + L)+1$. Under the assumption that $\Delta t$ is sufficiently small, we conclude
\begin{equation}
 \| \vartheta_j \|_H^2 \leq e^{2c_2j \Delta t} \left( \| \vartheta_1 \|_H^2 + \sum_{k=1}^j \Delta t_k (\| z_k\|_H^2 + c_1 \| \varrho_k \|_H^2 + c_1 \| \eta_k \|_H^2 + c_1 \| \zeta_k \|_H^2) \right).
\end{equation}
\noindent For more details on this, we refer to \cite{KV02} and \cite{Vol13}. We choose the initial condition for \eqref{ROM_timediscrete} such that $\vartheta_0 = \mathcal{P}^\ell y(t_0) - y_0^\ell = \mathcal{P}^\ell  g - y_0^\ell = 0$.\\

\noindent Next, we estimate the term involving $z_k$. It holds true
\begin{equation*}
\left.
 \begin{array}{rcl}
  \| z_k \|_H^2 & = & \| \mathcal{P}^\ell \left( \displaystyle\frac{y(t_k) - y(t_{k-1})}{\Delta t_k} \right) - \displaystyle\frac{y(t_k)- y(t_{k-1})}{\Delta t_k} \|_H^2\\
  & = & \| \mathcal{P}^\ell \left( \displaystyle\frac{y(t_k) - y(t_{k-1})}{\Delta t_k} \right) - \mathcal{P}^\ell \dot{y}(t_k) + \mathcal{P}^\ell \dot{y}(t_k) - \dot{y}(t_k) + \dot{y}(t_k) - \displaystyle\frac{y(t_k)- y(t_{k-1})}{\Delta t_k} \|_H^2\\
 & \leq &  2 \| \mathcal{P}^\ell \|_{\mathcal{L}(H)}^2 \| \displaystyle\frac{y(t_k)-y(t_{k-1})}{\Delta t_k} - \dot{y}(t_k) \|_H^2 + 2 \| \mathcal{P}^\ell \dot{y}(t_k) - \dot{y}(t_k) \|_H^2 + 2 \| \dot{y}(t_k) - \displaystyle\frac{y(t_k) - y(t_{k-1})}{\Delta t_k} \|_H^2 \\
  & \leq & c_3 \| w_k \|_H^2 + 2 \| \mathcal{P}^\ell \dot{y}(t_k) - \dot{y}(t_k) \|_H^2
  \end{array}
\right.
\end{equation*}
with $c_3=2+2 \| \mathcal{P}^\ell \|_{\mathcal{L}(H)}^2$ and $w_k := \dot{y}(t_k) - \displaystyle\frac{y(t_k) - y(t_{k-1})}{\Delta t_k} $, which can be estimated as 
\begin{equation*}
\sum_{k=1}^j \Delta t_k \| w_k \|_H^2 \leq \frac{(\Delta t)^2}{3} \| \ddot{y} \|_{L^2(0,t_j,H)}^2.
 \end{equation*}
For more details on this, we refer to \cite{KV02} and \cite{Vol13}.

\noindent Finally, we can summarize the estimation for the term involving $\vartheta_j$ by
\begin{equation*}
 \| \vartheta_j \|_H^2 \leq  c_4 \left(\displaystyle\sum_{k=1}^n \alpha_k ( \| \mathcal{P}^\ell 
 \dot{y}(t_k) - \dot{y}(t_k) \|_H^2 + \| \varrho_k \|_H^2 + \| \eta_k\|_H^2 + \| \zeta_k \|_H^2) + (\Delta t)^2 \| \ddot{y} \|_{L^2(0,T,H)}^2 
 \right)
\end{equation*}
with $c_4 := e^{2c_2 T} \max\{\frac{c_3}{3},4,2c_1\}$ and thus it is
\begin{equation}\label{vartheta}
\left.
 \begin{array}{rcl}
 \displaystyle\sum_{j=0}^n \alpha_j \| \vartheta_j \|_H^2 & \leq &  c_4 T (\displaystyle\sum_{j=0}^n \alpha_j \| \mathcal{P}^\ell 
 \dot{y}(t_j) - \dot{y}(t_j) \|_H^2 + (\Delta t)^2 \| \ddot{y} \|_{L^2(0,T,H)}^2 
 + \sum_{i=\ell +1}^d \lambda_i \\
 & & + 2T(1+\|\mathcal{P}^\ell\|_{\mathcal{L}(H)}^2)((\Delta t)^2 + \varepsilon_h^2))\\
   \end{array}
 \right.
\end{equation}

\noindent \textbf{Theorem 5.1.} Let $\{y(t_j)\}_{j=0}^n$ denote the solution to problem \eqref{P} at the time grid $\{t_j\}_{j=0}^n$ and $y_j^\ell$ is the solution to \eqref{ROM_timediscrete}. Let the nonlinear operator $\mathcal{N}$ be Lipschitz continuous with Lipschitz constant $L$ and the maximal time 
step $\Delta t:= \max_{j=0, \dotsc, n} \Delta t_j$ be sufficiently small. Furthermore, we assume $\ddot{y}(t)$ to be bounded on $[0,T]$. We choose the initial condition for \eqref{ROM_timediscrete} such that $ \mathcal{P}^\ell g = y_0^\ell $ is fulfilled. Then, there exists a constant $C =  C(T, c_v, \|P^\ell\|_{\mathcal{L}(H)}^2, \beta, L, \|\ddot{y}\|_{L^2(0,T,H)}^2) > 0$ such that
\begin{equation}\label{errorestimation}
\displaystyle\sum_{j=0}^n \alpha_j \| y(t_j) - y_j^\ell \|_H^2 \leq C \left( 
(\Delta t)^2 + \varepsilon_h^2 + \displaystyle\sum_{i=\ell+1}^d \lambda_i + 
\displaystyle\sum_{j=0}^n \alpha_j \| \mathcal{P}^\ell \dot{y}(t_j) - \dot{y}(t_j) 
\|_H^2 \right),
\end{equation}
where $c_v, \beta$ are from \eqref{norm_estimation},\eqref{bilinear}, the quantity $\varepsilon_h  := \max_{j=0, \dotsc, n} \| y_j - \bar{y}_j \|_H $ refers to the global spatial discretization error and $\bar{y}_j$ is the solution to \eqref{P_timediscrete} at time instance $t_j$.\\

\noindent \textit{Proof.} Utilizing the decomposition \eqref{decomposition} we infer
$$ \displaystyle\sum_{j=0}^n \alpha_j \| y(t_j) - y_j^\ell \|_H^2 \leq 2 \displaystyle\sum_{j=0}^n \alpha_j (\| \eta_j \|_H^2 + \| \varrho_j \|_H^2 + \|\zeta_j\|_H^2+ \| \vartheta_j\|_H^2).$$
Together with \eqref{rest_eta}, \eqref{zeta_j}, \eqref{varrho} and \eqref{vartheta}
 this leads to the claim. \hspace{3.5cm} $\square$\\
 
 \noindent \textbf{Remark 5.2.} i) The last term in \eqref{errorestimation} can be avoided by adding time derivatives to the snapshot set (more specifically, finite difference approximations of time derivatives), cf. \cite{KV01}. In the recent work \cite{Sin14}, a new error bound is proved which avoids the last term in \eqref{errorestimation} and does not need to include time derivative data in the POD snapshot set.\\
 ii) If we choose $V=H^1(\Omega)$ and $H=L^2(\Omega)$ and utilize a static piecewise linear finite element discretization with $h$ being the diameter of the triangles, then $\varepsilon_h =  \mathcal{O}(h^2)$. In the case of adaptively refined spatial grids $\varepsilon_h$ can be estimated by the prescribed error tolerance, if e.g. residual based a-posteriori error estimation is applied.\\

 \section{Numerical Realization and Examples}

 \noindent For all numerical examples, we choose $\Omega \subset \mathbb{R}^2$ as open and bounded domain and utilize conformal, piecewise linear and continuous finite elements for spatial discretization.\\ 
 In Example 6.1 and 6.2, we utilize an $h$-adaptive concept, which controls the spatial discretization error in the energy-norm. As error indicator for the adaptive strategy with respect to space we utilize at each time step the jump across the edges, which reflects the main contribution of the classical residual-based error estimation for elliptic problems, see e.g. (2.18)-(2.19) in \cite{AO11}. The error indicator is given by
 \begin{equation}\label{errorind}
  \eta_E = \sqrt{h_E} \| \left[\nabla y_j \right]_E \cdot \nu_E \|_ L^2(E) \quad \text{ for all } E \in \mathcal{E},
 \end{equation}
where $h_E$ is the length of the edge, $\left[ \; . \; \right]_E$ denotes the jump of the function across the edge $E$, $\nu_E$ is the outward normal derivative on the edge $E$ and $\mathcal{E}$ denotes the collection of all (interior) edges in the current triangulation. As refinement rule we use the bisection by newest vertex based on \cite{Che08}.\\

In Example 6.1 and 6.2, we utilize structured, hierarchical and nested grids. These numerical test cases illustrate our approach to set up and solve a POD reduced order model utilizing snapshots with adaptive spatial discretization, which is explained in Section 3 and 4 in the specific case in which a nested mesh structure is at hand.  Thus, the computations benefit from the fact that the intersection of two triangles coincides either with the smaller triangle, or is a common edge, or has no overlap. We compare our approach to the use of a uniform mesh, where the mesh size coincides with the fineness of the smallest triangle in the adaptive mesh. The aim of this comparison is to investigate numerically how the inclusion of spatial adaptivity for the snapshot discretization affects the accuracy of the POD reduced order solution compared to using a uniform mesh where no spatial adaptation is performed. Furthermore, we compare our approach to the use of a finest mesh following \cite{URL16} concerning numerical efficiency and accuracy. For this, the practical numerical concept works as follows. In a full dimensional simulation, $h$-adaptive snapshots are generated at the time instances $\{ t_j\}_{j=0}^n$. At the same time we carry along a reference grid with the simulation, which coincides with the computational grid at initial time and which is only refined in the same manner as the computational grid. In this way, the reference grid becomes the finest mesh (i.e. the overlay of all computational grids) at the end of the snapshot generation. Then, the snapshots are expressed with respect to the finite element basis functions corresponding to the finest mesh and the usual POD procedure is carried out, choosing $X=L^2(\Omega)$ as Hilbert space.\\
 In Example 6.3 we realize the numerical computation of the correlation 
 matrix $\mathcal{K}$ \eqref{mathcalK} as described in Section 3.2 for non-nested 
 meshes. In this case, the overlap of two meshes leads to cut finite elements, which are convex polygons of more than three node points. Thus, numerical 
 computations become more involved.\\
 All coding is done in C++ and we utilize FEniCS \cite{ABHJKLRW15, LMW12} for the solution of the differential equations and ALBERTA \cite{SS05} for dealing with hierarchical meshes. We run the numerical tests on a compute server with 24 CPU 
 kernels and 512 GB RAM.\\

  
 \noindent \textbf{Example 6.1: Linear heat equation.} We consider the Example 2.3 \eqref{heat} of a heat equation with homogeneous Dirichlet boundary condition and set $c\equiv 0$ such that the equation becomes linear. The spatial domain is chosen as $\Omega = [0,1] \times [0,1] \subset \mathbb{R}^2$, the time interval is $[0,T]=[0,1.57]$. We construct an example in such a way that we know the analytical solution. It is given by 
 $$ y(t,x) = r(t,x) \cdot \left( s_1(t,x) - s_2(t,x) \right)$$ 
 with\\

 \noindent $ r(t,x) = \frac{50000 \cdot x_0 \cdot (1-x_0) \cdot (0.5+\cos(t) \cdot(x_0-0.5)-\sin(t)\cdot (x_1-0.5))^4 \cdot \frac{1}{t+1} \cdot (1-(0.5+\cos(t) \cdot (x_0-0.5)-\sin(t) \cdot (x_1-0.5)))^4)}{1+1000 \cdot (\cos(t) \cdot (x_0-0.5)-\sin(t) \cdot (x_1-0.5))^2},$\\

\noindent $ s_1(t,x) = \frac{10000 \cdot x_1 \cdot (1-x_1) \cdot (0.5+\sin(t) \cdot (x_0-0.5)+\cos(t) \cdot (x_1-0.5))^2 \cdot (0.5-\sin(t) \cdot (x_0-0.5)-\cos(t) \cdot(x_1-0.5))^2}{1+100*((0.5+\sin(t)\cdot(x_0-0.5)+\cos(t) \cdot (x_1-0.5))-0.25)^2},$\\
 
\noindent  $ s_2(t,x) =  \frac{10000 \cdot x_1 \cdot (1-x_1) \cdot (0.5+\sin(t)\cdot (x_0-0.5)+\cos(t) \cdot (x_1-0.5))^2 \cdot (0.5-\sin(t) \cdot (x_0-0.5)-\cos(t)\cdot (x_1-0.5))^2}{(1+100*((0.5+\sin(t) \cdot (x_0-0.5)+\cos(t)\cdot (x_1-0.5))-0.75)^2}.$\\              
     
\begin{figure}[htbp]
 \centering
  \includegraphics[scale=0.12]{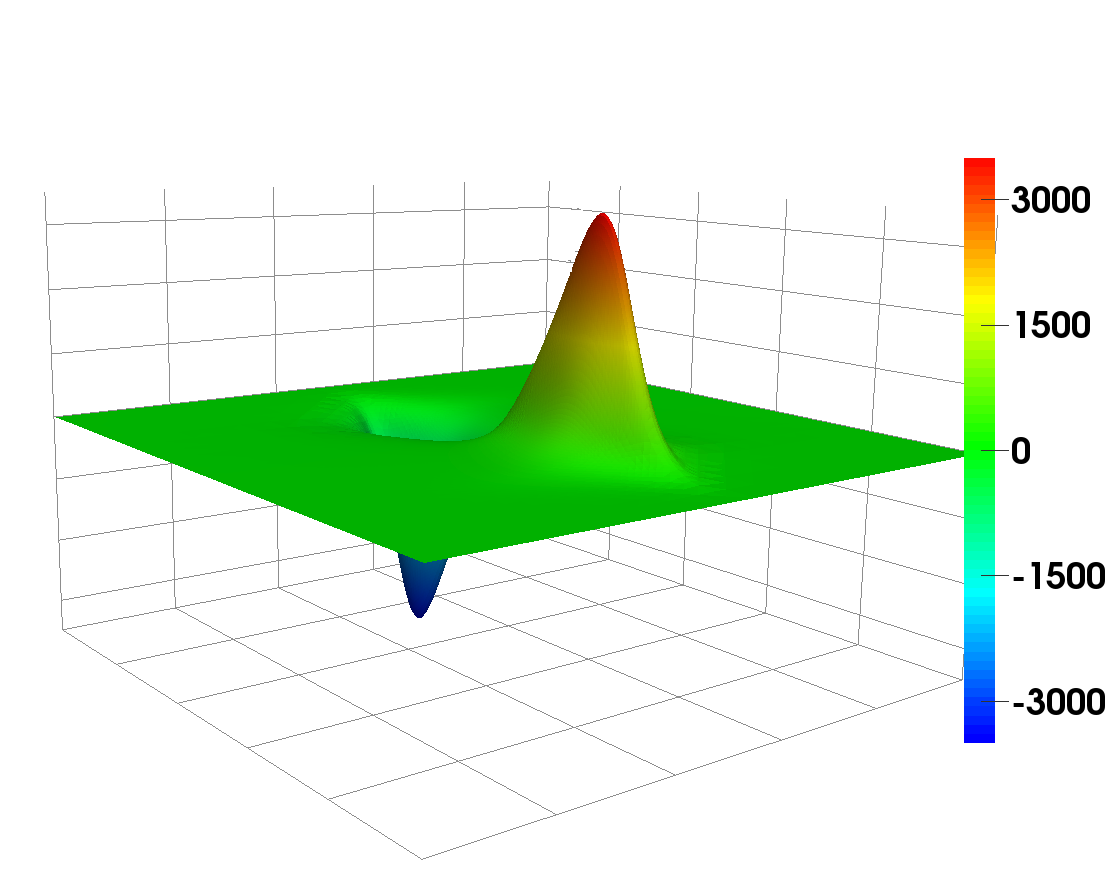} \hspace{0.1cm} \includegraphics[scale=0.12]{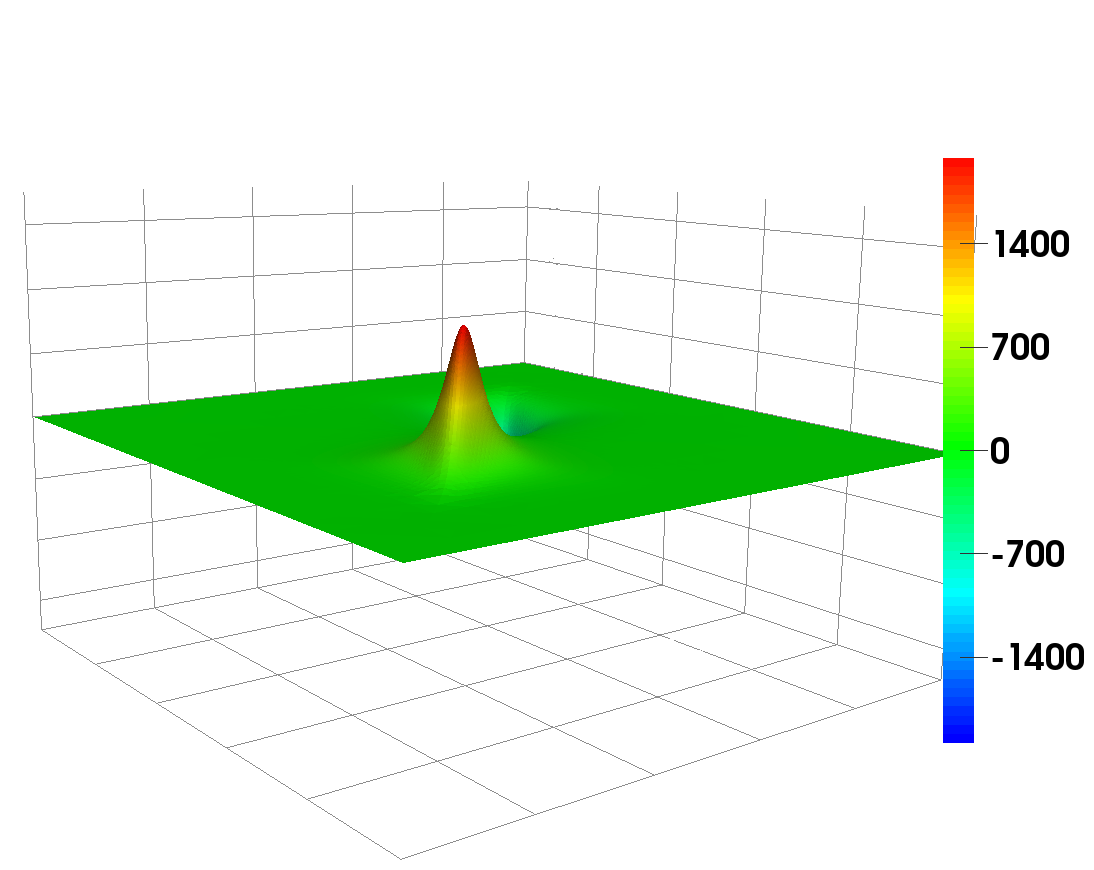}  \hspace{0.1cm} \includegraphics[scale=0.12]{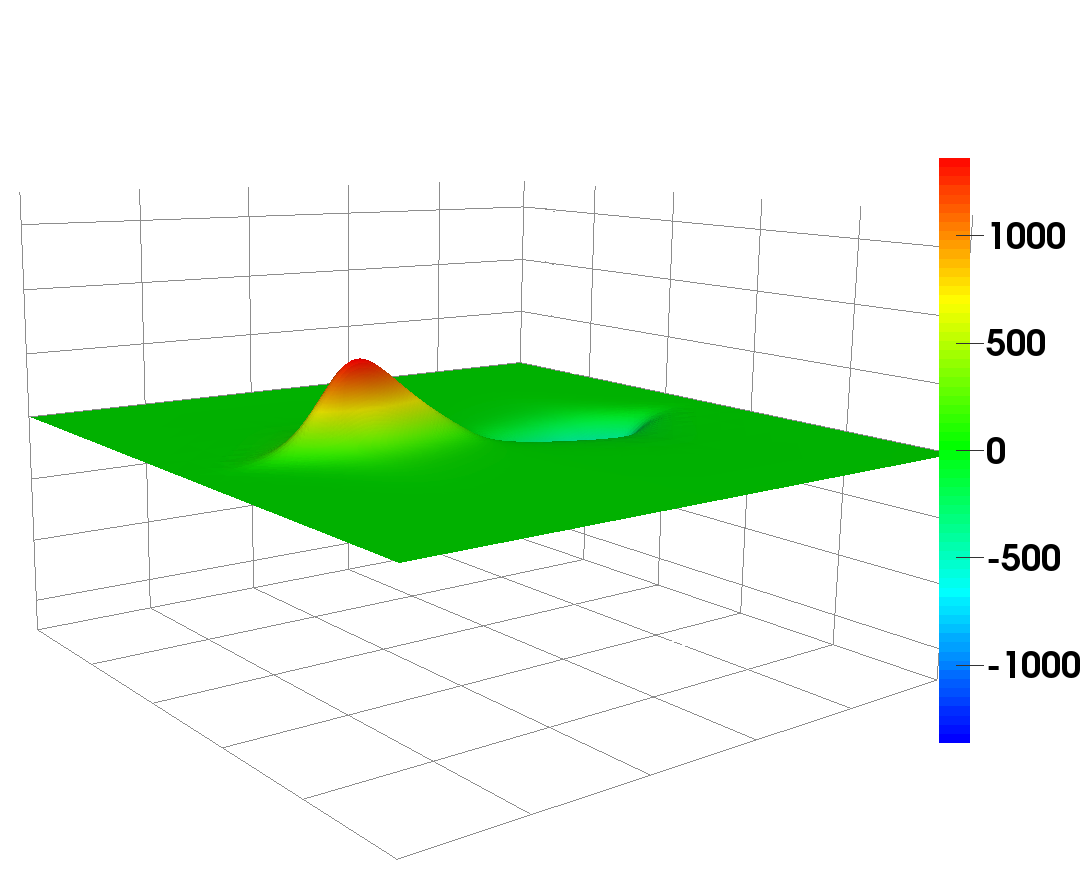} \\       
 \hspace{-0.5cm}  \includegraphics[scale=0.12]{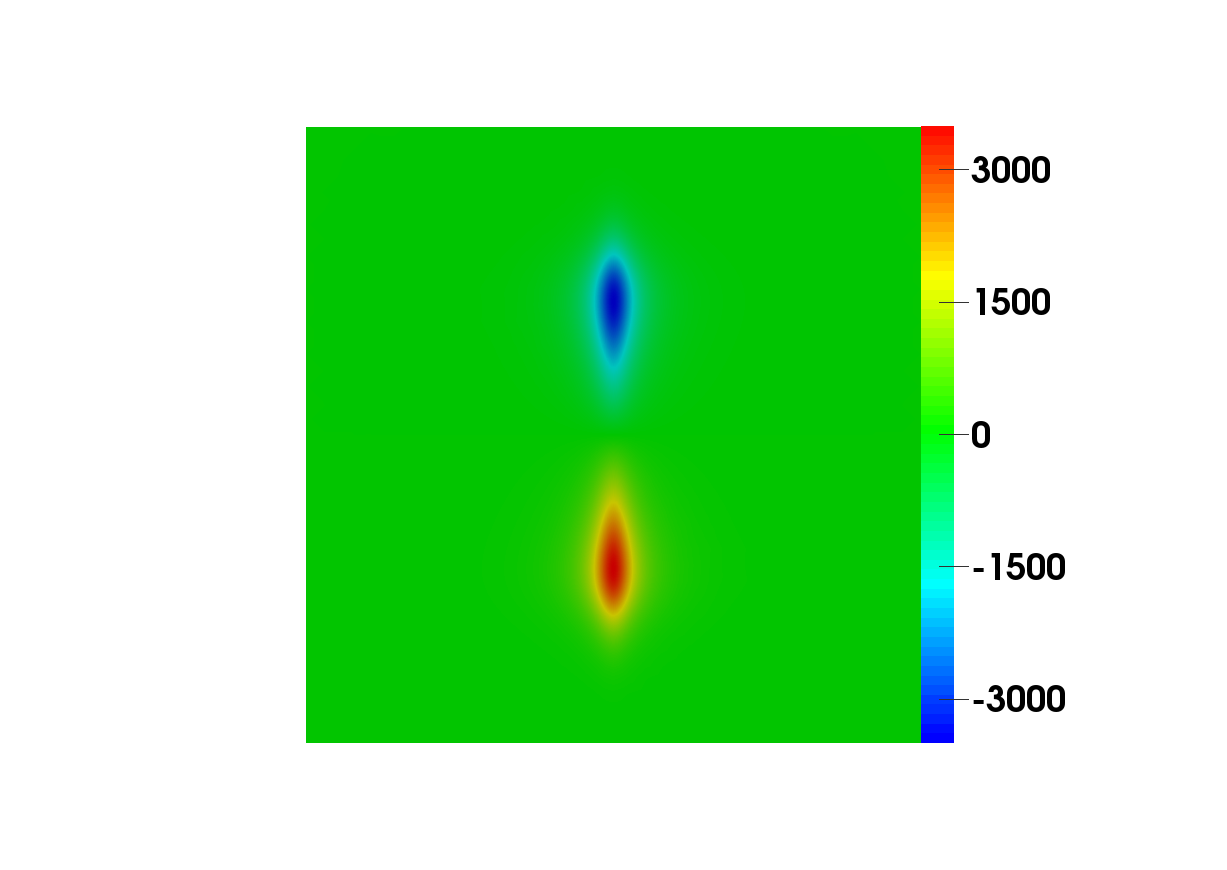} \hspace{-0.5cm} \includegraphics[scale=0.12]{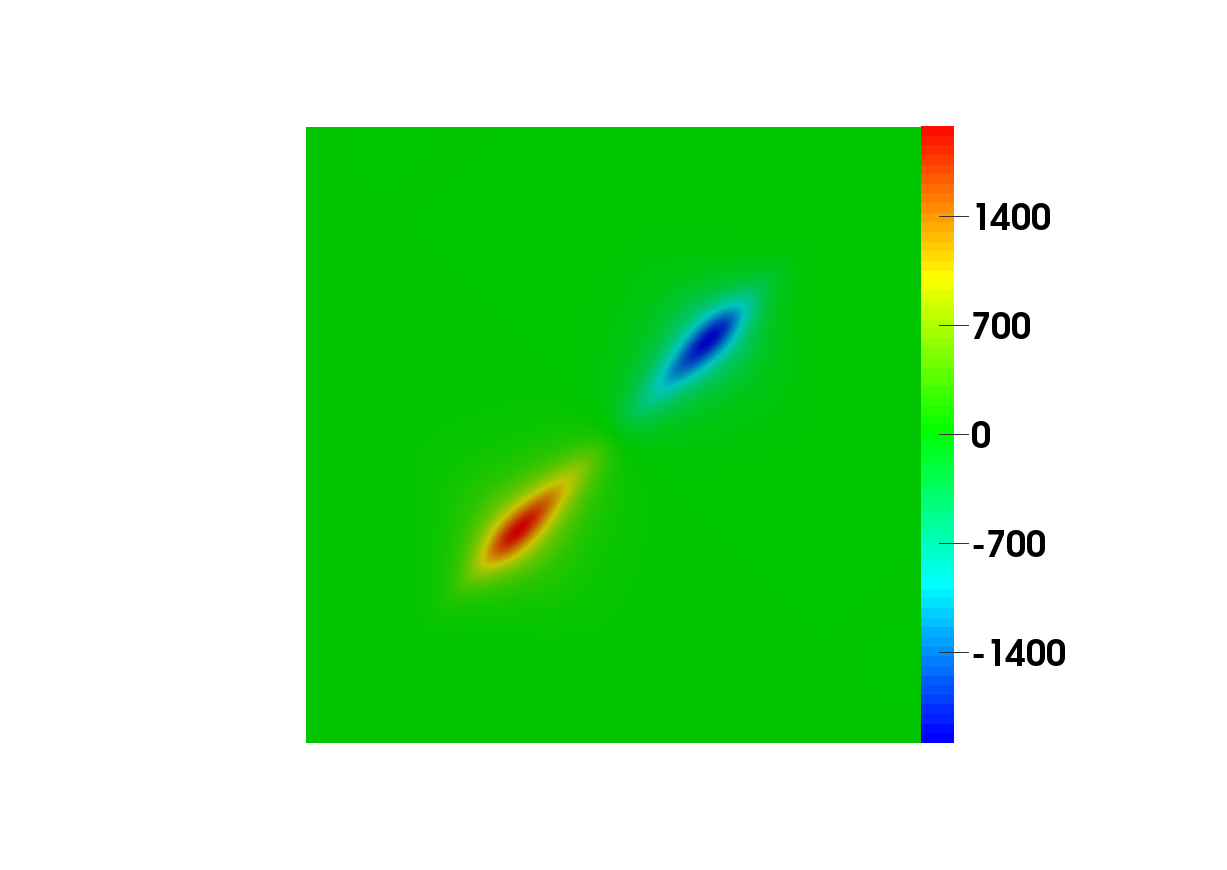} \hspace{-0.5cm} \includegraphics[scale=0.12]{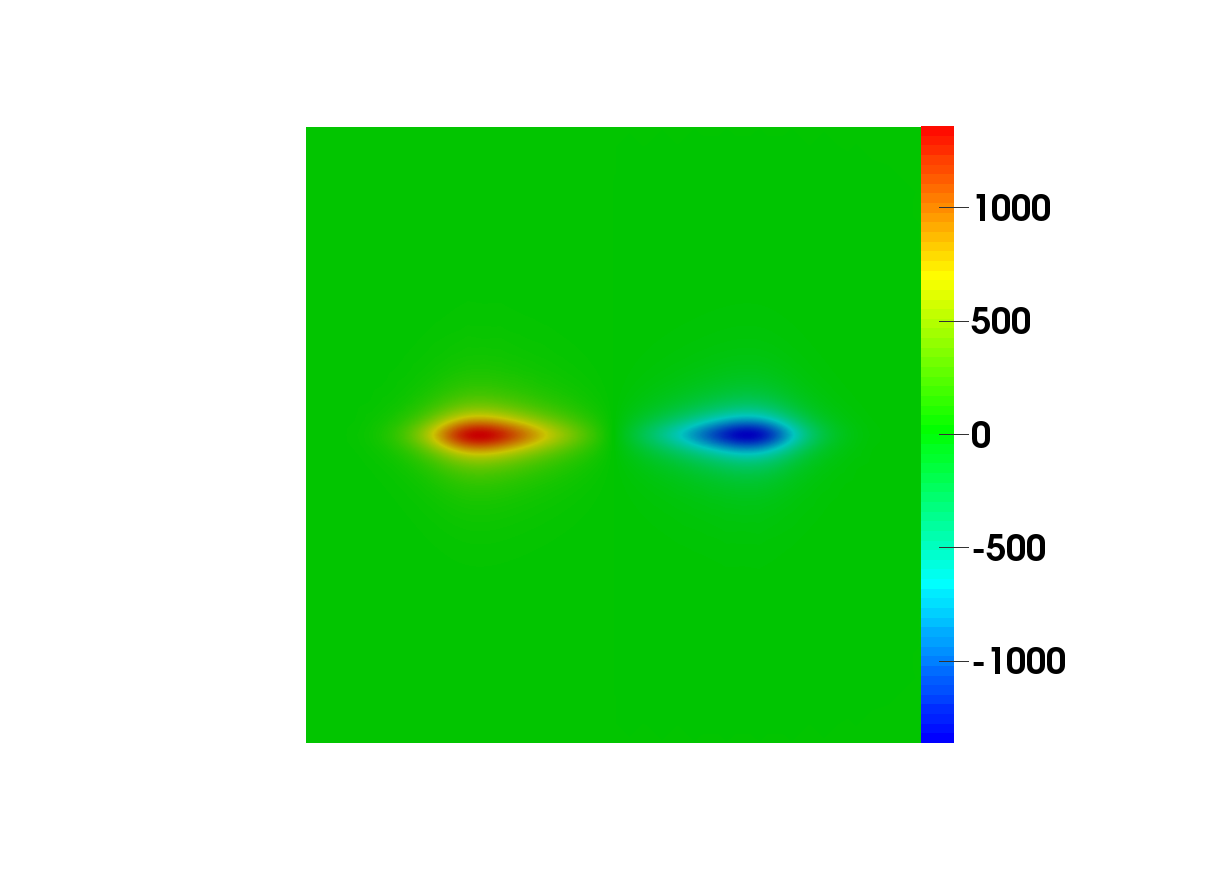}        
 \caption{\em Example 6.1: Surface plot (top) and view from above (bottom) of the analytical solution of \eqref{heat} at $t=t_0$ (left), $t=T/2$ (middle) and $t=T$ (right)}
 \label{fig:Heat_true}
 \end{figure}

\noindent The forcing term $f$ and the initial condition $g$ are chosen accordingly. For the temporal discretization we introduce the uniform time grid by $$t_j = j \Delta t$$ for $j=0, \dotsc, 1570$ with $\Delta t = 0.001$. The analytical solution at three different time points is shown in Figure \ref{fig:Heat_true}.

\begin{figure}[H]
 \centering
 \includegraphics[scale=0.12]{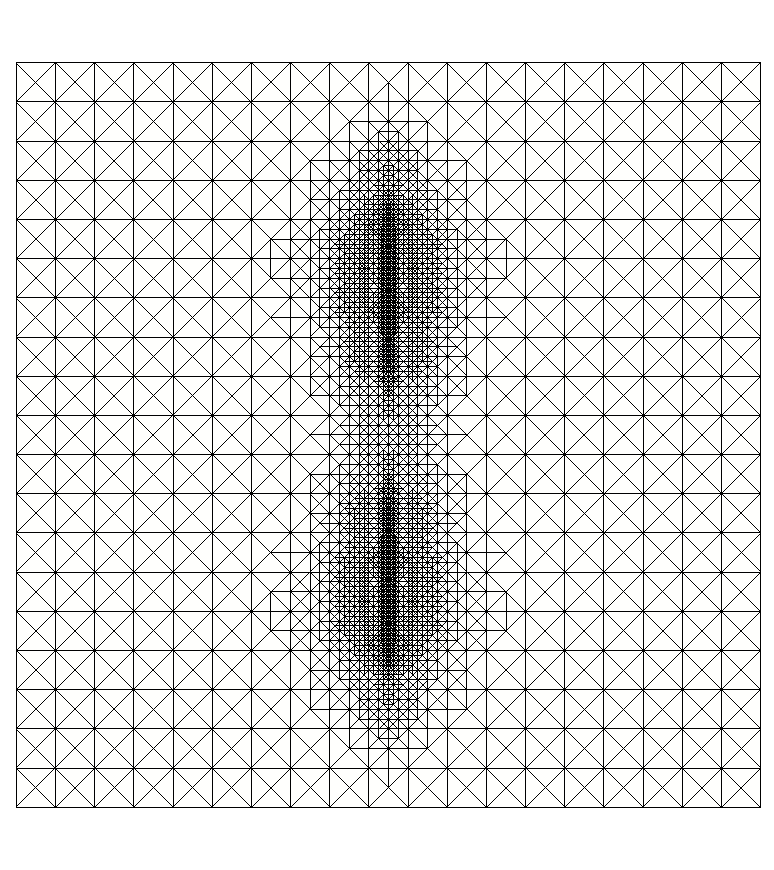} \hspace{0.3cm} \includegraphics[scale=0.12]{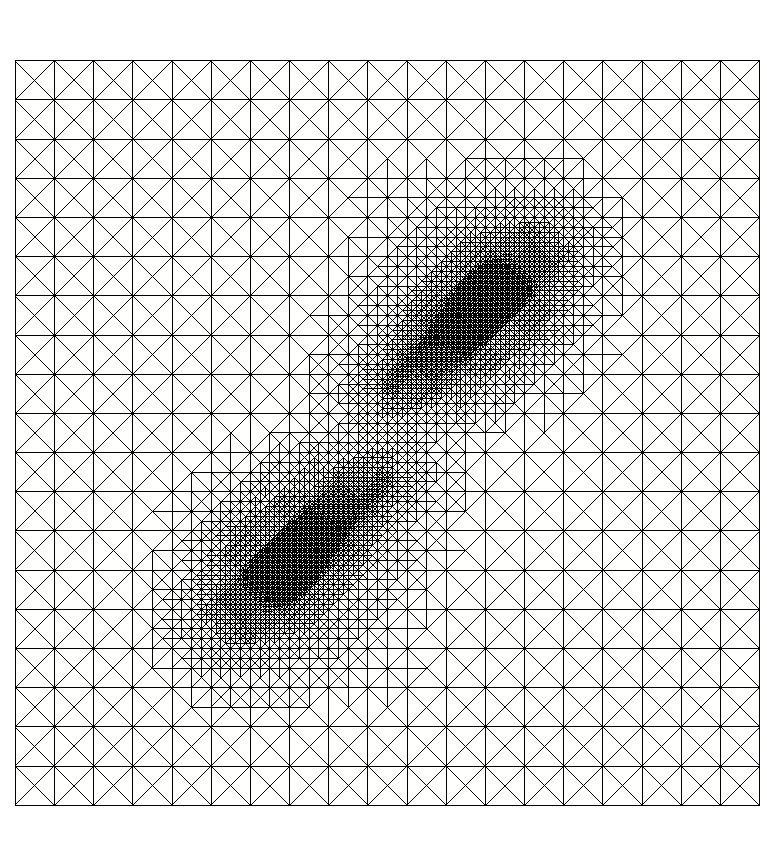} \hspace{0.3cm} \includegraphics[scale=0.12]{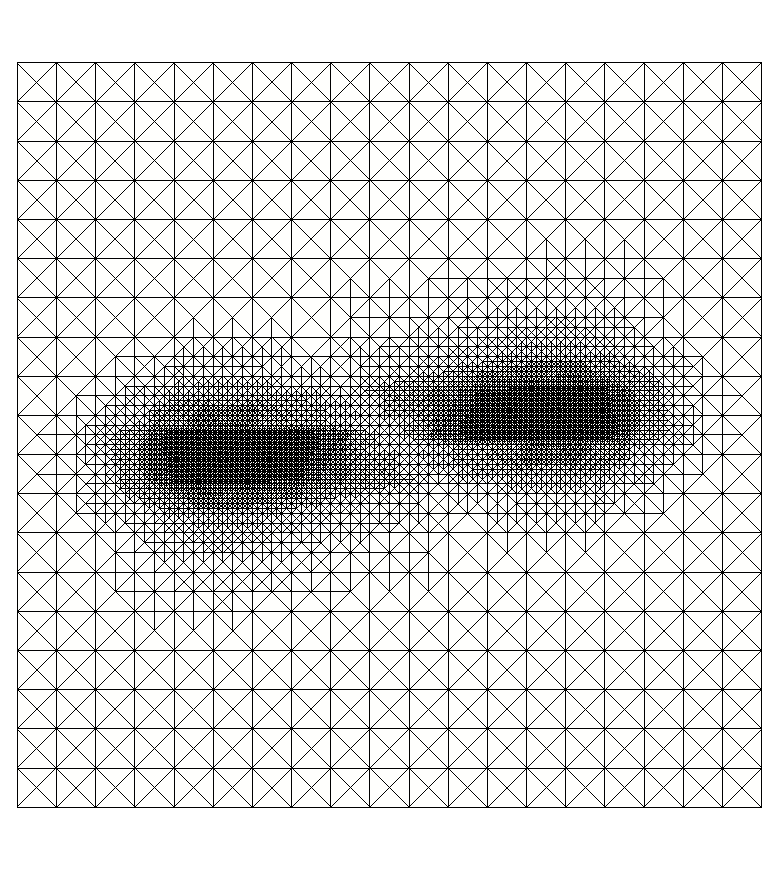}  \hspace{0.3cm} \includegraphics[scale=0.12]{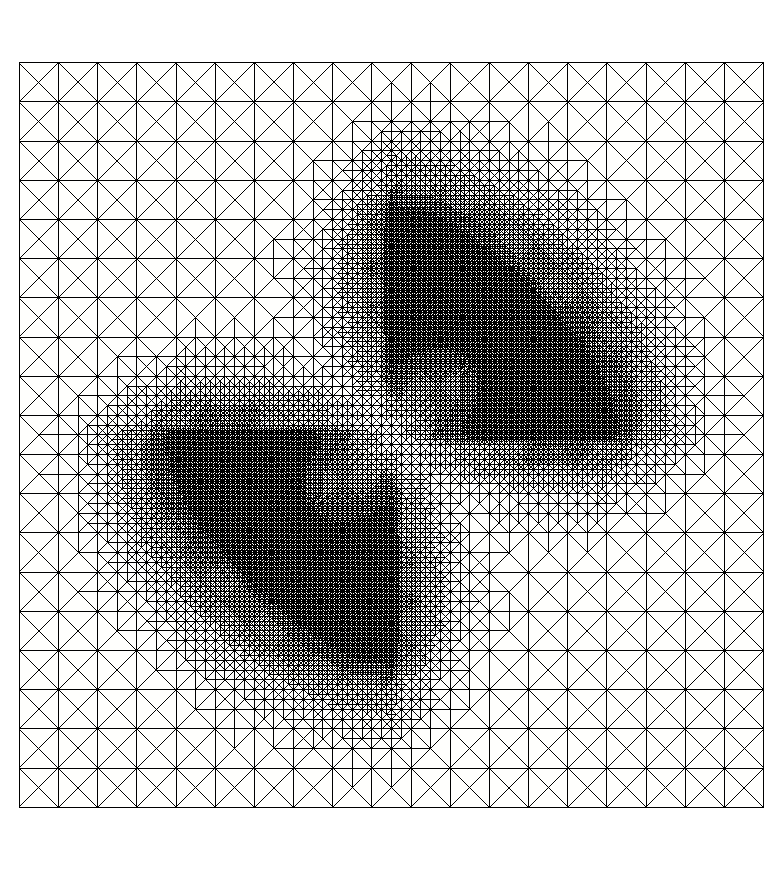} 
 \caption{\em Example 6.1: Adaptive finite element meshes at $t=t_0$ (left), $t=T/2$ (middle left), $t=T$ (middle right) and finest mesh (right)}
 \label{fig:Heat_meshes}
 \end{figure}

 \noindent Due to the steep gradients in the neighbourhood of the minimum and maximum, respectively, the use of an adaptive finite element discretization is 
 justified. The resulting computational meshes as well as the corresponding finest 
 mesh (reference mesh at the end of the simulation) are shown in Figure \ref{fig:Heat_meshes}. The number of node points of the adaptive meshes varies between 3637 and 7071 points. The finest mesh has 18628 node points. In contrary, a uniform mesh with the same discretization fineness as the finest triangle in the adaptive grids ($h_{\min}=0.0047$) would have 93025 node points. This clearly reveals the benefit of using adaptive meshes for snapshot generation. Particularly, the comparison of the computational times emphasizes the benefit of adaptive snapshot sampling: the snapshot generation on the adaptive mesh takes 944 seconds, whereas utilizing the uniform mesh it takes 8808 seconds. Therefore, we gain a speedup factor of 9 (see Table \ref{tab:CPU_times_heat}).

\begin{figure}[H]
 \centering
 \hspace{-0.5cm} \includegraphics[scale=0.32]{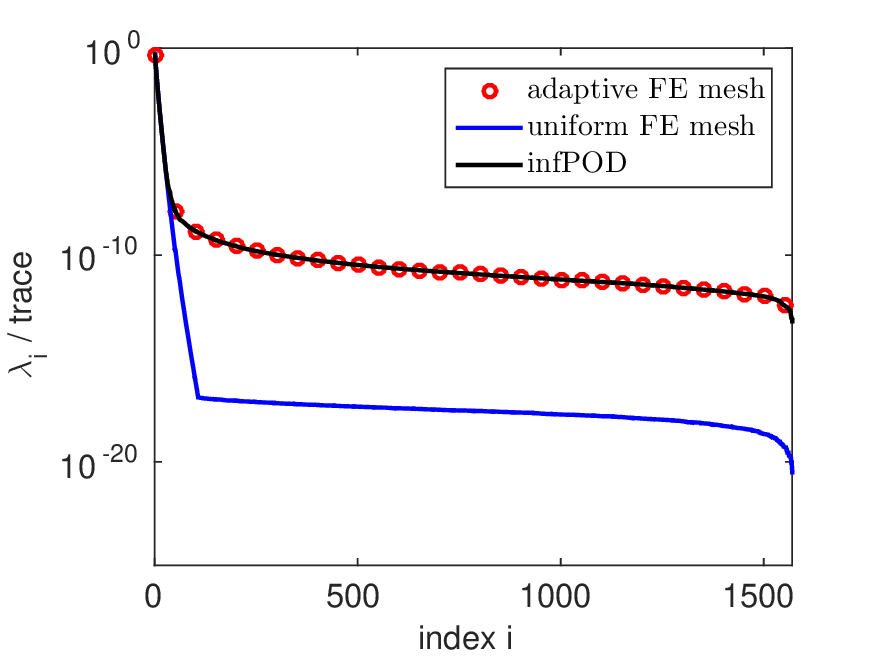}  \includegraphics[scale=0.32]{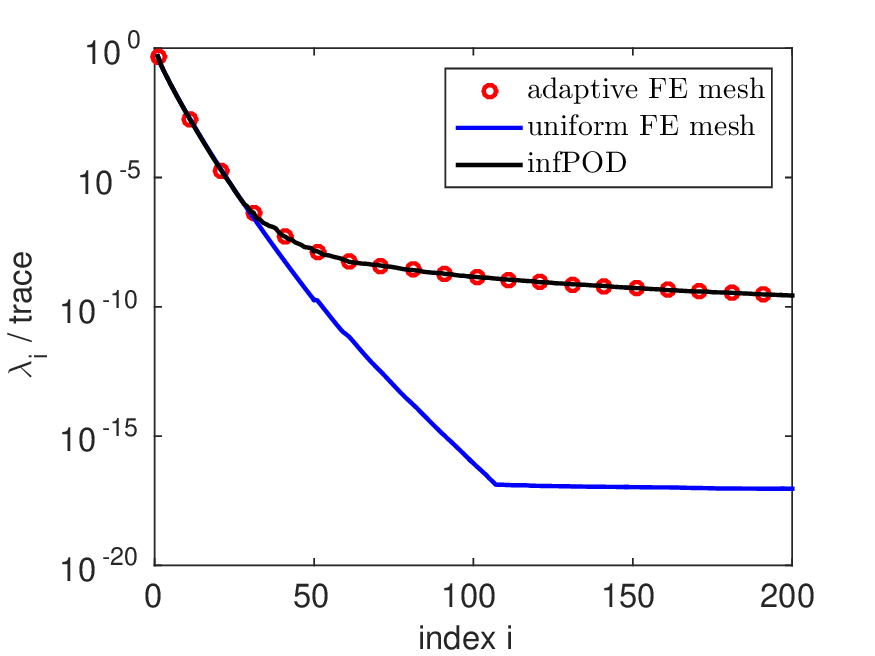} \includegraphics[scale=0.36]{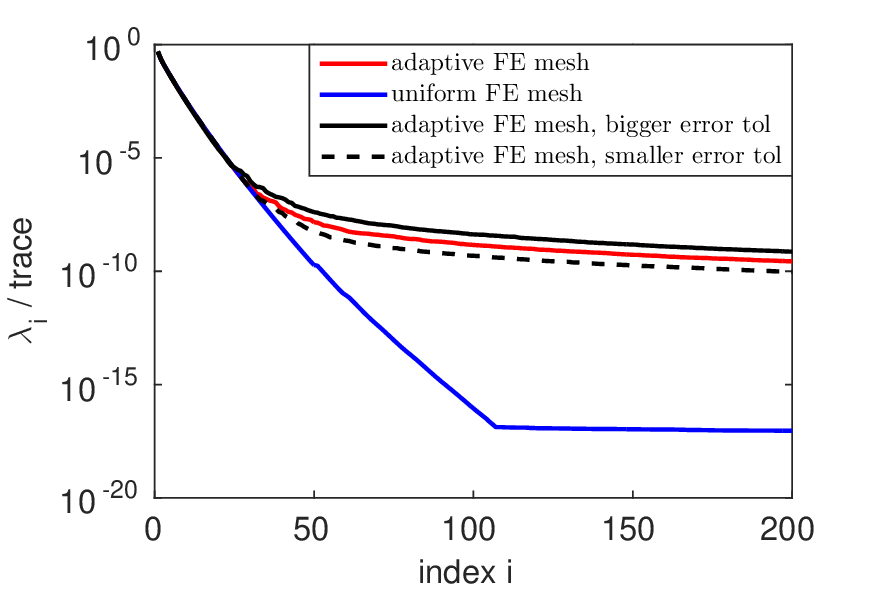}
 \caption{\em Example 6.1: Comparison of the normalized eigenvalues utilizing an adaptive and a uniform spatial mesh, respectively. Left: all eigenvalues, middle: first 200 largest eigenvalues, right: first 200 largest eigenvalues with different error tolerances for the adaptivity (1.5 times bigger and smaller error tolerances, respectively)}
 \label{fig:Heat_ev}
 \end{figure}

 \noindent In Figure \ref{fig:Heat_ev}, the resulting normalized eigenspectrum of the correlation matrix for uniform spatial discretization (``uniform FE mesh''), the normalized eigenspectrum of the matrix $\mathcal{K}$ \eqref{mathcalK} without 
 interpolation (``infPOD'') as well as the normalized eigenspectrum of the correlation matrix utilizing snapshots interpolated onto the finest mesh (``adaptive FE mesh'') is shown. We observe that the eigenvalues for both adaptive approaches coincide. This numerically validates what we expect from theory: the information content which is contained in the matrix $\mathcal{K}$ when we explicitly compute the entries without interpolation is the same as the information content contained within the eigenvalue problem which is formulated when using the finest mesh. No information is added or lost. Moreover, we recognize that about the first 28 eigenvalues computed corresponding to the adaptive simulation coincide with the simulation on a uniform mesh. From index 29 on, the methods deliver different results: for the uniform discretizations, the normalized eigenvalues fall below machine precision at around index 100 and stagnate. In contrary, the normalized eigenvalues for both adaptive approaches flatten in the order around $10^{-10}$. If the error tolerance for the spatial discretization error is chosen larger (or smaller), the stagnation of the eigenvalues in the adaptive method takes place at a higher (or lower) order (see Figure \ref{fig:Heat_ev}, right).
 
 Concerning dynamical systems, the magnitude of the eigenvalue corresponds to the 
 characteristic properties of the underlying dynamical system: the larger the 
 eigenvalue, the more information is contained in the corresponding eigenfunction. Since all adaptive meshes are contained in the uniform mesh, the difference in the amplitude of the eigenvalues is due to the interpolation errors during refinement and coarsening. This is the price we have to pay in order to get a fast snapshot generation utilizing adaptive finite elements. Moreover, the investigation of the decay of the eigenvalues can be interpreted as an analyzing tool for adaptivity 
 in the following sense: Using an adaptive mesh technique means that some parts of the domain are resolved coarsely according to the utilized error estimation, i.e. information gets lost. In the sense of a singular value analysis, this can be explained that adaptivity neglects the noise which is indicated by the singular values on the uniform spatial mesh at those places which are not resolved with the adaptive grid. We conclude that the overtones which get lost in the adaptive computations lie in the same space which is not considered by POD when using the adaptive finite element snapshots. This allows us to characterize the space which is not resolved by adaptivity. From this point of view, adaptivity can be interpreted as a smoother.

\begin{figure}[H]
 \centering
   \includegraphics[scale=0.16]{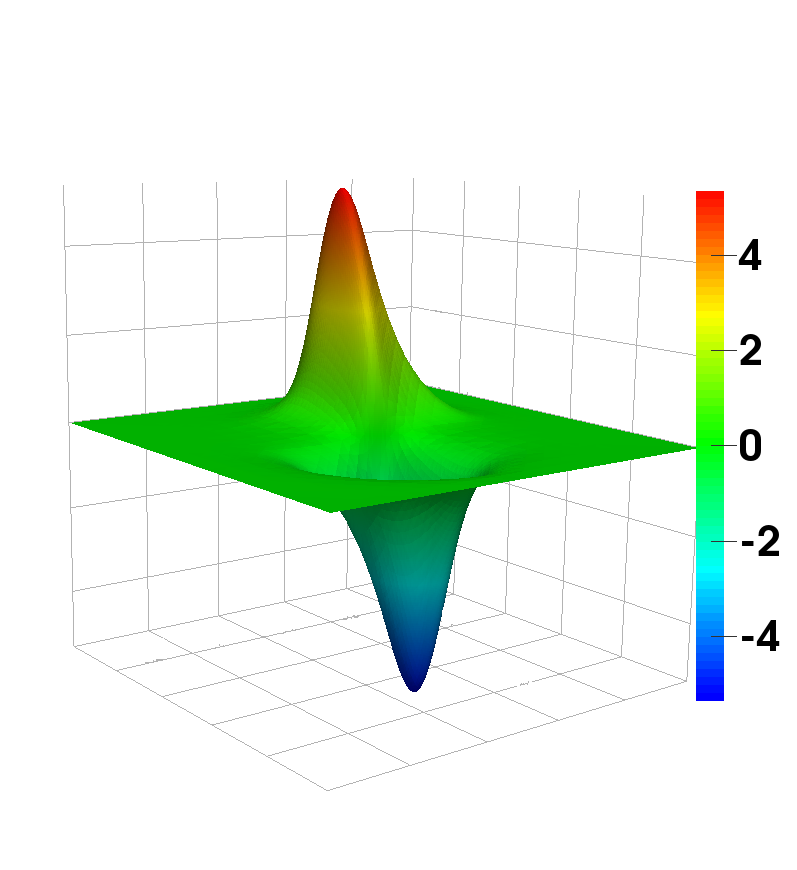} \hspace{0.2cm} \  \includegraphics[scale=0.16]{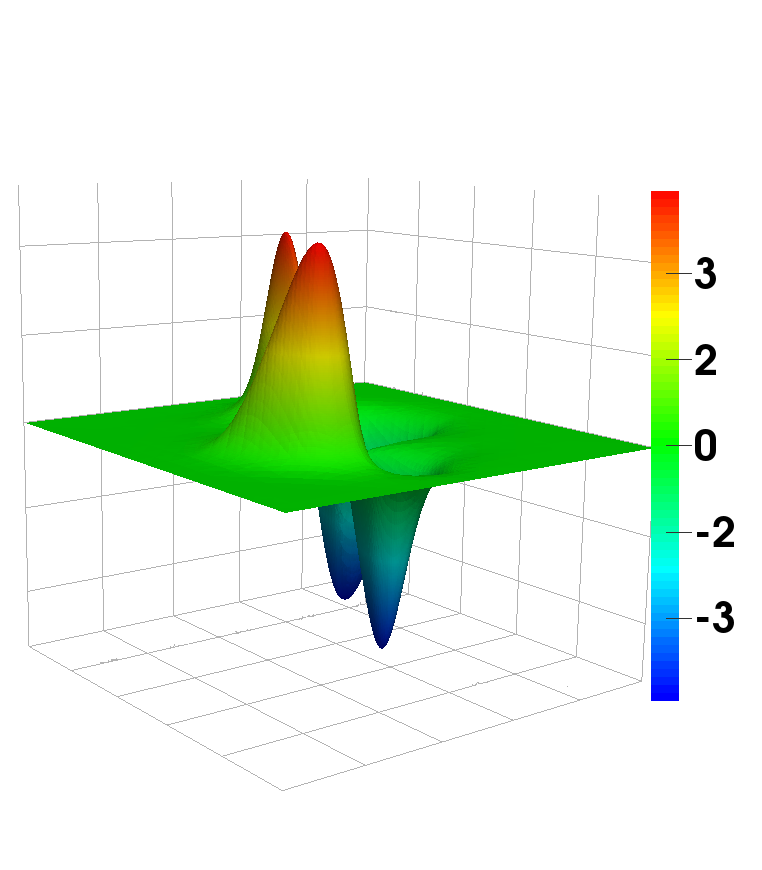} \hspace{0.2cm} \includegraphics[scale=0.16]{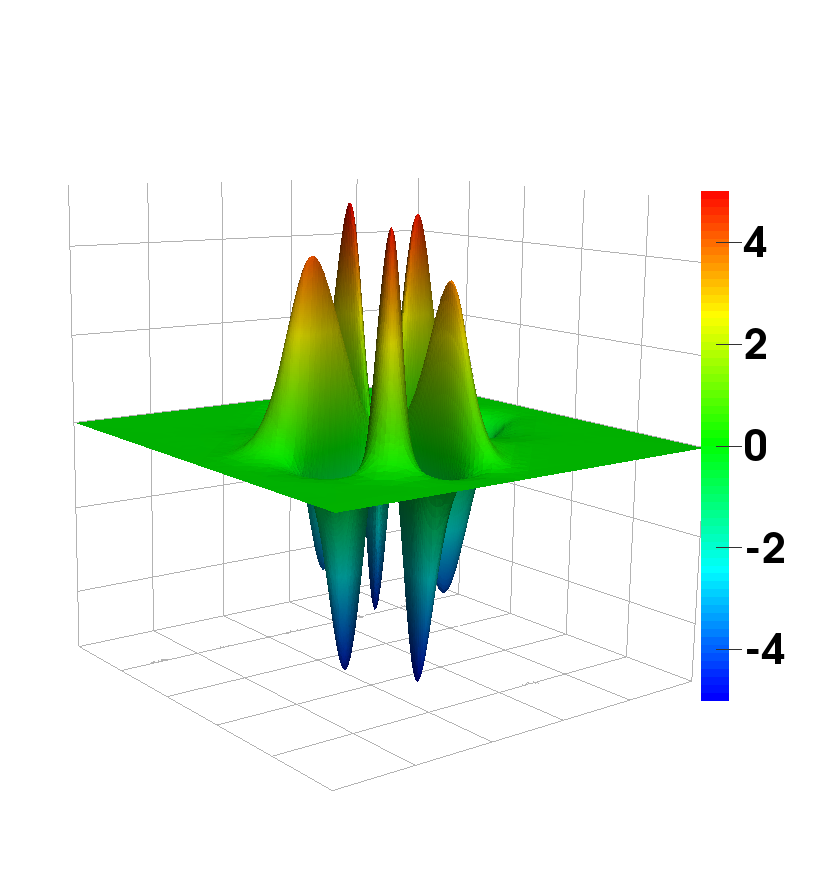} \\ 
 \hspace{-0.5cm} \includegraphics[scale=0.12]{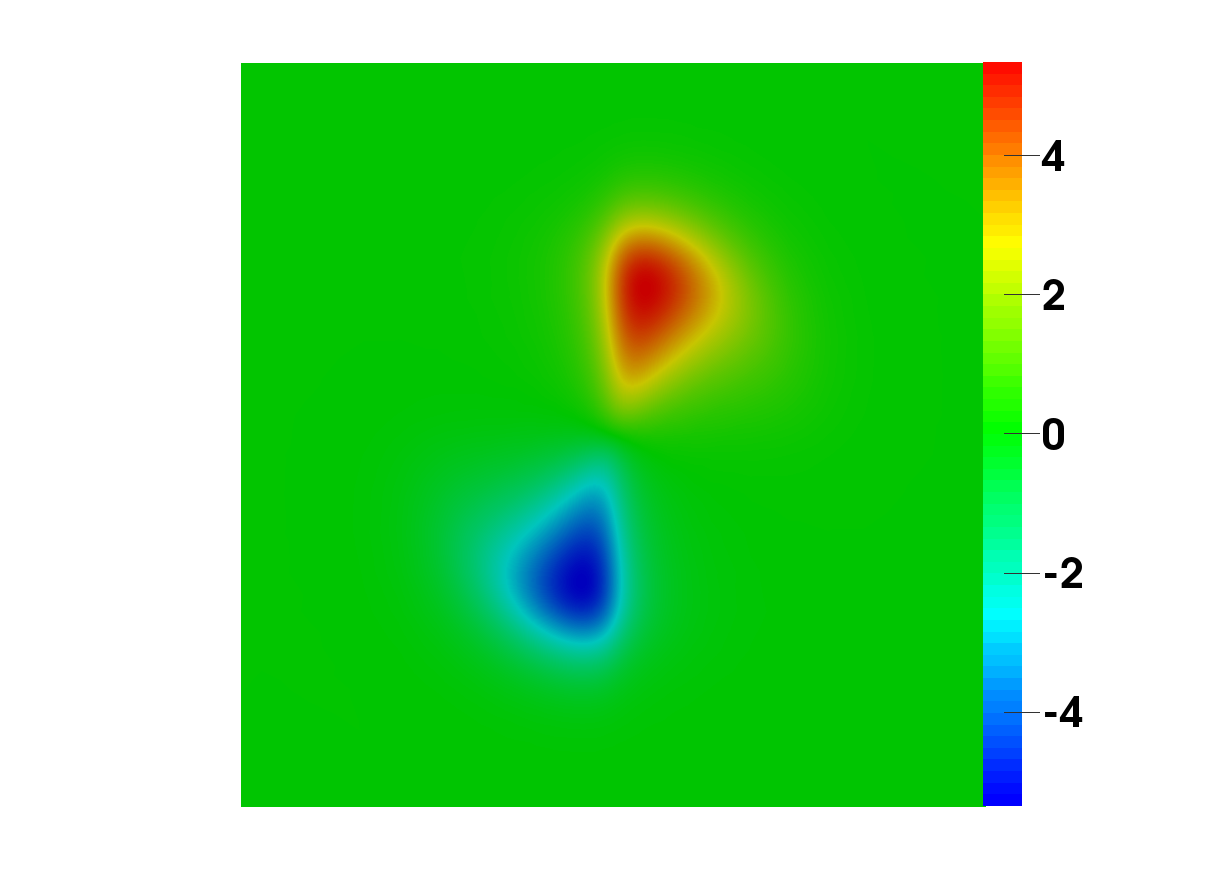} \hspace{-0.5cm} \includegraphics[scale=0.12]{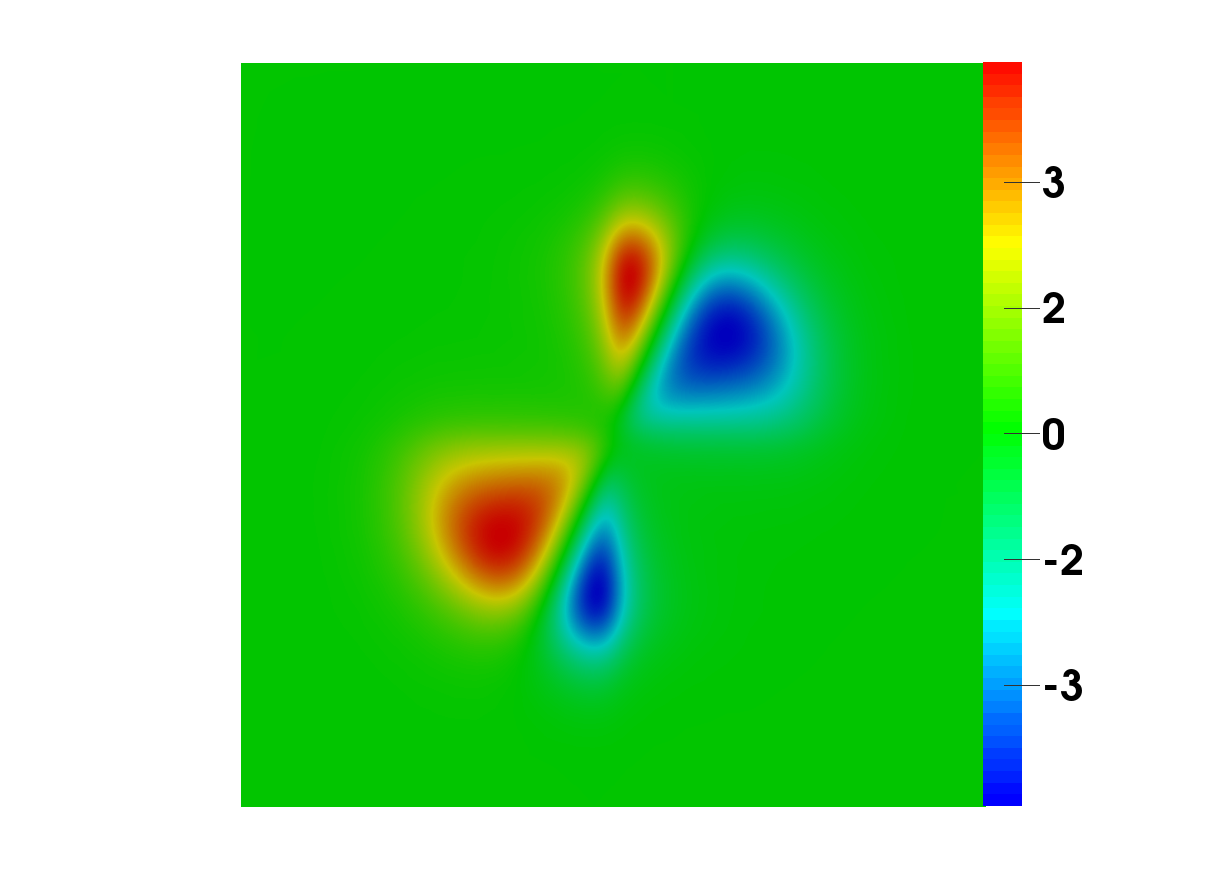} \hspace{-0.5cm} \includegraphics[scale=0.12]{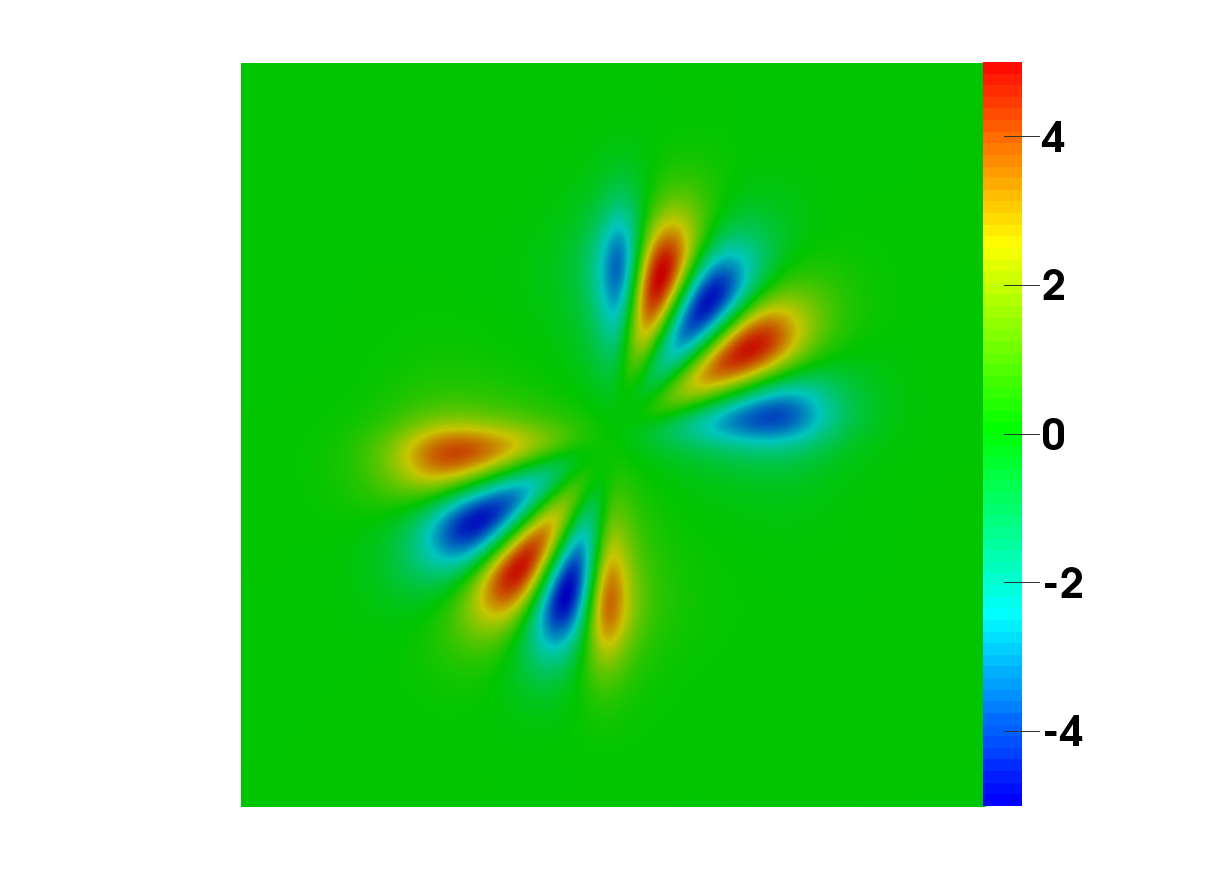} 
 \caption{\em Example 6.1: Surface plot (top) and view from above (bottom) of the POD basis functions
 $\psi_1$ (left), $\psi_2$ (middle) and $\psi_5$ (right)}
 \label{fig:Heat_PODmodes}
 \end{figure}

 \noindent Since the first few POD basis functions are the most important ones 
 regarding the captured information, we visualize $\psi_1, \psi_2$ and $\psi_5$ in 
 Figure \ref{fig:Heat_PODmodes}, which are computed corresponding to using an adaptive grid. The POD basis functions corresponding to the uniform spatial discretization have a similar appearance. Note that the POD bases are unique up to the sign. We can recognize the initial condition in the first POD basis function. Then, the index of the POD basis corresponds to the number of maxima and minima of the POD basis: $\psi_2$ has two minima and two maxima etc. This behaviour is similar to the increasing oscillations in higher frequencies in trigonometric approximations. The increasing number of oscillations is necessary in order to approximate the transport of the steep gradients of the solution with increasing accuracy.

\begin{figure}[htbp]
 \centering
  \includegraphics[scale=0.12]{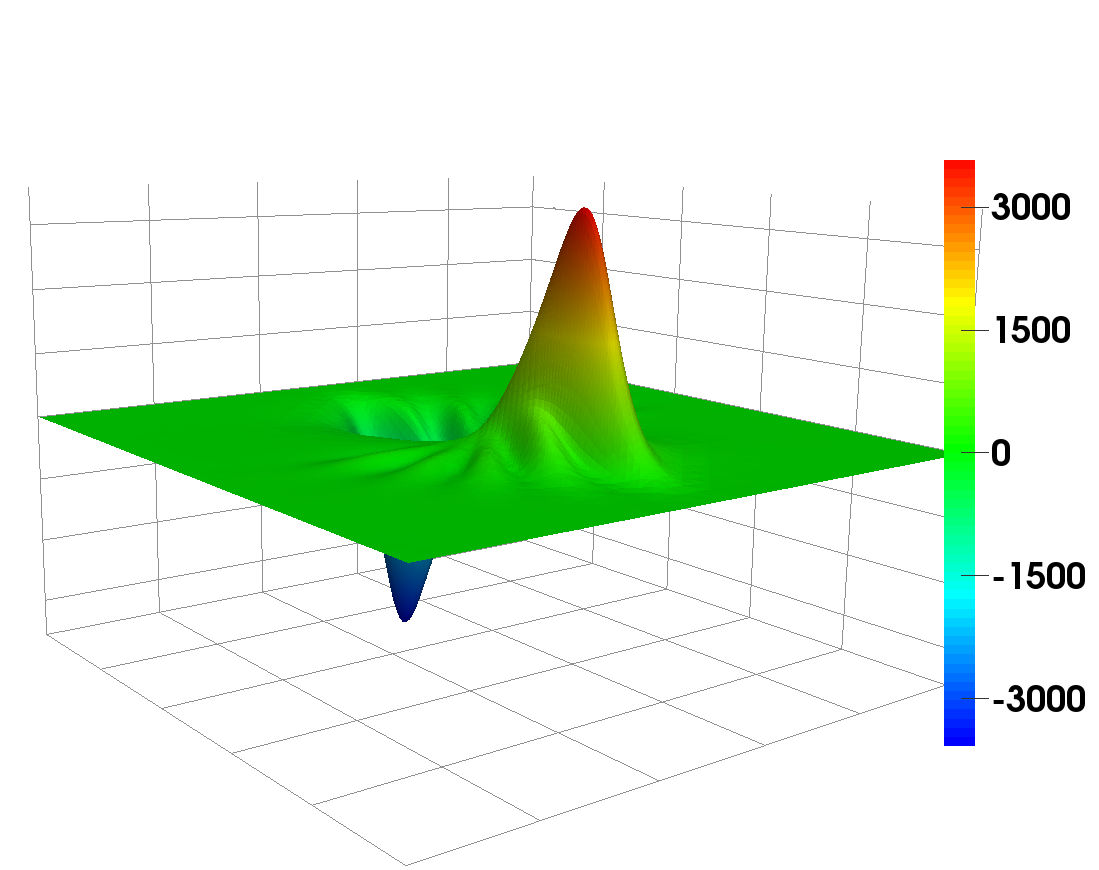} \hspace{0.1cm} \includegraphics[scale=0.12]{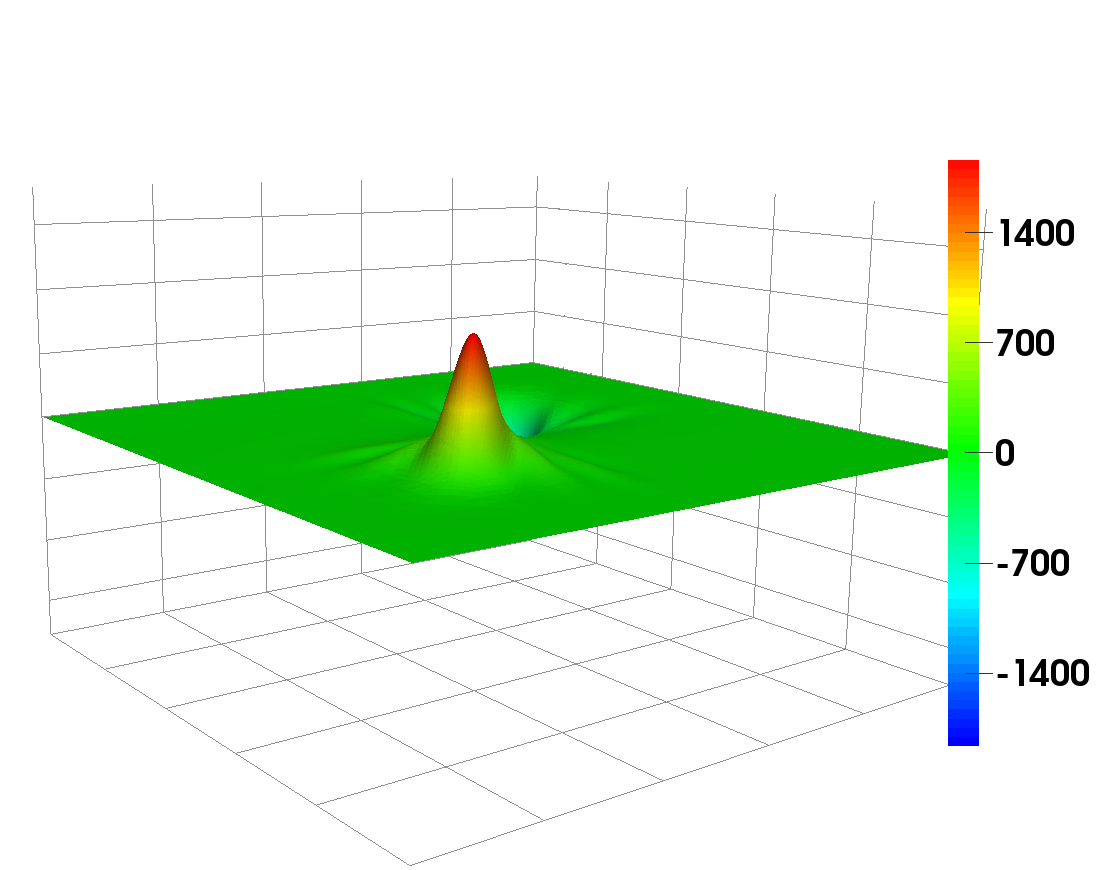} \hspace{0.1cm} \includegraphics[scale=0.12]{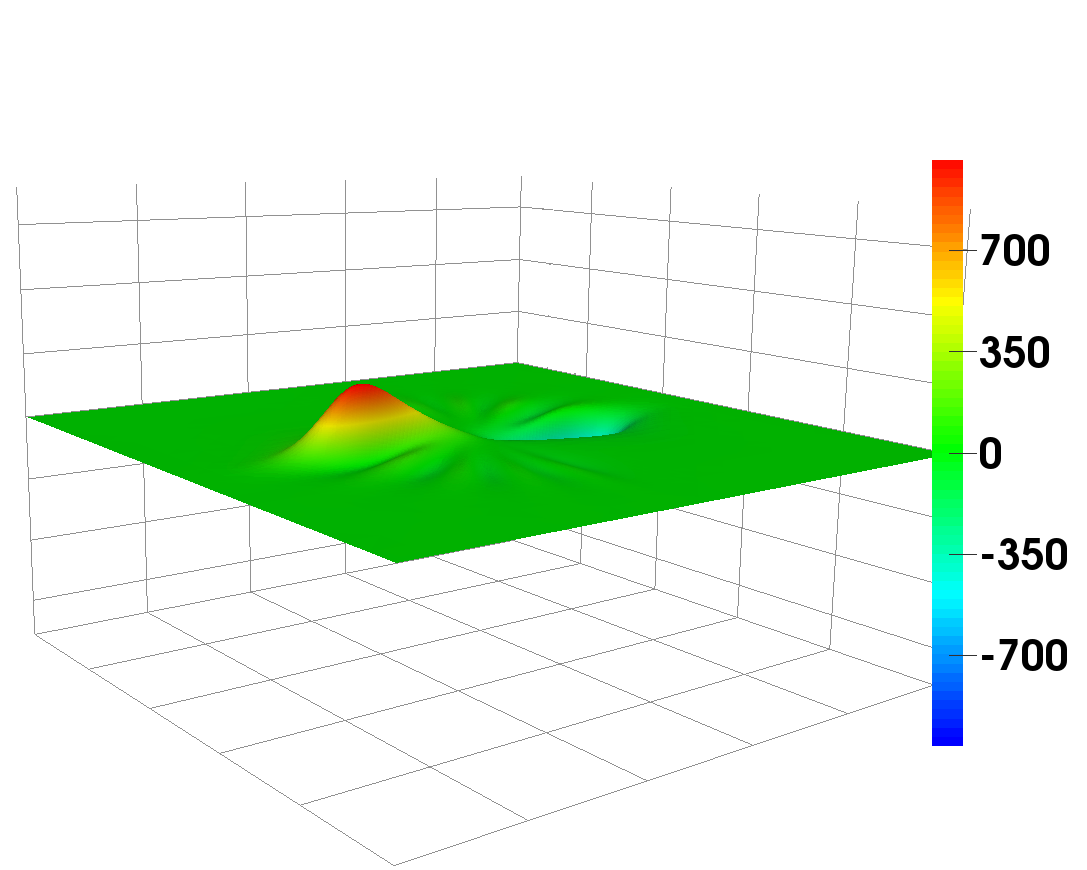}\\
 \includegraphics[scale=0.12]{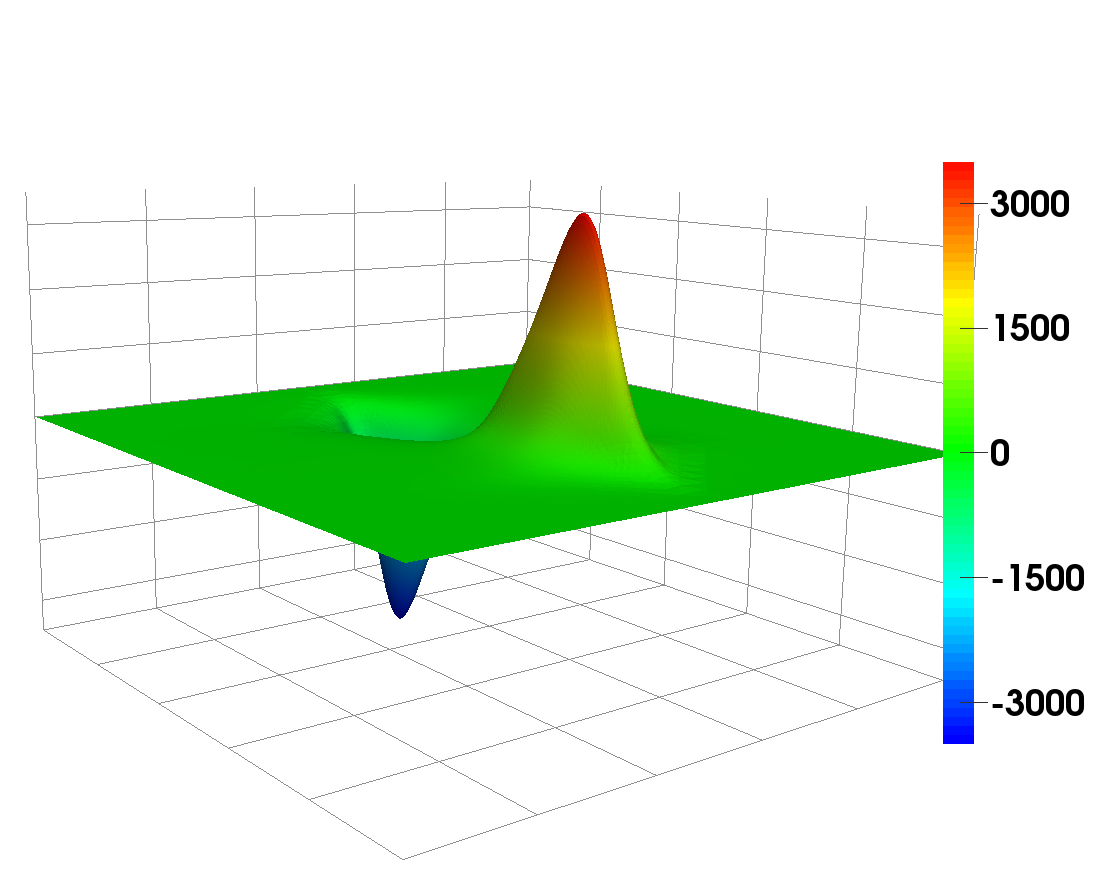} \hspace{0.1cm} \includegraphics[scale=0.12]{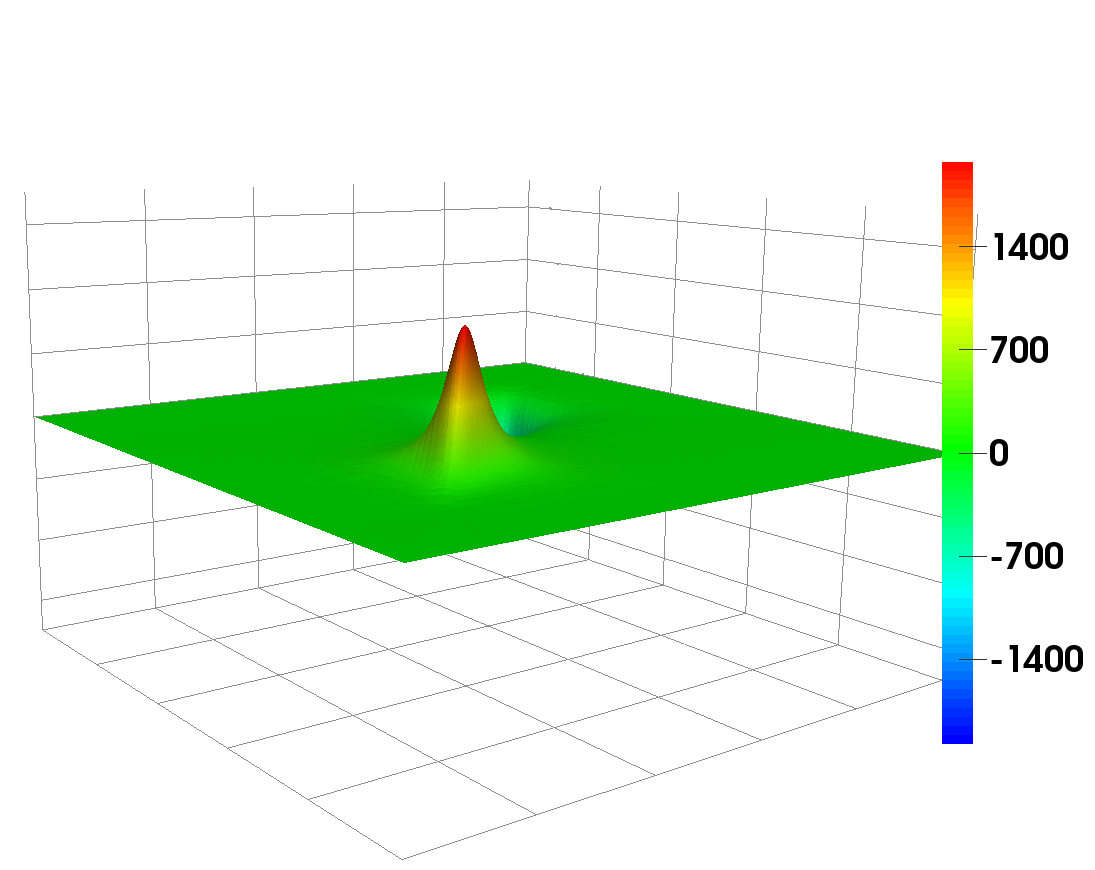}  \hspace{0.1cm} \includegraphics[scale=0.12]{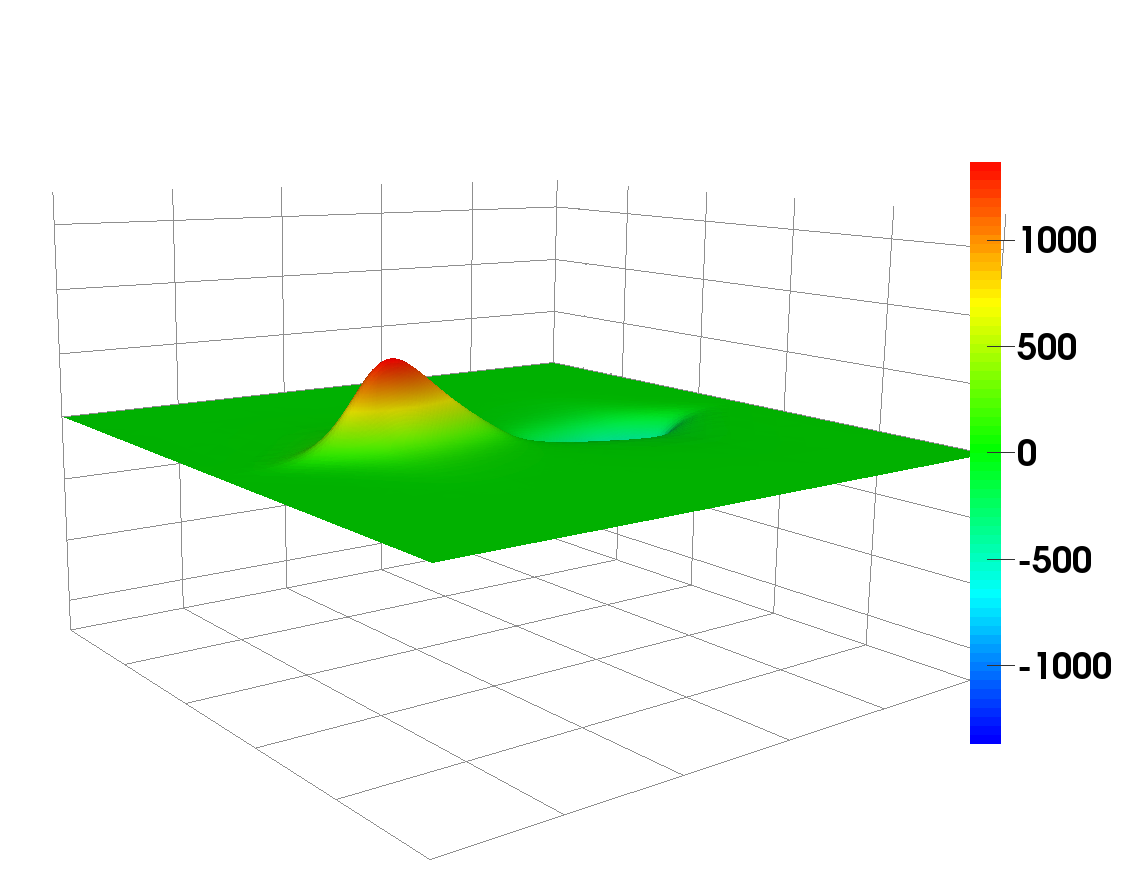}    
 \caption{\em Example 6.1: Surface plot of the POD solution utilizing $\ell = 10$ (top) and $\ell = 50$ (bottom) POD basis functions at $t=t_0$ (left), $t=T/2$ (middle) and $t=T$ (right)}
 \label{fig:Heat_POD}
 \end{figure}

 \noindent The POD solutions for $\ell = 10$ and $\ell = 50$ POD basis 
 functions utilizing spatial adaptive snapshots which are interpolated onto the finest mesh are shown in Figure \ref{fig:Heat_POD}. The visual comparison makes clear what influence the increase of the number of utilized POD basis functions has on the approximation quality. The more POD basis functions we use (until stagnation of the corresponding eigenvalues), the less oscillations appear in the POD solution and the better is the approximation.

Table \ref{tab:errors_heat} compares the approximation quality of the POD solution utilizing adaptively generated snapshots which are interpolated onto the finest mesh with snapshots of uniform spatial discretization depending on different POD basis lengths. Our approach from Sections 3 and 4 delivers very similar results as the use of adaptive finite element snapshots which are interpolated onto the finest mesh. For example, for $\ell = 20$ POD bases, we get the following errors: relative $L^2(0,T;L^2(\Omega))$-error between the POD solution and the finite element solution: $\varepsilon_{\text{FE}} = 3.07 \cdot 10^{-2}$, relative $L^2(0,T;L^2(\Omega))$-error between the POD solution and the true solution: $\varepsilon_{\text{true}} = 2.16 \cdot 10^{-2}$.\\

\begin{table}[H]
    \centering
    
  \begin{tabular}{ c   | c | c | c | c }
 $\ell$  &  $\varepsilon_{\text{FE}}^{\text{ad}}$ & $\varepsilon_{\text{FE}}^{\text{uni}}$   & $\varepsilon_{\text{true}}^{\text{ad}}$ &   $\varepsilon_{\text{true}}^{\text{uni}}$\\
 \hline
 1  & \hspace{-0.2cm}$1.30 \cdot 10^{0}$ & \hspace{-0.3cm} $1.30 \cdot 10^{0}$ & \hspace{-0.3cm}  $1.28 \cdot 10^{0}$  & \hspace{-0.3cm}  $1.30 \cdot 10^{0}$ \\
 3  & $7.49 \cdot 10^{-1}$  & $7.58 \cdot 10^{-1}$ &  $7.46 \cdot 10^{-1}$  &  $7.60 \cdot 10^{-1}$\\
 5   & $4.39 \cdot 10^{-1}$  & $4.45 \cdot 10^{-1}$ &  $4.39 \cdot 10^{-1}$  & $4.46 \cdot 10^{-1}$ \\
 10   & $1.37 \cdot 10^{-1}$   & $1.37 \cdot 10^{-1}$  & $1.36 \cdot 10^{-1}$  &  $1.38 \cdot 10^{-1}$ \\
 20  & $3.08 \cdot 10^{-2}$  & $1.56 \cdot 10^{-2}$  & $2.17 \cdot 10^{-2}$  & $1.60 \cdot 10^{-2}$   \\
 30  & $2.59 \cdot 10^{-2}$  &  $2.04 \cdot 10^{-3}$ & $1.49   \cdot 10^{-2}$  &  $3.00 \cdot 10^{-3}$ \\
 50  & $2.63 \cdot 10^{-2}$   &  $5.67 \cdot 10^{-5}$ & $1.41 \cdot 10^{-2}$  &    $ 2.07 \cdot 10^{-3}$  \\
 100  & $2.61 \cdot 10^{-2}$  &  $6.48 \cdot 10^{-8}$  & $1.40 \cdot 10^{-2}$  & $2.06 \cdot 10^{-3}$\\
 150  & $2.61 \cdot 10^{-2}$  &  $8.13 \cdot 10^{-7}$  & $1.39 \cdot 10^{-2}$   &   $2.07 \cdot 10^{-3}$\\
 \end{tabular}

  \caption{\em Example 6.1: Relative $L^2(0,T;L^2(\Omega))$-error between the POD solution and the finite element solution (columns 2-3) and the true solution (columns 4-5), respectively, utilizing adaptive finite element snapshots which are interpolated onto the finest mesh and utilizing a uniform mesh}
 \label{tab:errors_heat}
\end{table}

 \noindent We note that the error between POD solution and finite element solution utilizing a uniform mesh decays down to the order $10^{-8}$ ($\ell = 100$) and 
 then stagnates. This behaviour is clear, since the more POD basis we include (up 
 to stagnation of the corresponding eigenvalues), the better is the POD solution an 
 approximation for the finite element solution. In contrary, the error between the POD solution and the true solution starts to stagnate from $\ell = 29$. This is due to the fact that at this point the spatial discretization error dominates the modal error. This is in accordance to the decay of the eigenvalues shown in Figure \ref{fig:Heat_ev}. Due to the error estimation \eqref{errorestimation}, the error between the true solution to \eqref{P} and the POD reduced order solution is not only bounded by the sum of the neglected eigenvalues, which is small for sufficiently large number of utilized POD modes. It is also restricted by the spatial and temporal discretization error, which leads to a stagnation of the error in Table \ref{tab:errors_heat}, columns 4 and 5.\\

\noindent Finally, of particular interest is the computational efficiency of the 
POD reduced order modeling utilizing adaptive finite element discretizations. For this reason, the computational times for the full and the low order simulation utilizing uniform finite element discretizations and adaptive finite element snapshots, which are interpolated onto the finest mesh, respectively, are listed in Table \ref{tab:CPU_times_heat}.\\

\begin{table}[htbp]
\centering
 \begin{tabular}{ l | c | c | c}
  & adaptive FE mesh  & uniform FE mesh & speedup factor \\
 \hline
 FE simulation & 944 sec & 8808 sec & 9.3 \\
 POD offline computations & 264 sec & 1300 sec & 4.9 \\
 POD simulation & \multicolumn{2}{c|}{ 187 sec}  & --\\ 
 \hline
 speedup factor & 5.0 &  47.1 & --\\
 \end{tabular}
 \vspace{0.4cm} \caption{\em Example 6.1: CPU times for FE and POD simulation utilizing uniform finite element meshes and adaptive finite element snapshots which are interpolated onto the finest mesh, respectively, and utilizing $\ell = 50$ 
 POD modes}
 \label{tab:CPU_times_heat}
  \end{table}

\noindent Once the POD basis is computed in the offline phase, the POD simulation corresponding to adaptive snapshots is 5 times faster than the FE simulation utilizing adaptive finite element meshes. This speedup factor even gains greater importance, if we think of optimal control problems, where the repeated solving of several partial differential equations is necessary.  In the POD offline phase, the most expensive task is to express the snapshots with respect to the common finite element space, which takes 226 seconds. Since $\mathcal{K}$ is symmetric, it suffices to calculate the entries on and above the diagonal, which are $\sum_{k=1}^{n+1} k = \frac{1}{2}((n+1)^2 + n+1)$ entries. Thus, the computation of each entry in the correlation matrix $\mathcal{K}$ using a common finite element space takes around 0.00018 seconds. We note that in the approach explained in Section 3 and 4, the computation of the matrix $\mathcal{K}$ \eqref{mathcalK} is expensive. For each entry the calculation time is around 0.03 seconds, which leads to a computation time of around 36997 seconds for the matrix $\mathcal{K}$. The same effort is needed to build $\mathcal{Y}^\star \mathcal{A} \mathcal{Y}$. In this case, the offline phase takes therefore around 88271 seconds. For this reason, the approach to interpolate the adaptive generated snapshots onto the finest mesh is computationally more favorable. But since the computation of $\mathcal{K}$ and $\mathcal{Y}^\star \mathcal{A} \mathcal{Y}$ can be parallelized, the offline computation time can be reduced provided that the appropriate hardware is available.\\

 \noindent \textbf{Example 6.2: Cahn-Hilliard system.} We consider Example 2.4, 
 \eqref{CHcoupled} of the Cahn-Hilliard equations given in the coupled formulation for the phase field $c$ and the chemical potential $w$. The data is chosen as follows: we consider the rectangular domain $\Omega = (0,1.5) \times (0, 0.75)$, the end time $T = 0.025$, constant mobility $m \equiv 0.00002$ and a constant surface tension $\sigma \equiv 24.5$. The interface parameter $\varepsilon $ is set to $\varepsilon = 0.02$, which leads to an interface thickness of about $\pi \cdot \varepsilon \approx 0.0628$. We utilize the relaxed double obstacle free energy $W_s^{\text{rel}}$, \eqref{rdofe} with $s=10^4$. As initial condition, we choose a circle with radius $r = 0.25$ and center $(0.375,0.375)$. The initial condition 
 is transported horizontally with constant velocity $v = (30,0)^T$. Let us define the uniform time discretization 
 $$ t_j = j \Delta t $$
 for $j = 0, \dotsc, 1000$ with $\Delta t = 2.5 \cdot 10^{-5}$. We utilize a semi-implicit Euler scheme for temporal discretization. Let $c^{j-1} \in V$ and $c^j \in V$ denote the time-discrete solution at $t_{j-1}$ and $t_j$. Based on the variational formulation \eqref{CH-weak-abstract} we tackle the time-discrete version of \eqref{CHcoupled} in the form: given $c^{j-1}$, find $c^j$ with associated $w^j$ solving
 \begin{equation}\label{CH_num}
  \left\{
\begin{array}{rcll}
\langle \displaystyle\frac{c^j - c^{j-1}}{\Delta t}, v_1 \rangle_{L^2} + \langle v \cdot \nabla c^{j-1},v_1 \rangle_{L^2} + m\langle \nabla w^j, \nabla v_1 \rangle_{L^2}& = & 0 & \forall v_1 \in V,\\
- \langle w^j, v_2 \rangle_{L^2} + \sigma\varepsilon \langle \nabla c^j, \nabla v_2 \rangle_{L^2} +  \displaystyle\frac{\sigma}{\varepsilon} \langle W'_+(c^j) + W'_-(c^{j-1}) , v_2 \rangle_{L^2}  & = &  0 & \forall v_2 \in V,\\
\end{array}
\right.
\end{equation}
and $c^0 = c_0$. According to \eqref{CH-weak-abstract}, here it is $V=\{v \in H^1(\Omega), \frac{1}{|\Omega|} \int_\Omega v dx = 0\}$. Note that the free energy function $W$ is split into a convex part $W_+$ and a concave part $W_-$, such that $W = W_+ + W_-$ and $W'_+$ is treated implicitly with respect to time and $W'_-$ is treated explicitly with respect to time. This leads to an unconditionally energy stable time marching scheme, compare \cite{Eyr98}. The system \eqref{CH_num} is discretized in space utilizing piecewise linear and continuous finite elements and solved using a semi-smooth Newton method.\\

\begin{figure}[htbp]
 \centering
 \includegraphics[scale=0.15]{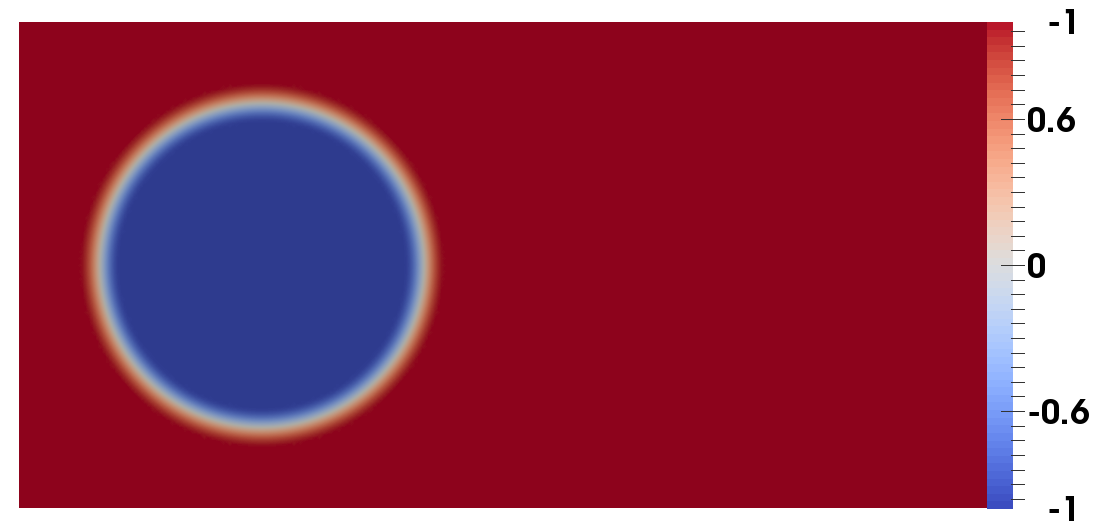} \hspace{0.2cm} \includegraphics[scale=0.15]{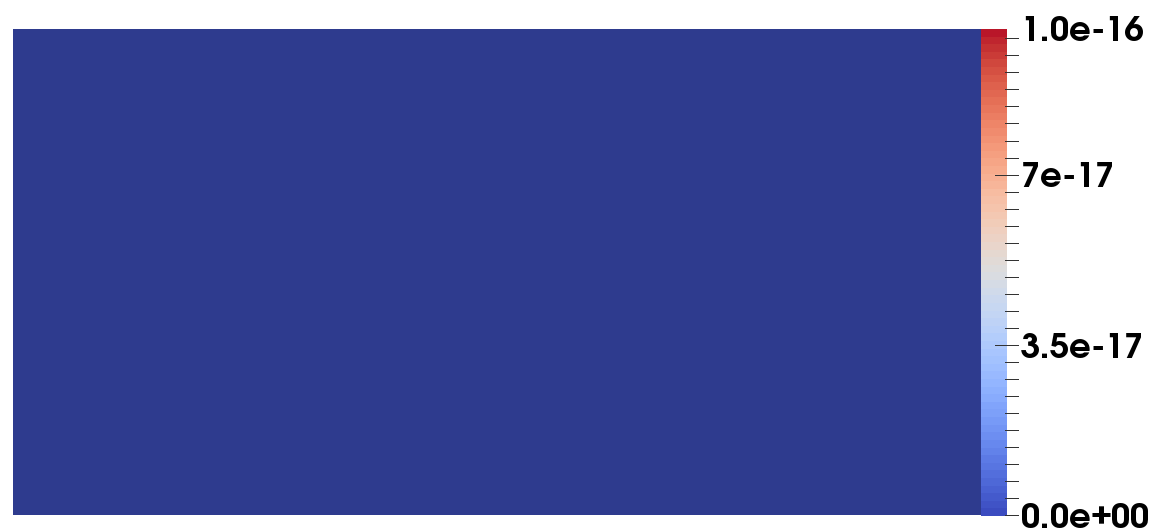}\\[0.1cm]
  \includegraphics[scale=0.15]{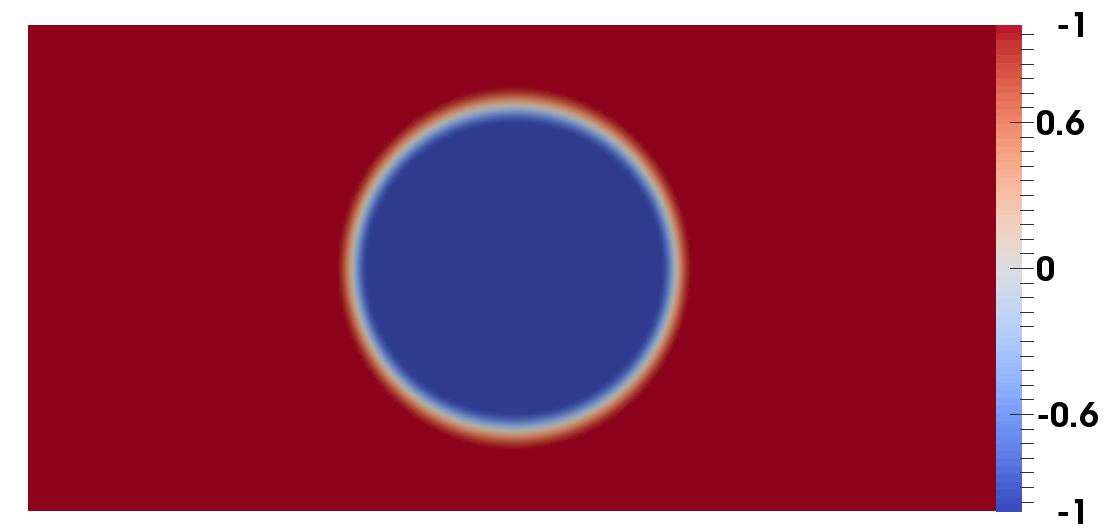} \hspace{0.2cm} \includegraphics[scale=0.15]{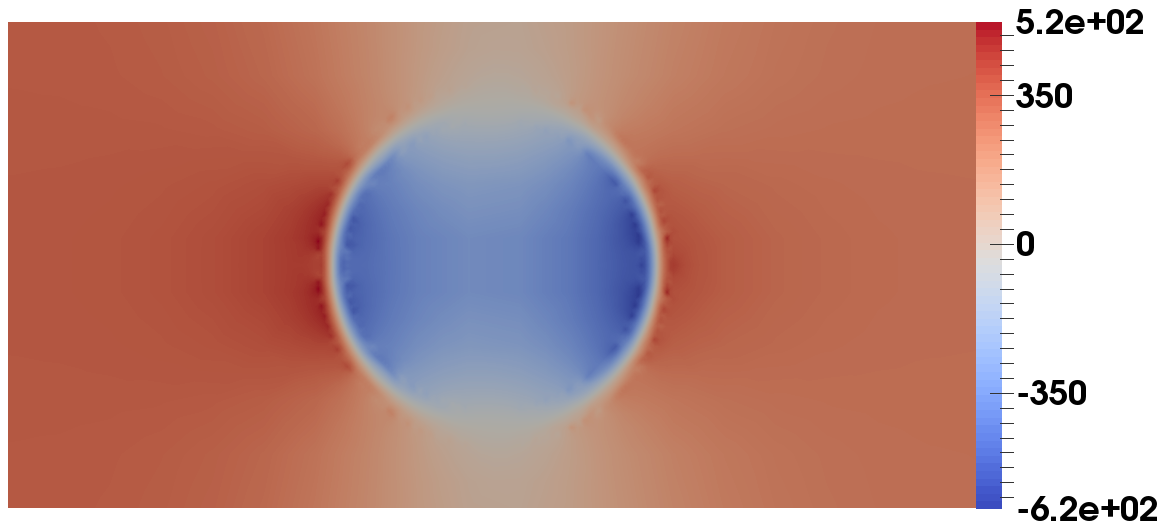}\\[0.1cm]
 \includegraphics[scale=0.15]{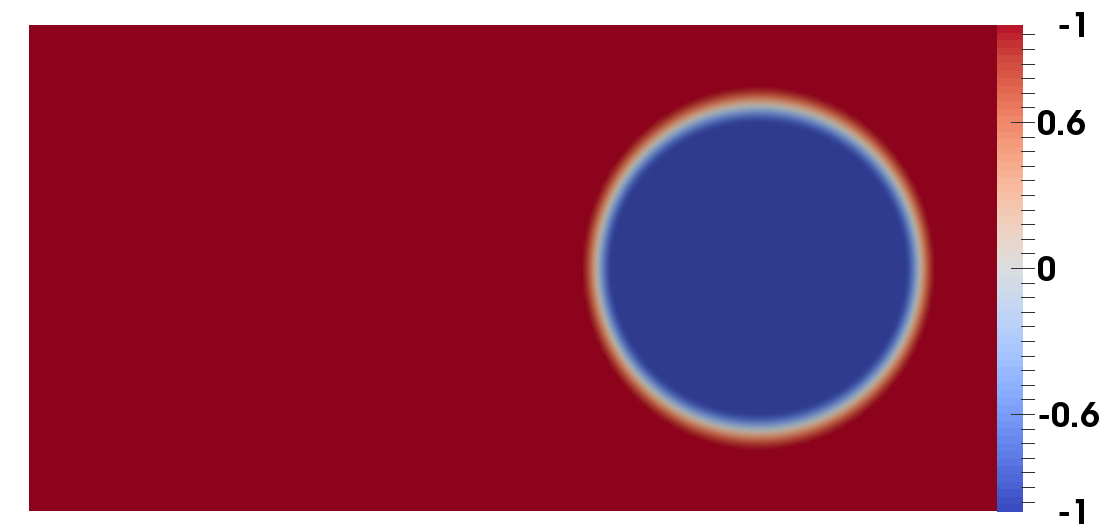} \hspace{0.2cm} \includegraphics[scale=0.15]{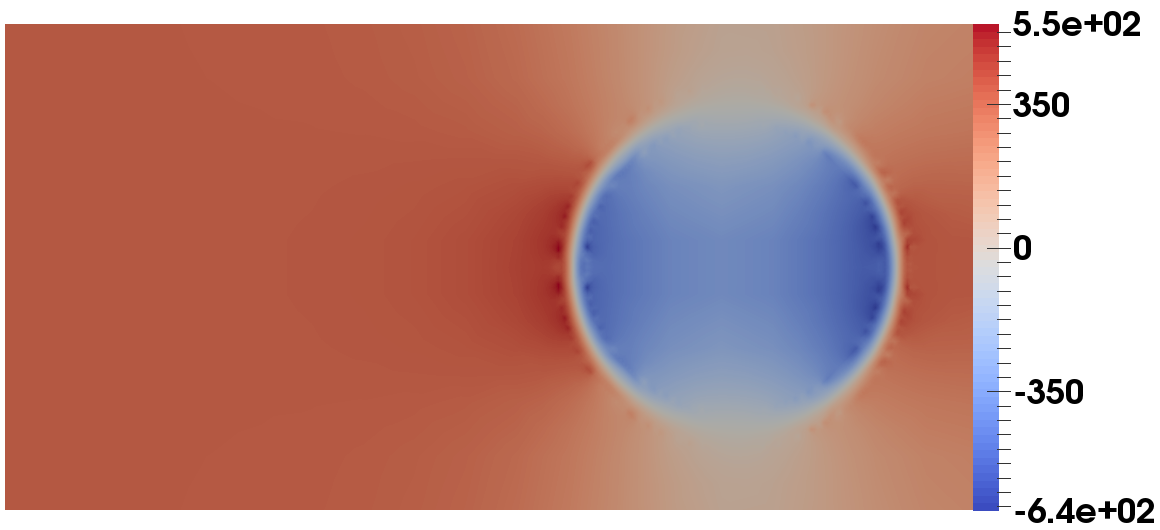}
 \caption{\em Example 6.2: Phase field $c$ (left) and chemical potential $w$ (right) computed on adaptive finite element meshes at $t=t_0$ (top), $t=T/2$ (middle) and $t=T$ (bottom)}
 \label{fig:CH_adapt_sim}
 \end{figure}

 \noindent Figures \ref{fig:CH_adapt_sim} shows the phase field (left) and the chemical potential (right) for the finite element simulation utilizing adaptive meshes. The initial condition $c_0$ is transported horizontally with constant velocity.\\

\begin{figure}[H]
 \centering
 \includegraphics[scale=0.15]{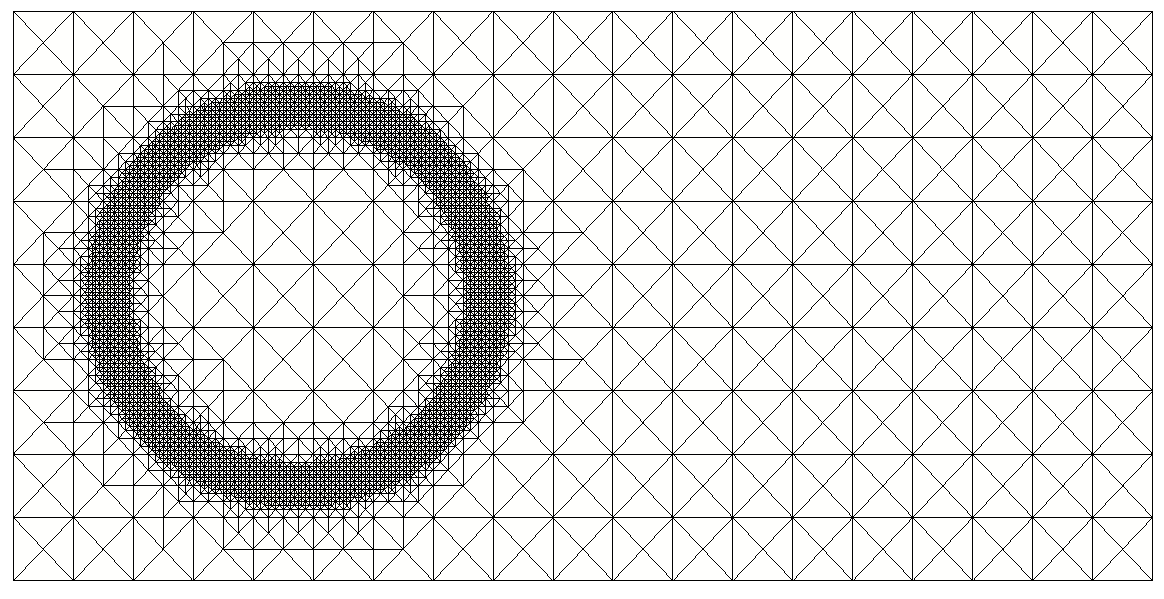} \hspace{0.2cm} \includegraphics[scale=0.15]{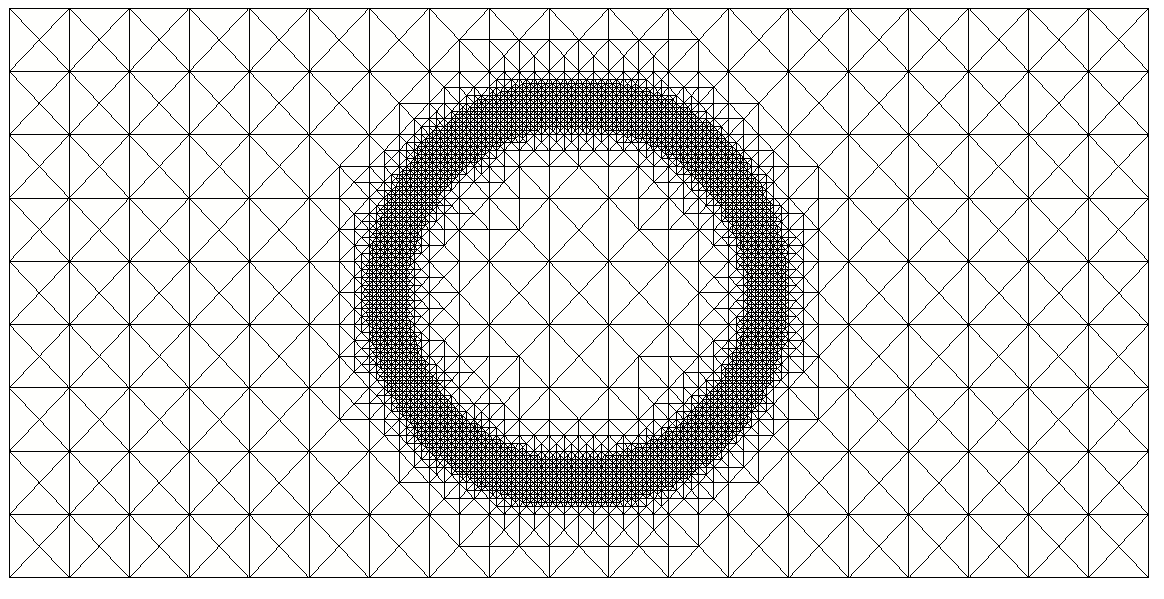}\\[0.2cm]
 \includegraphics[scale=0.15]{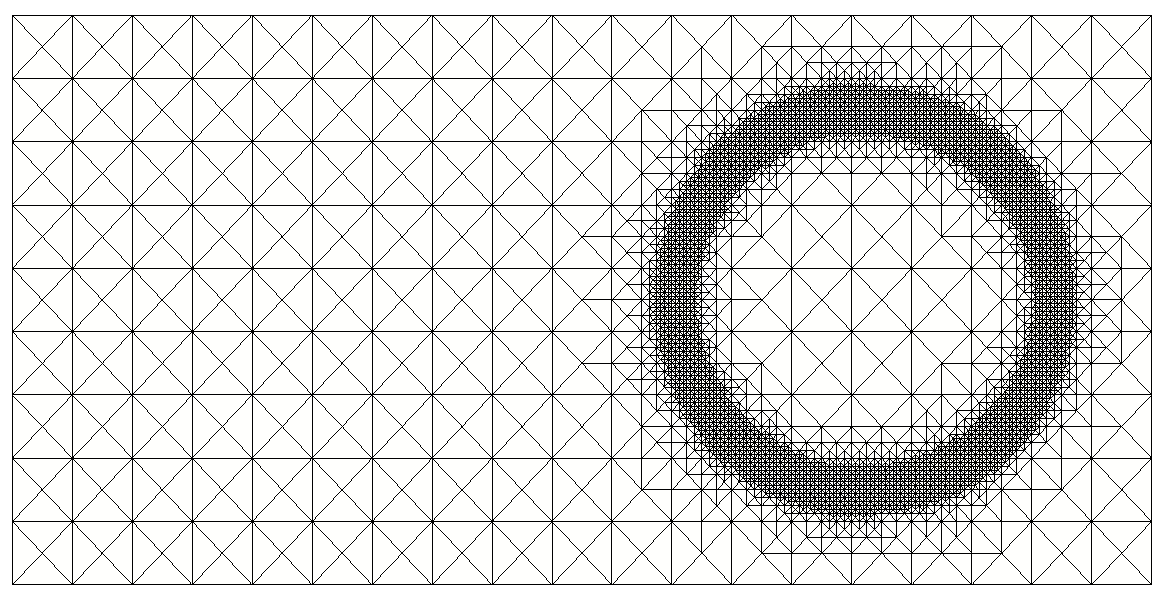} \hspace{0.2cm} \includegraphics[scale=0.15]{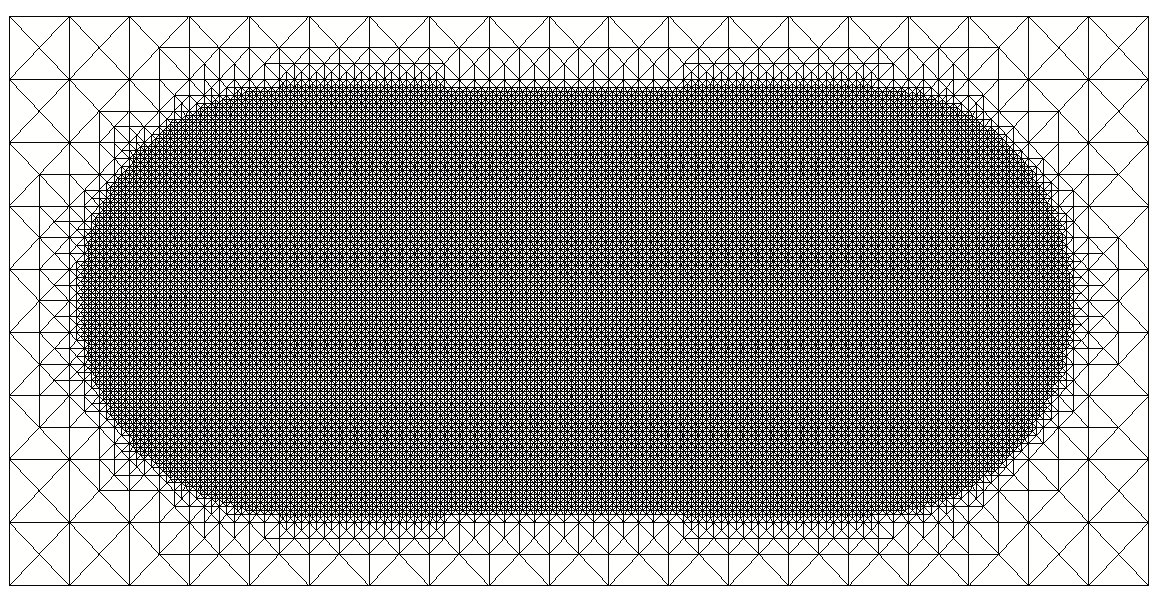}
 \caption{\em Example 6.2: Adaptive finite element meshes and finest mesh}
 \label{fig:CH_grids}
 \end{figure}

  \noindent The adaptive finite element meshes as well as the finest mesh which is generated during the adaptive finite element simulation are shown in Figure \ref{fig:CH_grids}. The number of degrees of freedom in the adaptive meshes varies between 6113 and 8795. The finest mesh (overlay of all adaptive meshes) has 54108 degrees of freedom, whereas a uniform mesh with discretization fineness as small as the smallest triangle in the adaptive meshes has 88450 degrees of freedom.

 \noindent In this example, we only compare the solution to the POD-ROM utilizing two kinds of snapshot discretizations: on the one hand we use adaptive finite elements and express these with respect to the finite element basis functions corresponding to the finest mesh. On the other hand we compute the solution to the POD-ROM with snapshots computed on a uniform finite element discretization, where the fineness is chosen to be of the same size as the smallest triangle in the adaptive meshes. We choose $X=L^2(\Omega)$ and compute a separate POD basis for each of the variables $c$ and $w$.

\begin{figure}[H]
 \centering
  \includegraphics[scale=0.4]{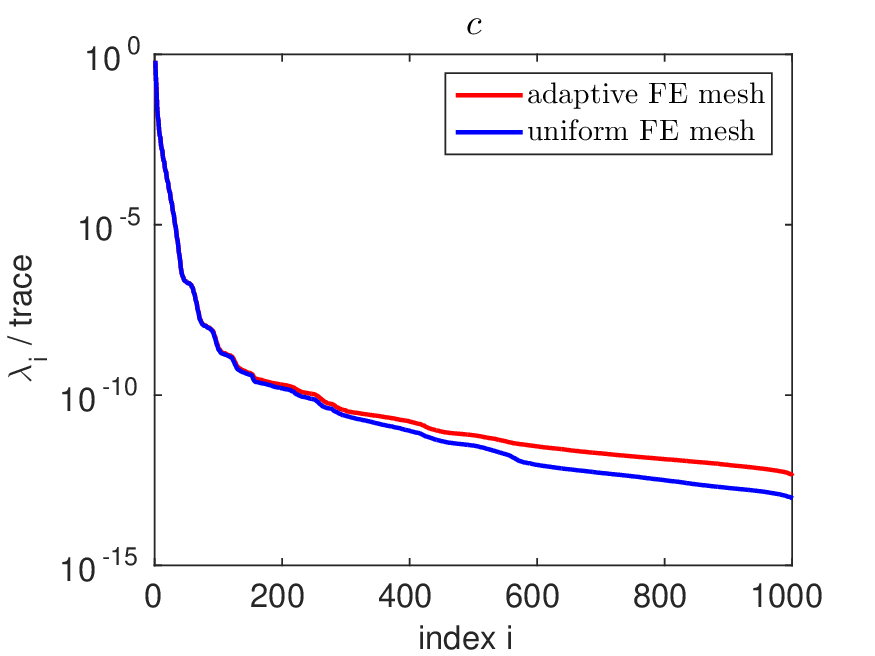} \includegraphics[scale=0.4]{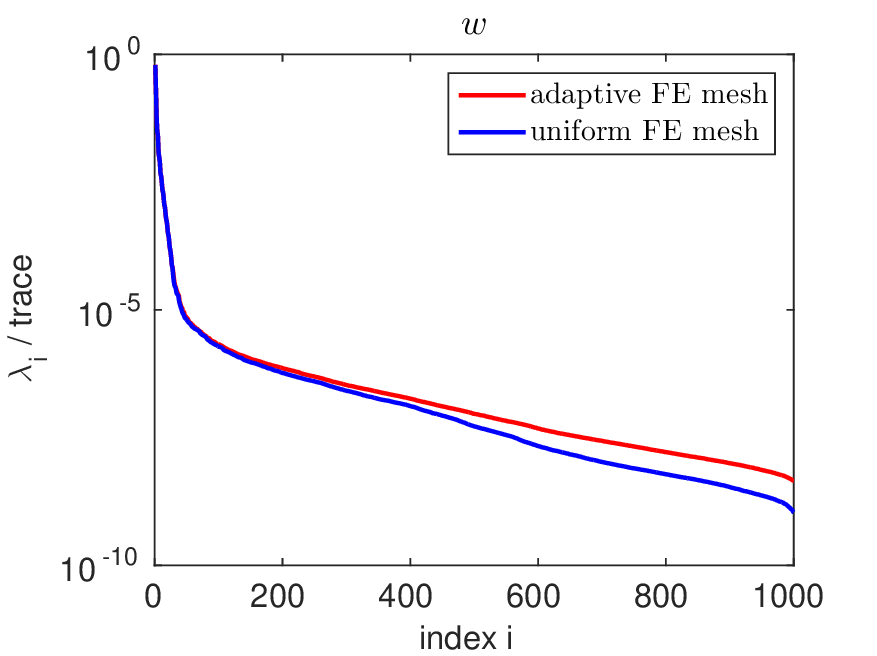}
 \caption{\em Example 6.2: Comparison of the normalized eigenvalues for the phase field $c$ (left) and the chemical potential $w$ (right) utilizing an adaptive and a uniform spatial mesh, respectively}
 \label{fig:CH_ev}
 \end{figure}

 In Figure \ref{fig:CH_ev}, a comparison is visualized concerning the normalized 
 eigenspectrum for the phase field $c$ and the chemical potential $w$ utilizing 
 uniform and adaptive finite element discretization. We note for the phase field 
 $c$ that about the first 180 eigenvalues computed corresponding to the adaptive simulation coincide with the eigenvalues of the simulation on the finest mesh. Then, the eigenvalues corresponding to the uniform simulation decay faster. Similar observations apply for the chemical potential $w$. 
 
 \noindent In order to control the accuracy of the POD-ROM, we 
 utilize the following criterion. The information content  of a POD basis of rank $\ell$ relatively to the amount of the information content of all snapshots is given by the ratio of modeled information and total information. It is defined by
 \begin{equation}\label{infocont}
  \Gamma(\ell) := \displaystyle\frac{\sum_{i=1}^\ell \lambda_i}{\sum_{i=1}^d \lambda_i}. 
  \end{equation}
 We will choose the POD basis length $\ell_c$ for the phase field $c$ and the number of POD modes $\ell_w$ for the chemical potential, such that
 $$ \ell_{\min} = \text{argmin} \{ \Gamma(\ell): \Gamma(\ell) > 1 - p \}, \quad  \text{with } \ell = \ell_c \text{ and } \ell_w, \text{ respectively,}$$
 for a given value $p $ representing the loss of information. Alternatively, one can choose the POD basis length such that the POD projection error \eqref{projection_error} is smaller or equal to $\min ((\Delta t)^2, \varepsilon_h^2)$, compare \eqref{errorestimation}. Table \ref{tab:loss_info} summarizes how to choose $\ell_c$ and $\ell_w$ in order to capture a desired amount of information. Moreover, it tabulates the POD projection error \eqref{projection_error} depending on the POD basis length, where $\lambda_i^c$ and $\lambda_i^w$ denote the eigenvalues for the phase field $c$ and the chemical potential $w$, respectively. The results in Table \ref{tab:loss_info} agree with our expectations: the smaller the loss of information $p$ is, the more POD modes are needed and the smaller is the POD projection error.

\begin{table}[H]
\centering
 \begin{tabular}{  c || c | c || c | c || c | c || c | c }
  $ p  $ & $\ell_c^{\text{ad}} $ &  $\sum_{i=\ell+1}^d \lambda_i^c$ & $\ell_w^{\text{ad}} $ &  $\sum_{i=\ell+1}^d \lambda_i^w$ &  $\ell_c^{\text{uni}} $ &  $\sum_{i=\ell+1}^d \lambda_i^c$ & $\ell_w^{\text{uni}} $ & $\sum_{i=\ell+1}^d \lambda_i^w$\\
   \hline
  $10^{-01}$ & 3 &  $2.0 \cdot 10^{-3}$ & 4 &  $156.9 \cdot 10^{0}$ & 3  &  $2.0 \cdot 10^{-3}$ & 4 &  $157.6 \cdot 10^{0}$\\
  $10^{-02}$ & 10 &  $2.1 \cdot 10^{-4}$ & 13 & \hspace{0.03cm} $15.8 \cdot 10^{0}$ & 10 &  $2.1 \cdot 10^{-4}$ & 13 &  \hspace{0.05cm} $15.6 \cdot 10^{0}$\\
  $10^{-03}$ & 19 &  $2.5 \cdot 10^{-5}$ & 26  & \hspace{0.18cm} $1.8 \cdot 10^{0}$ & 19 &  $2.5 \cdot 10^{-5}$ & 25 & \hspace{0.22cm} $1.8 \cdot 10^{0}$\\
  $10^{-04}$ & 29 &  $2.0 \cdot 10^{-6}$ & 211 & \hspace{0.4cm}  $1.8 \cdot 10^{-1}$ & 28 &  $2.6 \cdot 10^{-6}$ & 160 & \hspace{0.45cm} $1.9 \cdot 10^{-1}$\\
  $10^{-05}$ & 37 & $2.5 \cdot 10^{-7}$ & 644 & \hspace{0.4cm}  $1.1 \cdot 10^{-2}$ & 37 &  $2.4 \cdot 10^{-7}$ & 419  & \hspace{0.47cm} $2.5 \cdot 10^{-2}$\\
 \end{tabular} \hspace{0.5cm}
  \vspace{0.4cm} \caption{\em Example 6.2: Number of needed POD bases in order to achieve a loss of information below the tolerance $p$ utilizing adaptive finite element meshes (columns 2-5) and uniform finite element discretization (columns 6-9) and POD projection error}
 \label{tab:loss_info}
  \end{table}

\begin{figure}[H]
 \centering
 \includegraphics[scale=0.15]{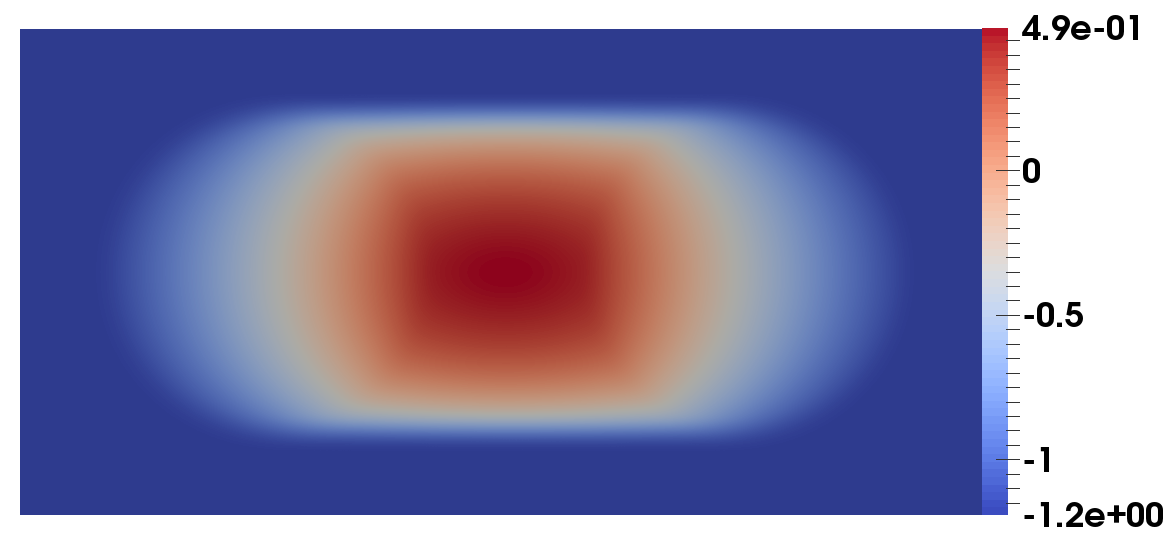} \includegraphics[scale=0.15]{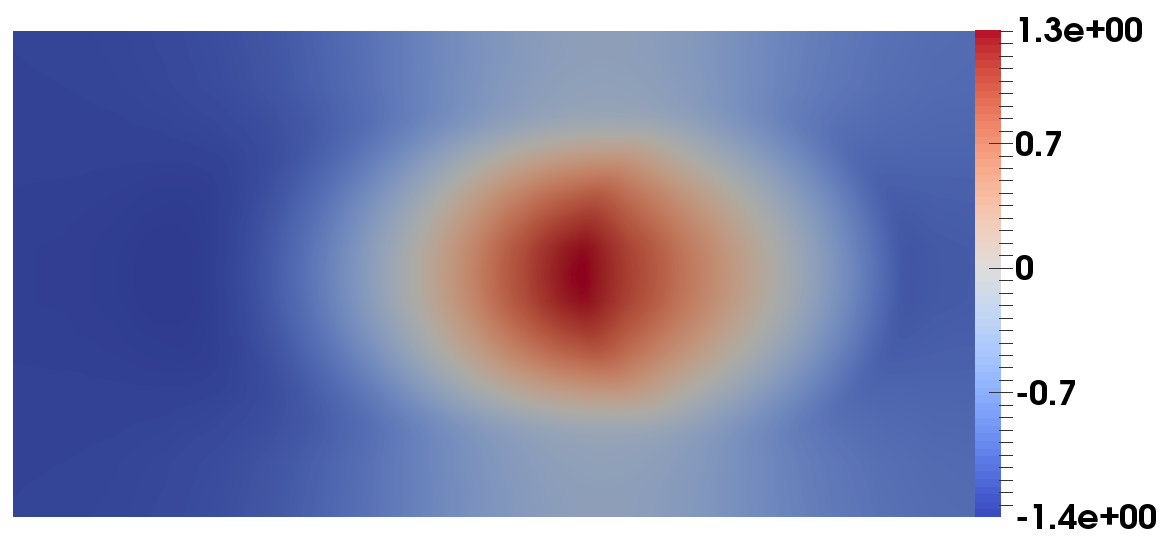}\\
  \includegraphics[scale=0.15]{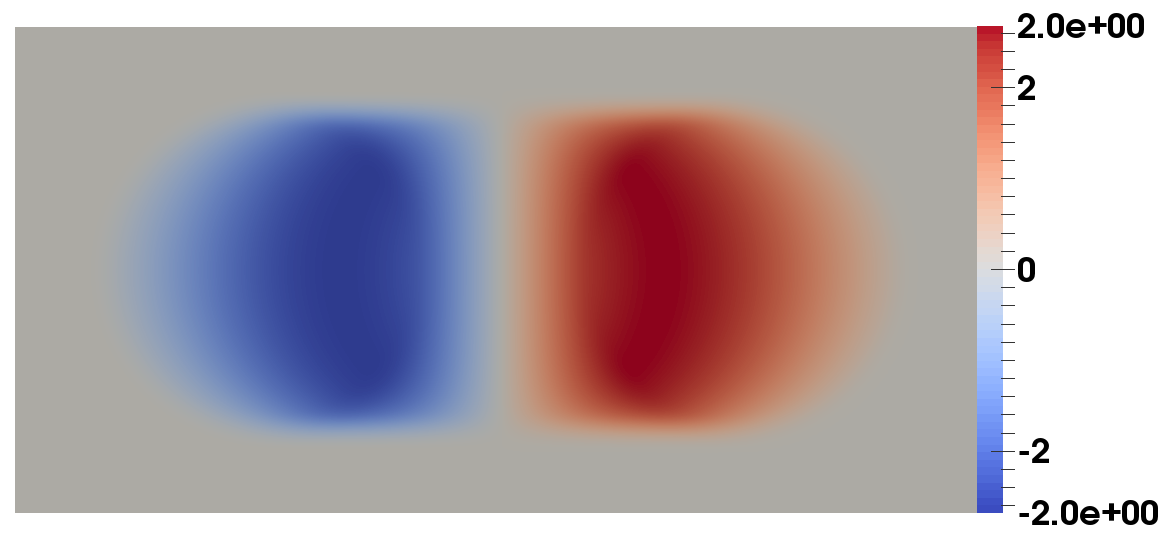} \includegraphics[scale=0.15]{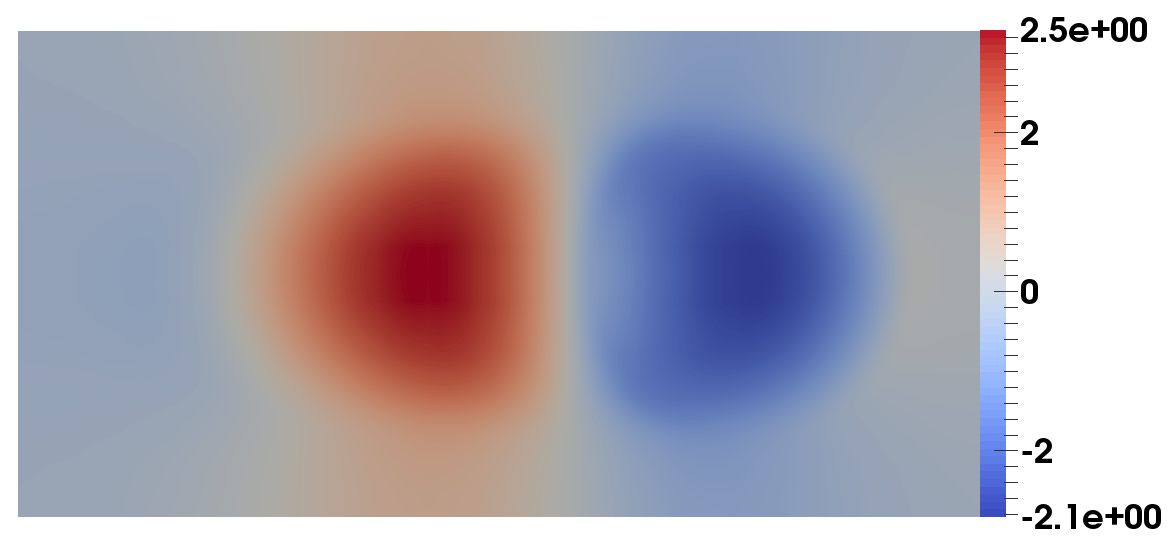}\\
 \includegraphics[scale=0.15]{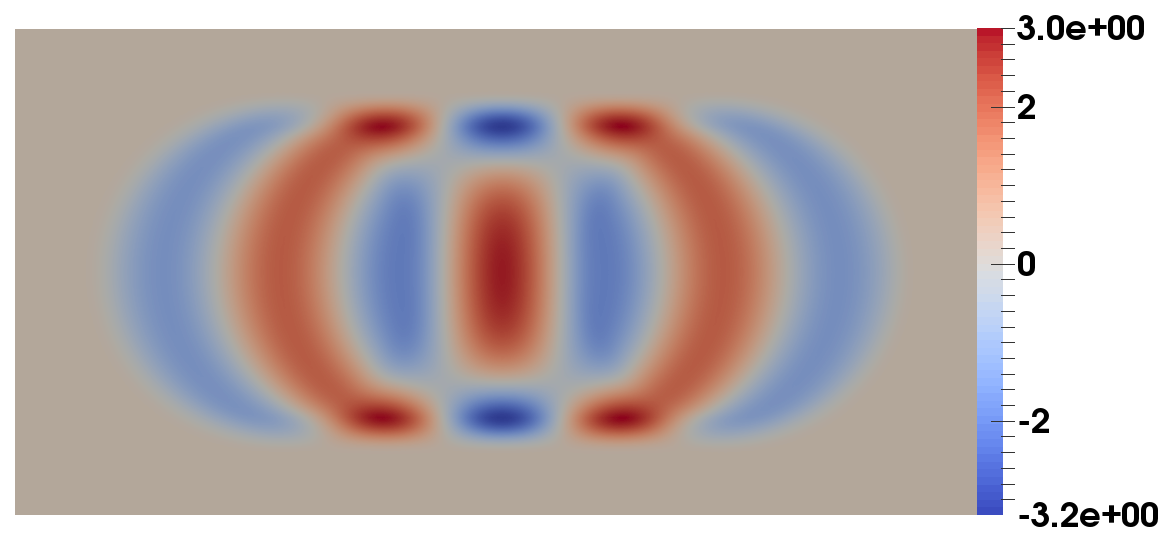} \includegraphics[scale=0.15]{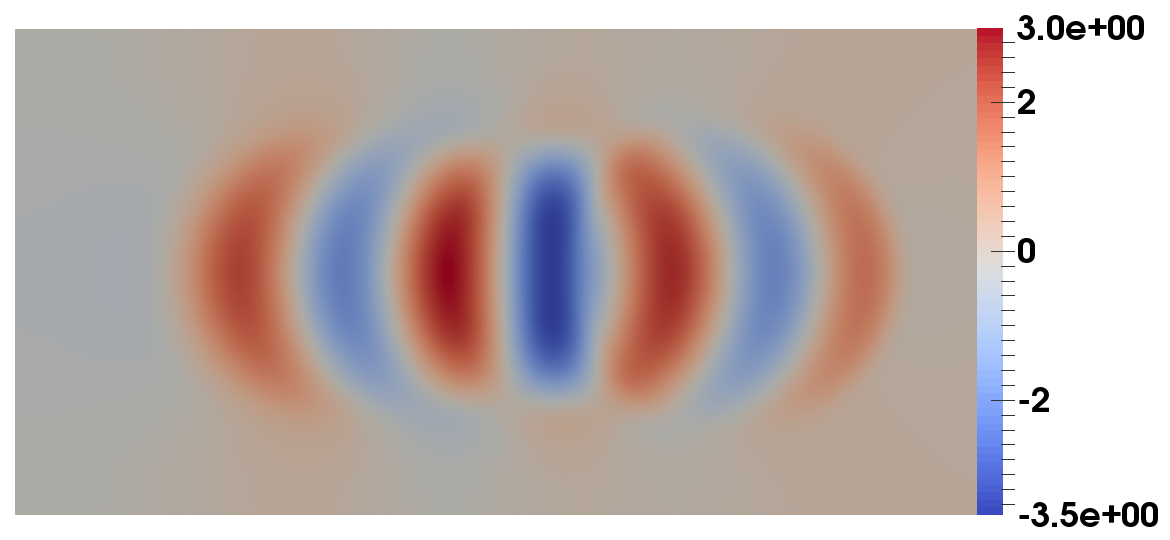}
 \caption{\em Example 6.2: First, second and fifth POD bases for $c$ (left) and $w$ (right).}
 \label{fig:CH_POD_bases}
 \end{figure}

  \noindent In the following, we run the numerical simulations for different combinations of numbers for $\ell_c$ and $\ell_w$ of Table \ref{tab:loss_info}. The approximation quality of the POD solution utilizing adaptive meshes is compared to the use of a uniform mesh in Table \ref{tab:CH_adaptiveFE}. As expected, Table \ref{tab:CH_adaptiveFE} shows that the error between the POD surrogate solution and the high-fidelity solution gets smaller for an increasing number of utilized POD basis functions. Moreover, a larger number of POD modes is needed for the chemical potential $w$ than for the phase field $c$ in order to get an error in the same order which is in accordance to the fact that the decay of the eigenvalues for $w$ is slower than for $c$ as seen in Figure \ref{fig:CH_ev}.\\
  Figure \ref{fig:CH_POD_bases} visualizes the first, second and fifth POD modes for the phase field $c$ and the chemical potential $w$. Analogue to Example 6.1, we observe a periodicity in the POD bases corresponding to its index number.

\begin{table}[H]
\centering
 \begin{tabular}{  c | c || c |  c  | c | c }
  $\ell^c$ & $\ell^w$ & $c: \varepsilon_{\text{FE}}^{\text{ad}}$ &    $w: \varepsilon_{\text{FE}}^{\text{ad}}$ & $c: \varepsilon_{\text{FE}}^{\text{uni}}$ &   $w: \varepsilon_{\text{FE}}^{\text{uni}}$  \\
 \hline
 3 & 4 & $8.44 \cdot 10^{-3}$ & \hspace{-0.3cm} $3.00 \cdot 10^{0}$ &  $8.44 \cdot 10^{-3}$  & \hspace{-0.3cm} $3.75 \cdot 10^{0}$\\
 10 & 13 & $3.30 \cdot 10^{-3}$  & $3.77 \cdot 10^{-1}$ & $ 3.30 \cdot 10^{-3} $ & $4.32 \cdot 10^{-1}$ \\
 19 & 26 & $1.57 \cdot 10^{-3}$  & $2.12 \cdot 10^{-1}$ & $1.57 \cdot 10^{-3}$ &  $2.39 \cdot 10^{-1}$ \\
 29 & 26 & $7.34 \cdot 10^{-4}$  & $1.09 \cdot 10^{-1}$ & $7.32 \cdot 10^{-4}$ &  $1.16 \cdot 10^{-1}$\\
 37 & 26 & $3.57 \cdot 10^{-4}$  & $4.82 \cdot 10^{-2}$ & $ 3.55 \cdot 10^{-4}$ &  $5.04 \cdot 10^{-2}$  \\
 50 & 50 & $1.88 \cdot 10^{-4}$  & $2.17 \cdot 10^{-2}$ & $1.86 \cdot 10^{-4}$ &  $2.33 \cdot 10^{-2} $ \\
 65 & 26 & $9.74 \cdot 10^{-5}$  & $1.11 \cdot 10^{-2}$ & $9.56 \cdot 10^{-5}$ &  $1.15 \cdot 10^{-2}$\\
 100 & 100 & $3.37 \cdot 10^{-5}$ & $3.56 \cdot 10^{-3}$ & $3.22 \cdot 10^{-5}$  & $3.42 \cdot 10^{-3}$ \\
 \end{tabular}
 \vspace{0.4cm} \caption{\em Example 6.2: Relative $L^2(0,T;L^2(\Omega))$-error between the POD solution and the finite element solution utilizing adaptive meshes (columns 3-4) and utilizing a uniform mesh (columns 5-6), respectively.}
 \label{tab:CH_adaptiveFE}
  \end{table}

\begin{figure}[H]
 \centering
  \includegraphics[scale=0.4]{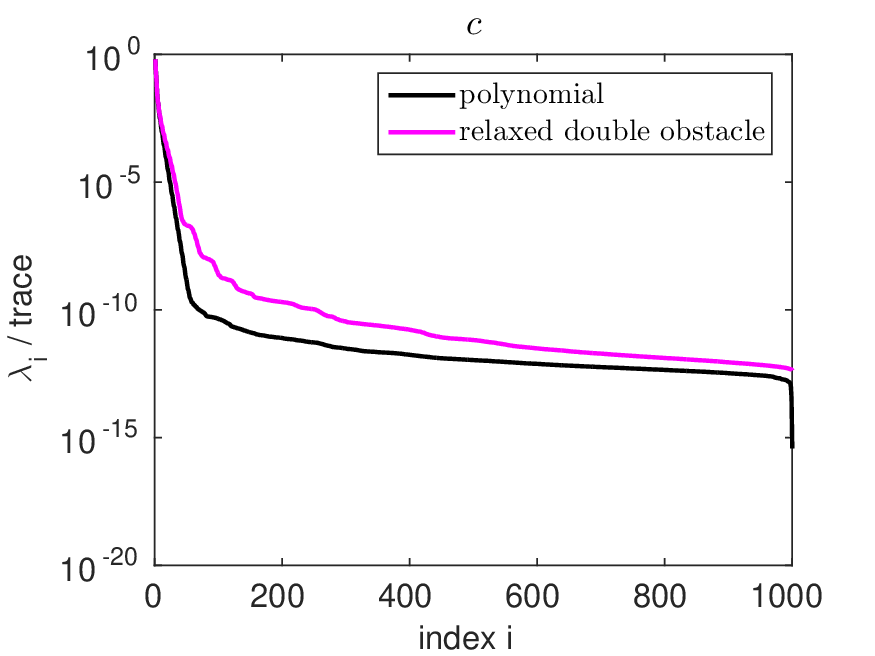} \includegraphics[scale=0.4]{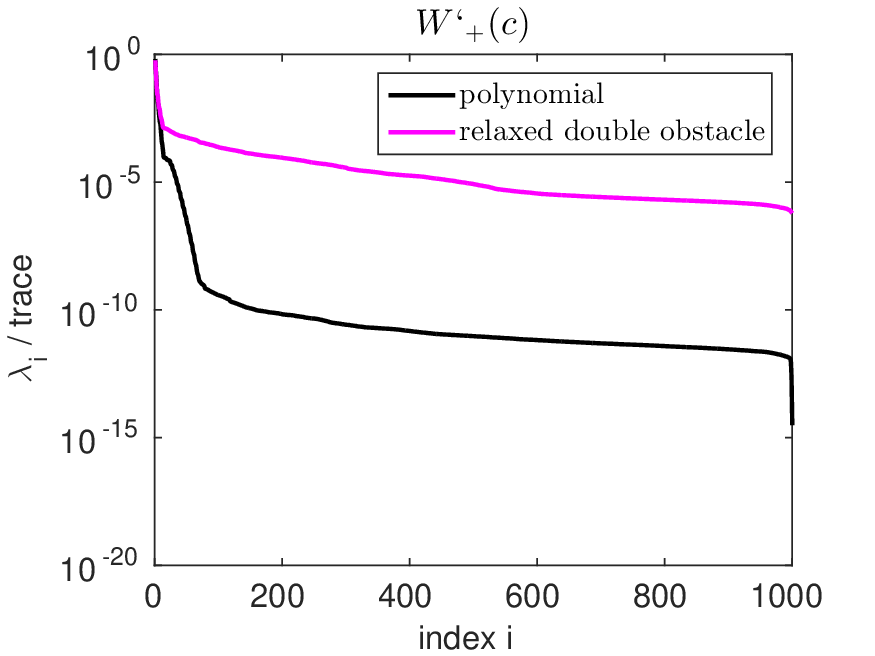} 
 \caption{\em Example 6.2: Comparison of the normalized eigenvalues for $c$ (left) and the first derivative of the convex part $W'_+$ of the free energy (right) utilizing polynomial and relaxed double obstacle energy, respectively}
 \label{fig:ev_psiprime}
 \end{figure}

    Finally, the treatment of the nonlinearity shall be discussed. Utilizing the 
  convex-concave splitting for $W$, we get for the Moreau-Yosida relaxed double obstacle free energy $W_-^{\text{rel}} (c) = \frac{1}{2}(1-c^2)$ for the concave part and $W_+^{\text{rel}}(c) = \frac{s}{2}(\max (c-1,0)^2 + \min (c+1,0)^2)$ for the convex part. This means that the first derivative of the concave part is linear with respect to the phase field variable $c$. The challenging part is the convex term which first derivative is non-smooth. For a comparison, we consider the smooth polynomial free energy which concave part is $W_-^p(c) = \frac{1}{4}(1-2c^2)$ and convex part is $W_+^p(c) = \frac{1}{4}c^4$.\\
  Figure \ref{fig:ev_psiprime} shows the decay of the normalized eigenspectrum for the phase field $c$ (left) and the first derivative of the convex part $W'_+(c)$ (right) for the polynomial and the relaxed double obstacle free energy. Obviously, in the non-smooth case more POD modes are needed for a good approximation than in the smooth case. This behaviour is similar to the decay of the Fourier coefficients in the context of trigonometric approximation, where the decay of the Fourier coefficients depends on the smoothness of the approximated object.\\  
   Table \ref{tab:CH_speedups} summarizes computational times for different finite element runs as well as reduced order simulations utilizing the polynomial and the relaxed double obstacle free energy, respectively. In addition, the approximation quality is compared. The computational times are rounded averages from various test runs. It turns out that the finite element simulation (row 1) using the smooth potential is around two times faster than using the non-smooth potential. This is due to

\begin{table}[H]
\centering
 \begin{tabular}{  l | r | r ||  r  | r  }
  & \multicolumn{2}{c||}{$W^p$} & \multicolumn{2}{c}{$W_s^{\text{rel}}$} \\
 \hline
  FE & \multicolumn{2}{c||}{ 1644 sec}  & \multicolumn{2}{c}{ 3129 sec} \\[0.2cm]
 \hline
  & $\ell_c=3, \ell_w=4$ & $\ell_c=19, \ell_w=26$ &  $\ell_c=3, \ell_w=4$ & $\ell_c=19, \ell_w=26$\\
\hline
 POD offline &  355 sec & 355 sec & 350 sec & 349 sec\\
 DEIM offline & 8 sec & 8 sec & 9 sec & 10 sec \\
 ROM & 183 sec & 191 sec & 2616 sec & 3388 sec\\
 ROM-DEIM & 0.05 sec & 0.1 sec & 0.04 sec & no conv. \\
 ROM-proj & 0.008 sec & 0.03 sec & 0.01 sec & 0.03 sec\\
 \hline
 speedup FE-ROM & 8.9 & 8.6 & 1.1 &  none \\
  speedup FE-ROM-DEIM &  32880 &  16440  &  78225 & --\\
 speedup FE-ROM-proj &  205500 &  54800 &  312900 & 104300 \\
 \hline
  rel $L^2(Q)$ error ROM & $5.46 \cdot 10^{-03}$ & $3.23 \cdot 10^{-04}$ &  $8.44 \cdot 10^{-03}$ &  $1.57 \cdot 10^{-03}$  \\
    rel $L^2(Q)$ error ROM-DEIM & $1.46 \cdot 10^{-02}$ & $ 3.83 \cdot 10^{-04}$  & $8.84 \cdot 10^{-03}$  & --\\
   rel $L^2(Q)$ error ROM-proj & $4.70 \cdot 10^{-02}$ & $ 4.18 \cdot 10^{-02}$  &  $8.72 \cdot 10^{-03}$  &  $ 9.80 \cdot 10^{-03}$ \\
 \end{tabular}
 \vspace{0.4cm} \caption{\em Example 6.2: Computational times, speedup factors and approximation quality for different POD basis lengths and using different free energy potentials}
 \label{tab:CH_speedups}
  \end{table}
  
 \noindent the fact that in the smooth case, 2-3 Newton steps are needed for convergence in each time step, whereas in the non-smooth case 6-8 iterations are needed in the semismooth Newton method.\\
  Utilizing the smooth polynomial free energy, the reduced order simulation is 8-9 times faster than the finite element simulation, whereas utilizing the relaxed double obstacle free energy only delivers a very small speedup. The inclusion of DEIM (we use $\ell_{\text{deim}}=\ell_c$) in the reduced order model leads to immense speedup factors for both free energy functions (row 8). This is due to the fact that the evaluation of the nonlinearity in the reduced order model is still dependent on the full spatial dimension and hyper reduction methods are necessary for useful speedup factors. Note that the speedup factors are of particular interest in the context of optimal control problems. At the same time, the relative $L^2(0,T;L^2(\Omega))$-error between the finite element solution
  and the ROM-DEIM solution is close to the quality of the reduced order model solution (row 10-11).\\
  However, in the case of the non-smooth free energy function utilizing $\ell_c=19$ POD modes for the phase field and $\ell_w=26$ POD modes for the chemical potential, the inclusion of DEIM has the effect that the semismooth Newton method does not converge. For this reason, we treat the nonlinearity by applying the technique explained in Section 4.2., i.e. we project the finite element snapshots for $W'_+(c)$ (which are interpolated onto the finest mesh) onto the POD space. Since this leads to linear systems, the computational times are very small (row 6). The error between the finite element solution and the reduced order solution utilizing projection of the nonlinearity lies in the area $10^{-02}/10^{-03}$. 
  Depending on the motivation, this approximation quality might be sufficient. Nevertheless, we note that that for large numbers of POD modes, utilizing the projection of the nonlinearity onto the POD space leads to a large increase of the error.\\

  \noindent \textbf{Example 6.3: Linear heat equation (revisited).}
  
\begin{figure}[H]
 \centering
  \includegraphics[scale=0.1]{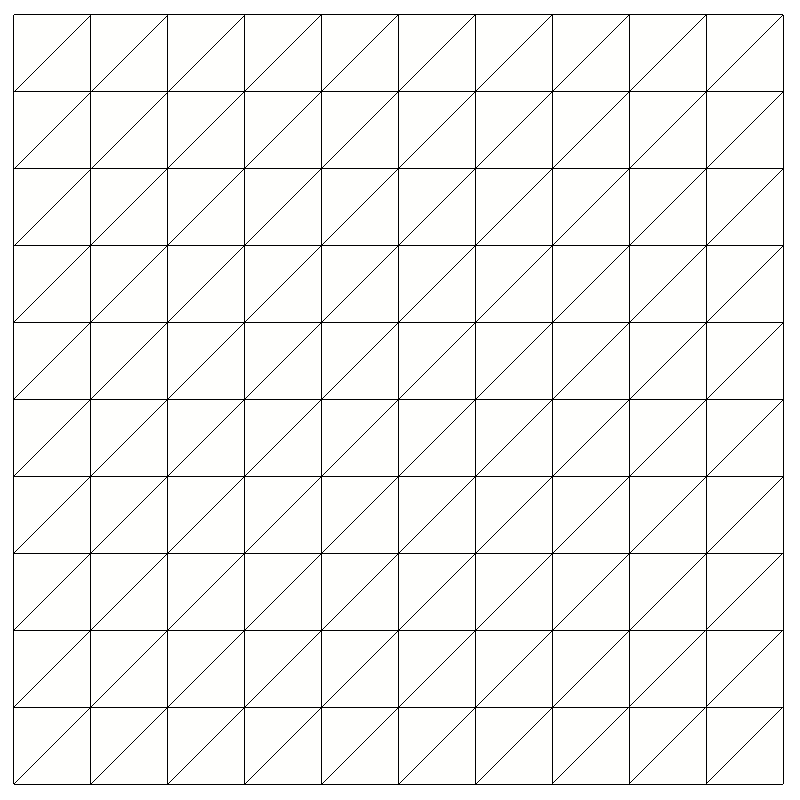} 
    \includegraphics[scale=0.1]{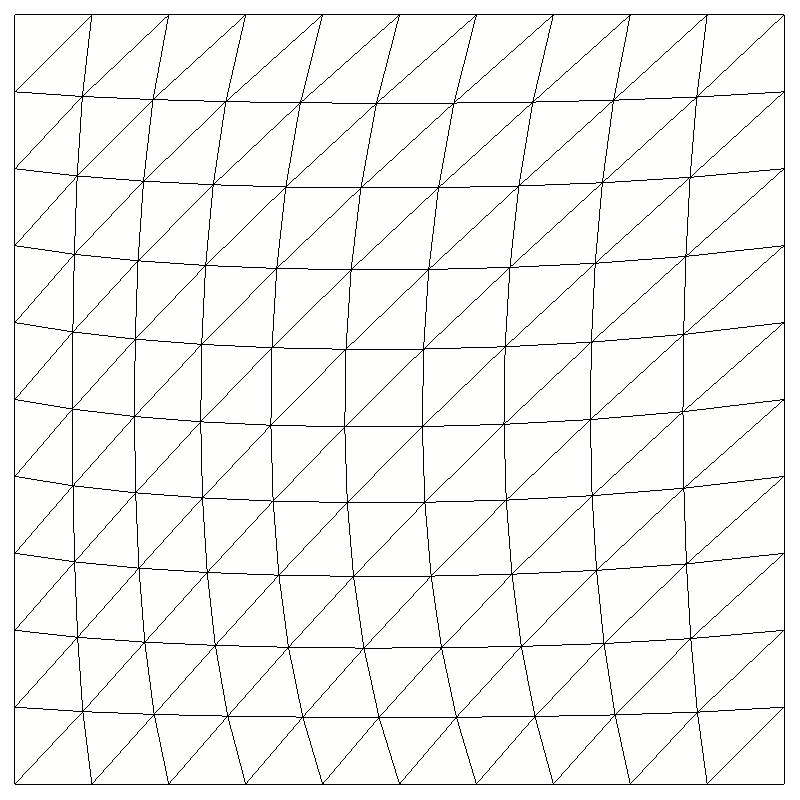} 
      \includegraphics[scale=0.1]{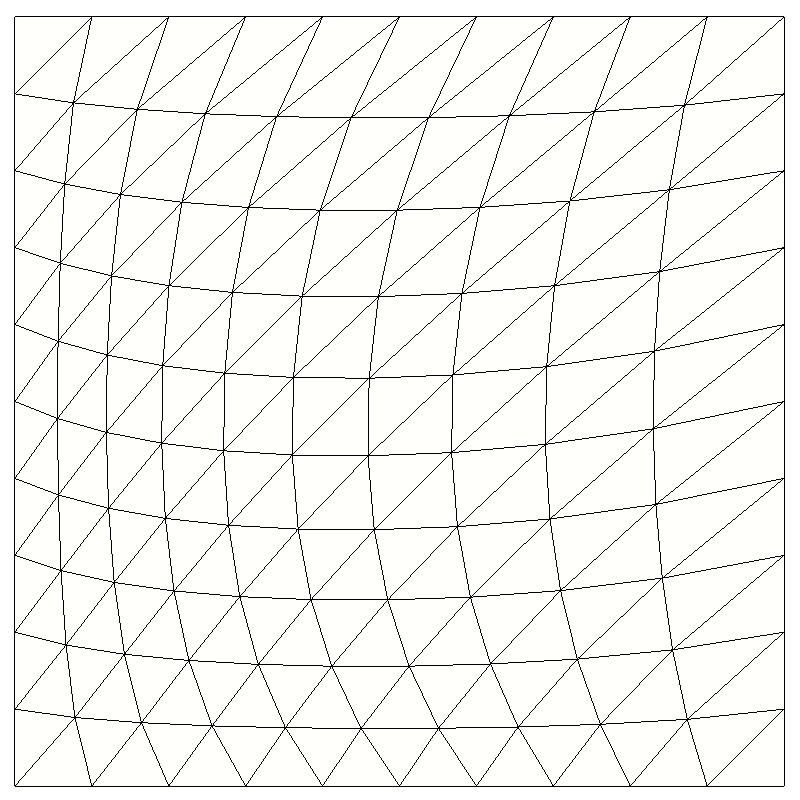} 
        \includegraphics[scale=0.1]{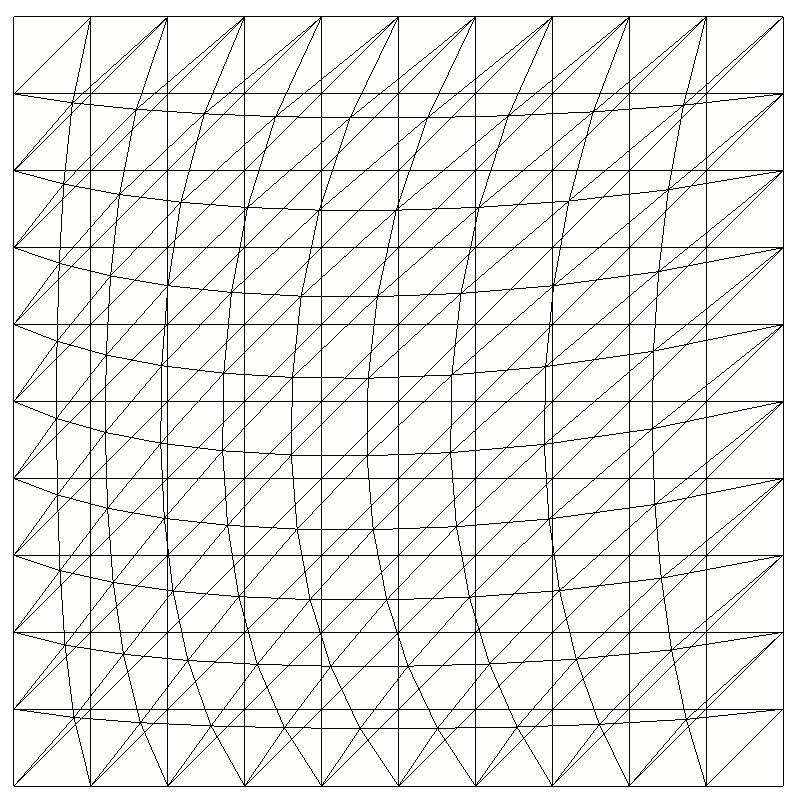} 
 \caption{\em Example 6.3: Uniform mesh (left), disturbed meshes at $t=0.5$ and $t=1.0$ (middle left, middle right), overlap of the mesh at $t=0$ with 
 the mesh at $t=1.0$ (right). Here, we use $\theta=10$.}
 \label{fig:ex3_meshes}
 \end{figure}

 \noindent Like in Example 6.1, let us consider again Example 2.3 \eqref{heat} 
of a heat equation with $c\equiv 0$ . The purpose of this example is to  
 confirm the numerical applicability of the strategy described in Section 3.2. We set up the matrix $\mathcal{K}$ for snapshots given on non-nested spatial discretization which requires the integration over cut elements. The data is chosen as follows: we consider homogeneous Dirichlet boundary conditions. As spatial domain we choose $\Omega = [0,1]\times [0,1] \subset \mathbb{R}^2$, the time interval is $[0,T]=[0,1]$, and we utilize a uniform temporal discretization with time step size $\Delta t = 0.01$. We construct an example such that the analytical solution is known. It is given by
  $$ y(t,x) = \sin(\pi x_0)\cdot \sin(\pi x_1) \cdot \cos(2 \pi t x_0).$$
  The source term $f$ and the initial condition $g$ are chosen accordingly. The initial condition is discretized using piecewise linear and continuous finite 
  elements on a uniform spatial mesh which is shown in Figure \ref{fig:ex3_meshes} (left). Then, at each time step, the mesh is disturbed by relocating each mesh 
  node according to the assignment
    \begin{equation*}
  \begin{array}{r c l}
   x_0 & \leftarrow & x_0 + \theta \cdot x_0 \cdot (x_0-1) \cdot (\Delta t / 10),\\
    x_1 & \leftarrow & x_1 + \theta \cdot 0.5 \cdot x_1 \cdot (x_1-1) \cdot (\Delta t / 10),
   \end{array}
  \end{equation*}
  where $\theta \in \mathbb{R}_+$ is sufficiently small such that all coordinates of the inner node points fulfill $0<x_0<1$ and $0<x_1<1$. After relocating the mesh nodes, the heat equation is solved on this mesh for the next time instance. For this, we use the Lagrange interpolation in order to interpolate the finite element solution of the previous time step onto the new mesh. The disturbed meshes at $t=0.5$ and $t=1.0$ as well as an overlap of two meshes are shown in Figure \ref{fig:ex3_meshes}. We follow the strategy explained 
  in Section 3.2 and compute the matrix $\mathcal{K}$ from \eqref{mathcalK} by evaluating the inner products of the snapshots, where we need to integrate over cut elements.

\begin{figure}[htbp]
 \centering
  \includegraphics[scale=0.4]{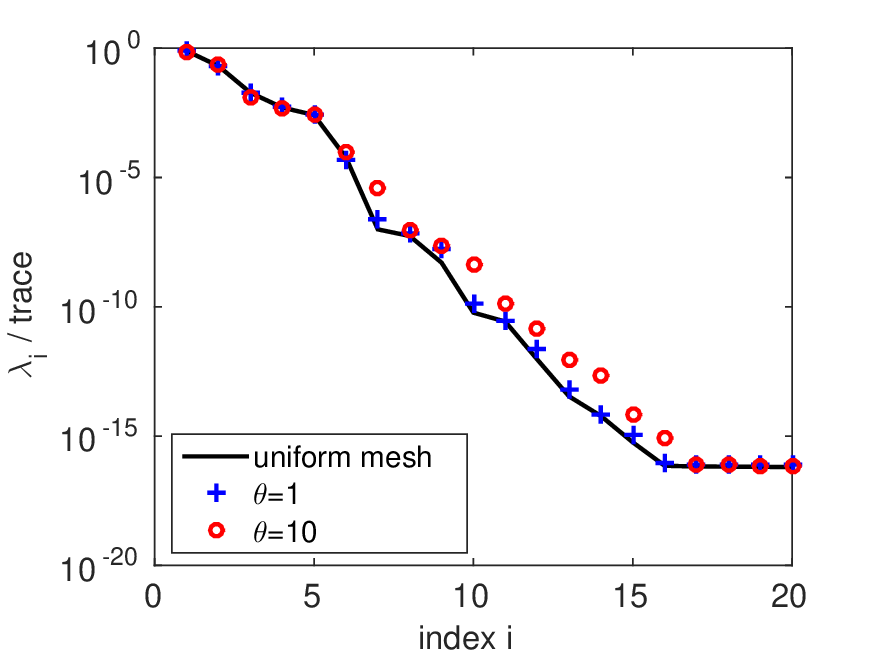}   \includegraphics[scale=0.4]{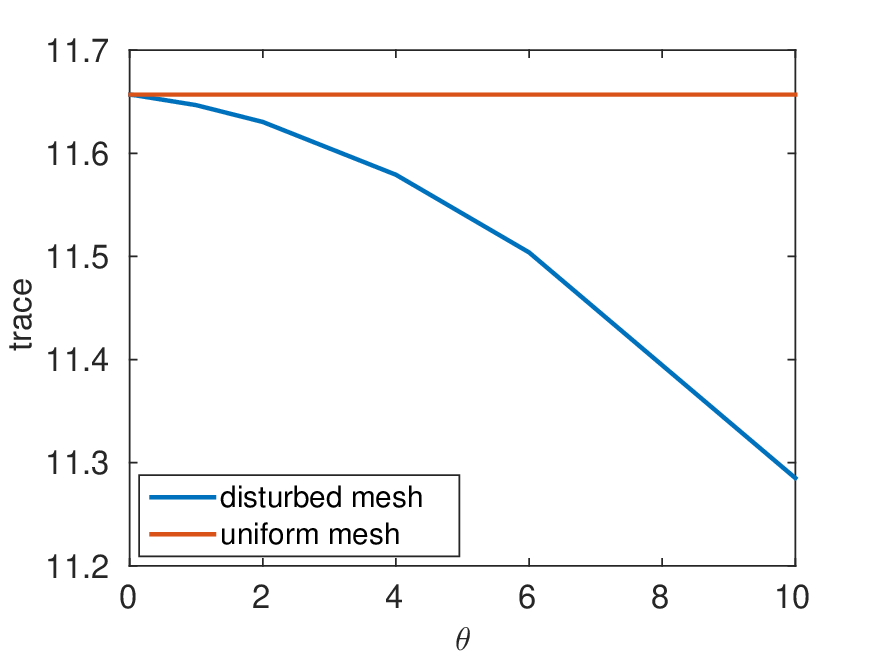} 
 \caption{\em Example 6.3: Decay of eigenvalues of matrix $\mathcal{K}$ with different meshes}
 \label{fig:ex3_ev}
 \end{figure}

 \noindent We compute the eigenvalue decomposition of $\mathcal{K}$ for different 
 values of $\theta$ and compare the results with a uniform mesh (i.e. $\theta = 0$) in Figure \ref{fig:ex3_ev}. We note that the eigenvalues of the disturbed mesh are converging to the eigenvalues of the uniform mesh for $\theta \to 0$. As expected, the eigenvalue spectrum depends only weakly on the underlying mesh given that the mesh size is sufficiently small. The POD-ROM follows along the lines of Examples 6.1 and 6.2.

  Our last remark concerns the computational complexity. Solving the heat equation takes 2.1sec on the disturbed meshes and 1.8sec on the uniform mesh. The computational time for each entry of the matrix $\mathcal{K}$ is 0.022sec and computing the eigenvalue decomposition for $\mathcal{K}$ takes 0.0056sec. Note that the cut element integration problem for each matrix entry takes a fraction of time required to solve the finite element problem.

 \section{Conclusion}

 \noindent In this work, a POD reduced order model is proposed which can be 
 set up and solved for snapshots which are discretized utilizing arbitrary finite elements. The method is applicable for $h$-, $p$- and $r$-adaptive finite elements. The approach is motivated from an infinite-dimensional perspective. Using the method of snapshots we are able to set up the correlation matrix by evaluating the inner products of snapshots which live in different finite element spaces. For non-nested meshes, this requires the detection of cell collision and integration over cut finite elements. A numerical strategy how to implement this practically is elaborated and numerically tested. Utilizing the eigenvalues and eigenvectors of this correlation matrix, we are able to set up and solve a POD surrogate model that does not need the expression of the snapshots with respect to the basis of a common finite element space or the interpolation onto a common reference mesh. Moreover, an error bound for the error between the true solution and the solution to the POD-ROM using spatially adapted snapshots is derived. The error estimation contains an additional term according to the spatial discretization error compared to existing error bounds. The numerical tests show that the POD projection error decreases if the number of utilized POD basis functions is increased. However, the error between the POD solution and the true solution stagnates when the spatial discretization error dominates. Moreover, the numerics show that utilizing the correlation matrix calculated explicitly without interpolation in order to build a POD-ROM gives the same results as the approach in which the snapshots are interpolated onto the finest mesh. From a computational point of view, sufficient hardware should be available in order to compute the correlation matrix in parallel and make the offline computational time competitive. For semilinear evolution problems, the nonlinearity is treated by linearization. This is of interest in view of optimal control problems, in which a linearized state equation has to be solved in each SQP iteration level. In future work, we intend to study the combination of adaptive finite elements and POD reduced order modeling in the context of multi-phase flow and optimal control.

\end{document}